\documentclass[12pt]{article}
\usepackage{mathrsfs}
\usepackage{amsmath, amsthm,amssymb}
\usepackage{amsfonts}
\usepackage{color,xcolor}
\usepackage{graphicx}
\usepackage{geometry}
\usepackage{manfnt}
\usepackage[pagewise]{lineno}
\usepackage[colorlinks,citecolor=blue,urlcolor=blue]{hyperref}
\allowbreak
\newtheorem{theorem}{Theorem}[section]

\newtheorem{lemma}[theorem]{Lemma}

\newtheorem{definition}[theorem]{Definition}
\newtheorem{remark}[theorem]{Remark}

\def\geq{\geqslant}\def\leq{\leqslant}

\def\div{\mathord{{\rm div}}}

\begin{document}
\title{\bf Stochastic equations with low regularity growing drifts}

\author{Jinlong Wei$^a$, Junhao Hu$^b$ and Chenggui Yuan$^c$
\\{\small \it $^a$School of Statistics and Mathematics, Zhongnan University of Economics}\\
{\small \it  and Law, Wuhan 430073, China}  \\ {\small \tt  weijinlong.hust@gmail.com}
\\ {\small \it $^b$School of Mathematics and Statistics, South-Central Minzu University} \\ {\small \it Wuhan 430074, China} \\ {\small \tt junhaohu74@163.com} \\ {\small \it $^c$Department of Mathematics, Swansea University, Bay Campus} \\ {\small \it Swansea SA1 8EN, United Kingdom} \\ {\small \tt C.Yuan@swansea.ac.uk}}

\date{}

\maketitle
\noindent{\hrulefill}
\vskip1mm\noindent{\bf Abstract} By using the It\^{o}-Tanaka trick, we prove the unique strong solvability as well as the gradient estimates for stochastic differential equations with irregular growing drifts in low regularity Lebesgue-H\"{o}lder space $L^q(0,T;
{\mathcal C}^{\frac{2}{q}-1}\cap{\mathcal C}^\alpha({\mathbb R}^d))$ with $\alpha\in(0,1)$ and $q\in (2/(1+\alpha),2$). As applications, we show the unique weak and strong solvability for stochastic transport equations driven by the low regularity growing drift with $q\in (4/(2+\alpha),2$) as well as the local Lipschitz estimate for stochastic strong solutions.


\vskip1mm\noindent{\hrulefill}

\section{Introduction}\label{sec1}\setcounter{equation}{0}
Let $T>0$ be a given real number. We are concerned with the following stochastic differential equation (SDE for short) in ${\mathbb R}^d$:
\begin{eqnarray}\label{1.1}
dX_{s,t}=b(t,X_{s,t})dt+\sigma(t)dW_t, \ t\in (s,T], \
X_{s,t}|_{t=s}=x\in{\mathbb R}^d,
\end{eqnarray}
where $\{W_t\}_{0\leq t\leq T}=\{(W_{1,t}, \ldots, W_{d,t})^\top\}_{0\leq t\leq T}$ is
a $d$-dimensional standard Wiener process defined on a given stochastic
basis ($\Omega, {\mathcal F},{\mathbb P},\{{\mathcal F}_t\}_{0\leq t\leq T}$) and $b: [0,T]\times{\mathbb R}^d\rightarrow{\mathbb R}^d$, $\sigma: [0,T]\rightarrow{\mathbb R}^{d\times d}$
are Borel measurable functions.

\medskip
When $b$  and $\sigma$ are bounded, Veretennikov \cite{Ver}
first proved the unique strong solvability for SDE (\ref{1.1}).
Since then, Veretennikov's result was strengthened in different forms, see \cite{Dav07,MNP,WLW}. When $\sigma=I_{d\times d}$ and $b$ is more regular, i.e., $b\in L^\infty(0,T;{\mathcal C}_b^\alpha({\mathbb R}^d;{\mathbb R}^d))$ with $\alpha\in (0,1)$, Flandoli, Gubinelli and Priola \cite{FGP1} proved that the unique strong solution forms a ${\mathcal C}^{1,\alpha^\prime} \ (\alpha^\prime\in (0,\alpha)$) stochastic  flow of diffeomorphisms. This result was then generalized by Wei, Duan, Gao and Lv \cite{WDGL} to the case of $b\in L^q(0,T;{\mathcal C}_b^\alpha({\mathbb R}^d;{\mathbb R}^d))$ with $q>2/\alpha$.

\medskip
When $\sigma=I_{d\times d}$ and $b$ is not bounded but only integrable, which is in the Krylov-R\"{o}ckner class:
\begin{eqnarray}\label{1.2}
b\in L^q(0,T;L^p({\mathbb R}^d;{\mathbb R}^d))
\end{eqnarray}
with some $p,q\in [2,+\infty]$ such that
\begin{eqnarray}\label{1.3}
-\frac{2}{q}-\frac{d}{p}>-1,
\end{eqnarray}
the unique strong solvability for SDE (\ref{1.1}) was first obtained by Krylov and
R\"{o}ckner \cite{KR}. Recently, by using the It\^{o}-Tanaka trick, Fedrizzi
and Flandoli \cite{FF13} proved further that the unique strong solution forms a ${\mathcal C}^{\alpha^\prime}$ ($\alpha^\prime\in (0,1)$) stochastic  flow of homeomorphisms. Some further extensions for non-constant diffusion coefficients can be found in Zhang \cite{Zha05,Zha11}, Zhang and Yuan \cite{ZY}. More recently, Xia, Xie, Zhang and Zhao \cite{XX} studied the weak differentiability of the unique strong solution with respect to the starting point, and proved the Bismut-Elworthy-Li derivative formula for the strong solution.

\medskip
It is known that solutions of Navier-Stokes equations can be analyzed by
probabilistic representations based on SDEs with irregular coefficient $b$, see e.g., Rezakhanlou \cite{Rez}, Constantin and Iyer \cite{CI}, and from the viewpoint of Navier-Stokes equations $b$ can be taken in the critical case, i.e. the greater than sign in (\ref{1.3}) is replaced by the following equal sign:
\begin{eqnarray}\label{1.4}
-\frac{2}{q}-\frac{d}{p}=-1.
\end{eqnarray}
Therefore, the study of the unique  solvability for (\ref{1.1}), (\ref{1.2}) and (\ref{1.4}) is of very high importance. When $\sigma=I_{d\times d}$ and $p<+\infty$, the unique strong solvability for the critical case has been established by R\"{o}ckner and Zhao \cite{RZ2}, and  for more details we also refer to \cite{BFGM,Kry21,Krylov21,Krylov23,Nam,WLWu}. When $\sigma=I_{d\times d}$ and $p=+\infty$, the weak existence can also be derived by using a compactness argument \cite{RZ1}. However, it is still unknown whether the pathwise uniqueness is true or not under the critical case $q=2,~p=+\infty$, and this problem is a long-standing open problem.

\medskip
In the following, we  shall interpret the above critical condition by using a different philosophy through the notion of `degree'. Before giving the definition of degree, we introduce some notions.
\subsection{Lebesgue-H\"{o}lder spaces}\label{sec1.1}
Let $\gamma\in(0,1)$. We define the H\"{o}lder space ${\mathcal C}^\gamma({\mathbb R}^d)$, as the set consisting of all continuous functions $h:{\mathbb R}^d\rightarrow{\mathbb R}$ for which
\begin{eqnarray*}
[h]_\gamma=\sup_{x,y\in{\mathbb R}^d, x\neq y}\frac{|h(x)-h(y)|}{|x-y|^\gamma}<+\infty.
\end{eqnarray*}
The set ${\mathcal C}^\gamma({\mathbb R}^d)$ becomes
a Banach space with respect to the norm
\begin{eqnarray*}
\|h\|_{{\mathcal C}^\gamma({\mathbb R}^d)}=\sup_{x\in {\mathbb R}^d}\frac{|h(x)|}{1+|x|^\gamma}+[h]_\gamma=: \|(1+|\cdot|^\gamma)^{-1}h(\cdot)\|_0+[h]_\gamma.
\end{eqnarray*}
For $\gamma^\prime\in (0,\gamma)$, if we define the norm on ${\mathcal C}^{\gamma^\prime}\cap{\mathcal C}^\gamma ({\mathbb R}^d)$ by
\begin{eqnarray*}
\|h\|_{{\mathcal C}^{\gamma^\prime}\cap{\mathcal C}^\gamma ({\mathbb R}^d)}= \|(1+|\cdot|^{\gamma^\prime})^{-1}h(\cdot)\|_0+[h]_{\gamma^\prime}+[h]_\gamma,
\end{eqnarray*}
then ${\mathcal C}^{\gamma^\prime}\cap{\mathcal C}^\gamma ({\mathbb R}^d)$ is a Banach space as well. We then define ${\mathcal C}_b^\gamma({\mathbb R}^d)$ as the subset of ${\mathcal C}^\gamma({\mathbb R}^d)$, in which all elements are bounded, and for $h\in {\mathcal C}_b^\gamma({\mathbb R}^d)$, we define
\begin{eqnarray*}
\|h\|_{{\mathcal C}_b^\gamma({\mathbb R}^d)}=\sup_{x\in{\mathbb R}^d}|h(x)|+\sup_{x,y\in{\mathbb R}^d, x\neq y}\frac{|h(x)-h(y)|}{|x-y|^{\gamma}}= \|h\|_0+[h]_\gamma.
\end{eqnarray*}
Moreover, if $h\in {\mathcal C}_b^\gamma({\mathbb R}^d)$ and $\nabla^i h$ ($i$-th order gradient with respect to the space variable), $i=1,2,\ldots,k\in {\mathbb N}\backslash\{0\}$, are bounded and continuous, and $[\nabla^k h]_\gamma$ is finite, then we say $h\in {\mathcal C}_b^{k,\gamma}({\mathbb R}^d)$.  For $h\in {\mathcal C}_b^{k,\gamma}({\mathbb R}^d)$, we define
\begin{eqnarray*}
\|h\|_{{\mathcal C}_b^{k,\gamma}({\mathbb R}^d)}=\sum_{i=0}^k\sup_{x\in{\mathbb R}^d}|\nabla^ih(x)|+\sup_{x,y\in{\mathbb R}^d, x\neq y}\frac{|\nabla^kh(x)-\nabla^kh(y)|}{|x-y|^{\gamma}}
=\sum_{i=0}^k\|\nabla^ih\|_0+[\nabla^kh]_\gamma.
\end{eqnarray*}
We then define the Lebesgue-H\"{o}lder space $L^q(0,T;{\mathcal C}^{\gamma^\prime}\cap{\mathcal C}^\gamma({\mathbb R}^d))$ and $L^q(0,T;{\mathcal C}_b^{k,\gamma}({\mathbb R}^d))$ in a similar way for $q\in[1,+\infty]$, $k\in {\mathbb N}$ and $0<\gamma^\prime<\gamma<1$. For  $f\in L^q(0,T;{\mathcal C}^{\gamma^\prime}\cap{\mathcal C}^\gamma({\mathbb R}^d))$ and $g\in L^q(0,T;{\mathcal C}_b^{k,\gamma}({\mathbb R}^d))$, we set
\begin{eqnarray*}
\|f\|_{L^q(0,T;{\mathcal C}^{\gamma^\prime}\cap{\mathcal C}^\gamma({\mathbb R}^d))}=\Big[\int_0^T \|f(t)\|_{{\mathcal C}^{\gamma^\prime}\cap{\mathcal C}^\gamma({\mathbb R}^d)}^qdt\Big]^{\frac{1}{q}}
\end{eqnarray*}
and
\begin{eqnarray*}
&&\|g\|_{L^q(0,T;{\mathcal C}_b^{k,\gamma}({\mathbb R}^d))}\nonumber\\ &=& \Big[\int_0^T \|g(t)\|_{{\mathcal C}_b^{k,\gamma}({\mathbb R}^d)}^qdt\Big]^{\frac{1}{q}}
\nonumber\\&=&
\Big[\sum_{i=0}^k\int_0^T\|\nabla^{i}g(t)\|_0^qdt+\int_0^T[\nabla^kg(t)]^q_\gamma dt\Big]^{\frac{1}{q}}=:\Big[\sum_{i=0}^k\|\nabla^{i}g\|_{q,0}^q +[\nabla^kg]_{q,\gamma}^q\Big]^{\frac{1}{q}},
\end{eqnarray*}
respectively, where the integrals in the above identities are interpreted as the essential supermum when $q=+\infty$. ${\mathcal C}_b^{k,\gamma}({\mathbb R}^d)$ is regarded as ${\mathcal C}_b^\gamma({\mathbb R}^d)$ when $k=0$.

\medskip
For $h\in L^\infty({\mathbb R}^d)$, we define its Poisson integral by
\begin{eqnarray}\label{1.5}
P_\xi h(x)=\frac{\Gamma(\frac{d+1}{2})}{\pi^{\frac{d+1}{2}}}\int_{{\mathbb R}^d}\frac{\xi h(x-z)}{(\xi^2+|z|^2)^{\frac{d+1}{2}}}dz, \quad \forall \ \xi\in {\mathbb R}_+.
\end{eqnarray}
By \cite[Proposition 7, p.142]{Stein}, $h\in{\mathcal C}_b^\gamma({\mathbb R}^d)$ if and only if $h\in{\mathcal C}_b({\mathbb R}^d)$ and there exists a positive constant $A$ such that
\begin{eqnarray*}
\|\partial_\xi P_\xi h\|_0=\sup_{x\in{\mathbb R}^d}|\partial_\xi P_\xi h(x)|\leq A\xi^{-1+\gamma},  \quad \forall \ \xi\in {\mathbb R}_+.
\end{eqnarray*}
Moreover, if $h\in {\mathcal C}_b^\gamma({\mathbb R}^d),$ then
$\|h\|_0+\sup\limits_{\xi>0}[\xi^{1-\gamma}\|\partial_\xi P_\xi h\|_0]$
and $\|h\|_{{\mathcal C}^\gamma_b({\mathbb R}^d)}$ are equivalent norms. Therefore, for $g\in L^q(0,T;{\mathcal C}^{k,\gamma}_{b}({\mathbb R}^d))$ ($k\in {\mathbb N}$), we can define its equivalent norm by
\begin{eqnarray}\label{1.6}
\|g\|_{L^q(0,T;{\mathcal C}^{k,\gamma}_b({\mathbb R}^d))}
&=&\Big[\sum_{i=0}^k\|\nabla^ig\|_{q,0}^q+\int_0^T\sup_{\xi>0}\|\xi^{1-\gamma}\partial_\xi P_\xi \nabla^kg(t)\|_0^q dt\Big]^{\frac{1}{q}}\nonumber\\&=&\Big[\sum_{i=0}^k\|\nabla^ig\|_{q,0}^q+\int_0^T\sup_{(\xi,x)\in {\mathbb R}_+\times {\mathbb R}^d }|\xi^{1-\gamma}\partial_\xi P_\xi \nabla^kg(t,x)|^q dt\Big]^{\frac{1}{q}}.
\end{eqnarray}

\subsection{Definitions}\label{sec1.2a}
Let $k\in {\mathbb N}$ and $p\in [1,+\infty]$, and let $W^{k,p}({\mathbb R}^d)$ be the Sobolev space consisting of all locally integrable functions $h: {\mathbb R}^d\rightarrow {\mathbb R}$ such that for every $0\leq i\leq k$, $\nabla^i h$ exists in the weak sense and belongs to $L^p({\mathbb R}^d)$ $(:=W^{0,p}({\mathbb R}^d))$. For $h_1\in W^{k,p}({\mathbb R}^d)$ and $h_2\in {\mathcal C}^{\gamma}({\mathbb R}^d)$ ($\gamma\in (0,1))$,  by the scaling transformation, it yields that
\begin{eqnarray*}
\|\nabla^k(h_1(l\cdot))\|_{L^p({\mathbb R}^d)}=
l^{k-\frac{d}{p}}\|\nabla^kh_1\|_{L^p({\mathbb R}^d)} \quad {\rm and} \quad [h_2(l\cdot)]_\gamma=l^\gamma [h_2]_\gamma,  \quad \forall \ l>0.
\end{eqnarray*}
We then define their `degree' by ${\rm deg}(h_1)=k-d/p$ and ${\rm deg}(h_2)=\gamma$, respectively. Notice that for a second order parabolic equation we can `trade' space-regularity against time-regularity at a cost of one time derivative
for two space derivatives, we define that ${\rm deg}(f_1)=k-2/q-d/p$ and ${\rm deg}(f_2)=\gamma-2/q$ if $f_1\in L^q(0,T;W^{k,p}({\mathbb R}^d))$ and $f_2\in L^q(0,T;{\mathcal C}^{\gamma}({\mathbb R}^d))$, respectively. More general, we introduce the following definition (see \cite[Definition 1.1]{TWT}).
\begin{definition}\label{def1.1} Let $p,q\in [1,+\infty]$, $k\in {\mathbb N}$ and $\gamma\in (0,1)$. If $f\in L^q(0,T;W^{k,p}_{loc}({\mathbb R}^d))$, we define its degree by
${\rm deg}(f)=k-2/q-d/p$, and if $f\in L^q(0,T;{\mathcal C}^\gamma({\mathbb R}^d))$, we define ${\rm deg}(f)=\gamma-2/q$.
\end{definition}

From the point of `degree', the critical condition (\ref{1.4}) on the drift coefficient $b$ for SDE (\ref{1.1}) can be restated as ${\rm deg}(b)=-1$. On the other hand, by the classical Cauchy-Lipschitz theorem, if $b\in L^1(0,T;Lip({\mathbb R}^d;{\mathbb R}^d))\subset L^1(0,T;W^{1,\infty}_{loc}({\mathbb R}^d;{\mathbb R}^d))$ ($\sigma\in L^2(0,T)$ is enough), It\^{o} \cite{Ito} proved that there exists a unique strong solution to SDE (\ref{1.1}), and when $b\in L^1(0,T;Lip({\mathbb R}^d;{\mathbb R}^d))$ we also have $\deg(b)=-1$. Hence, we give a unified view for SDEs between the classical It\^{o} theory and the modern Krylov-R\"{o}ckner theory. However, from the viewpoint of classical It\^{o} theory, the drift can be taken in a low regularity Banach space for time variable (such as $L^1$) if it has `good' regularity in space variable (such as Lipschitz continuity), and thus we could establish the unique strong solvability for SDE (\ref{1.1}) if the drift is in this low regularity Banach space. This is our motivation to open this study and our main result can make a bridge that connects classical It\^{o}'s theory and the modern Krylov-R\^{o}ckner theory. Since the difficulty for the critical case, we focus our discussion on the sub-critical case (i.e. ${\rm deg}(b)>-1$).

\medskip
Recently, for the sub-critical drift which is square-integrable in time variable, and bounded and H\"{o}lder continuous in space variable, Tian, Ding and Wei \cite{TDW} proved the unique strong solvability for SDE (\ref{1.1}) for Sobolev differentiable diffusion. More recently, Galeati and Gerencs\'{e}r \cite{GG} studied SDE (\ref{1.1}) for a general fractional Brownian noise with the Hurst index $H\in {\mathbb R}_+\setminus{\mathbb N}$, in which the drift is in $L^q(0,T;{\mathcal C}_b^\alpha({\mathbb R}^d;{\mathbb R}^d))$ such that $q\in (1,2]$ and $\alpha\in (1-(q-1)/(qH),1)$, by developing some new stochastic sewing lemmas, they established the unique strong solvability as well as some other properties for solutions, such as stability, continuous differentiability of the flow and its inverse and Malliavin differentiability. In particular, these results are true for the Brownian noise ($H=1/2$) with $b\in L^q(0,T;{\mathcal C}^{\alpha}_b({\mathbb R}^d;{\mathbb R}^d))$ for $\alpha\in (0,1)$ and $q\in (2/(1+\alpha),2)$ (${\rm deg}(b)>-1$). However, from It\^{o}'s theory, the drift does not need to be bounded for the space variable, and this problem has been studied by
Flandoli, Gubinelli and Priola in \cite{FGP2} for the time-independent case.
There are relatively few works to discuss SDE (\ref{1.1}) when the drift is only $q$-th integrable in time variable. It is still unknown whether SDE (\ref{1.1}) is well-posed or not when $b\in L^q(0,T;{\mathcal C}^{\alpha}({\mathbb R}^d;{\mathbb R}^d))$ for $\alpha\in (0,1)$ and $q\in (1,2)$ such that ${\rm deg}(b)>-1$. In this paper, by applying the It\^{o}-Tanaka trick and combining the regularity estimates of solutions for Kolmogorov equations, we will establish the strong well-posedness to SDE (\ref{1.1}) for a class of low regularity growing drifts. Before giving the main result, we need another definition.

\begin{definition} (\cite[p.114]{Kun90})\label{def1.2}
A stochastic flow of homeomorphisms  on the stochastic
basis
$(\Omega, \mathcal{F},{\mathbb P}, (\mathcal{F}_t)_{0\leq t\leq T})$ associated to SDE (\ref{1.1}) is a
map $(s,t,x,\omega) \rightarrow X_{s,t}(x,\omega)$, defined for
$0\leq s \leq t \leq T, \ x\in {\mathbb R}^d, \ \omega \in \Omega$ with
values in ${\mathbb R}^d$, such that

\medskip
(i) the process $\{X_{s,\cdot}(x)\}= \{X_{s,t}(x), \ t\in [s,T]\}$ is a continuous
$\{\mathcal{F}_{s,t}\}_{s\leq t\leq T}$-adapted solution of SDE (\ref{1.1}) for  every  $s\in [0,T]$ and $x\in{\mathbb R}^d$;

\medskip
(ii) ${\mathbb P}$-a.s., $X_{s,t}(\cdot)$ is a  homeomorphism, for all $0\leq s\leq t\leq T$, and the functions $X_{s,t}(x)$ and $X_{s,t}^{-1}(x)$ are continuous in $(s,t,x)$, where $X^{-1}_{s,t}(\cdot)$ is the inverse of $X_{s,t}(\cdot)$;

\medskip
(iii) ${\mathbb P}$-a.s., $X_{s,t}(x)=X_{\tau,t}(X_{s,\tau}(x))$  for all
$0\leq s\leq \tau \leq t \leq T$, $x\in {\mathbb R}^d$  and $X_{s,s}(x)=x$.
\end{definition}

Now, let us give our main results.
\subsection{Main results}\label{sec1.2}
\begin{theorem}\label{the1.3} Let $b\in L^q(0,T;{\mathcal C}^{\frac{2}{q}-1}\cap{\mathcal C}^\alpha({\mathbb R}^d;{\mathbb R}^d))$ with $\alpha\in (0,1)$ and $q\in (2/(1+\alpha),2)$, and let $\sigma\in L^\infty(0,T;{\mathbb R}^{d\times d})$.  We assume further that $(a_{i,j})_{d\times d}=(\sigma_{i,k}\sigma_{j,k})_{d\times d}$ is uniformly elliptic, i.e. for every $t\in [0,T]$, there is a constant $\Theta>1$ such that
\begin{eqnarray}\label{1.7}
\Theta^{-1} |\vartheta|^2\leq \vartheta^\top a(t)\vartheta \leq \Theta |\vartheta|^2, \quad
\forall \ \vartheta=(\vartheta_1,\vartheta_2,\ldots,\vartheta_d)\in{\mathbb R}^d.
\end{eqnarray}
Then it holds that:

\medskip
(i) (\textbf{Stochastic flow of homeomorphisms}) There exists a unique stochastic flow of homeomorphisms $\{X_{s,t}(x), \ t\in [s,T ]\}$ to SDE
(\ref{1.1}).

\medskip
(ii) (\textbf{Gradient and H\"{o}lder estimates}) $X_{s,t}(x)$ and $X^{-1}_{s,t}(x)$ are differentiable in $x$ in $L^2(\Omega)$ for every $0\leq s\leq t\leq T$, and for every $p\geq 2$,
\begin{eqnarray}\label{1.8}
\sup_{x\in{\mathbb R}^d}\sup_{0\leq s\leq
T}{\mathbb E}[\sup_{s\leq t \leq T}\|\nabla X_{s,t}(x)\|^p]+\sup_{x\in{\mathbb R}^d}\sup_{0\leq t\leq
T}{\mathbb E}[\sup_{0\leq s\leq t}\|\nabla X_{s,t}^{-1}(x)\|^p]
<+\infty.
\end{eqnarray}
Moreover,  $\nabla X_{s,t}(\cdot)$ and $\nabla X^{-1}_{s,t}(\cdot)$ have continuous realizations (denoted by themselves), which are locally $\beta$-H\"{o}lder continuous in $x$ for every $\beta\in(0,1+\alpha-2/q)$, and for every $p\geq 2$ and $R>0$,
\begin{eqnarray}\label{1.9}
&&\sup_{0\leq s\leq
T}{\mathbb E}\Big[\sup_{s\leq t \leq T}\Big(\sup_{x,y\in B_R,x\neq y}\frac{|\nabla X_{s,t}(x)-\nabla X_{s,t}(y)|}{|x-y|^\beta}\Big)^p\Big]\nonumber \\
&&+\sup_{0\leq t\leq
T}{\mathbb E}\Big[\sup_{0\leq s\leq t}\Big(\sup_{x,y\in B_R,x\neq y}\frac{|\nabla X_{s,t}^{-1}(x)-\nabla X_{s,t}^{-1}(y)|}{|x-y|^\beta}\Big)^p\Big]<+\infty,
\end{eqnarray}
here a random field $\nabla \tilde{X}_{s,t}(\cdot)$ is called a realization of $\nabla X_{s,t}(\cdot)$ if there exists $\Omega_0\subset\Omega$ such that ${\mathbb P}(\Omega_0)=1$ and for each $\omega\in \Omega_0$, $\nabla X_{s,t}(x,\omega)=\nabla \tilde{X}_{s,t}(x,\omega)$ for all $x\in {\mathbb R}^d$.

\medskip
(iii) (\textbf{Stability}) We assume further that $q>4/(2+\alpha)$. Let $\rho$ be a symmetric regularizing kernel, that is
\begin{eqnarray}\label{1.10}
\rho(x)=\rho(-x) \ \ {\rm with} \ \ 0\leq \rho \in {\mathcal C}^\infty_0({\mathbb R}^d) , \ \ {\rm supp}(\rho)\subset B_1, \ \ \int_{{\mathbb R}^d}\rho(x)dx=1.
\end{eqnarray}
For $n\in {\mathbb N}$, we set $\rho_n(x)=n^d \rho(nx)$ and
\begin{eqnarray}\label{1.11}
b_n(t,x)=\int_{{\mathbb R}^d}b(t,x-y)\rho_n(y)dy=:b\ast\rho_n(t,x).
\end{eqnarray}
Let $X^n$ be the  stochastic flow corresponding to the vector field $b_n$ and $X^{n,-1}$ be its inverse. Then for every $p\geq 2$,
\begin{eqnarray}\label{1.12}
&&\lim _{n\rightarrow+\infty}\sup_{x\in {\mathbb R}^d} \sup_{0\leq s\leq
T}{\mathbb E}[\sup_{s\leq t \leq T}|X^n_{s,t}(x)-X_{s,t}(x)|^p]\nonumber \\
&&+\lim _{n\rightarrow+\infty}\sup_{x\in {\mathbb R}^d}\sup_{s\leq t \leq T}
{\mathbb E}[\sup_{0\leq s\leq t}|X^{n,-1}_{s,t}(x)-X_{s,t}^{-1}(x)|^p]=0
\end{eqnarray}
and
\begin{eqnarray}\label{1.13}
&&\lim_{n\rightarrow+\infty}\sup_{x\in {\mathbb R}^d} \sup_{0\leq s\leq
T}{\mathbb E}[\sup_{s\leq t \leq T}\|\nabla X^n_{s,t}(x)-\nabla X_{s,t}(x)\|^p]\nonumber \\
&&+\lim _{n\rightarrow+\infty}\sup_{x\in {\mathbb R}^d}\sup_{s\leq t \leq T}
{\mathbb E}[\sup_{0\leq s\leq t}\|\nabla X^{n,-1}_{s,t}(x)-\nabla X_{s,t}^{-1}(x)\|^p]=0.
\end{eqnarray}
\end{theorem}

\begin{remark} \label{rem1.4} From the proof in Sec. \ref{sec3}, we also get the H\"{o}lder continuity of  $X_{s,t}(x)$ in $(s,t)$. In fact, by (\ref{3.12}) and the Kolmogorov-Chentsov continuity criterion, the random field $Y_{s,t}(y)$ has a continuous realization (denoted by itself), which is locally $(\beta_1,\beta_2,\beta_3)$-H\"{o}lder continuous in $(s,t,y)$, for every $\beta_1,\beta_2\in (0,1/2)$ and $\beta_3\in (0,1)$, i.e., there exists $\Omega_0\subset\Omega$ with ${\mathbb P}(\Omega_0)=1$ such that $Y_{s,t}(y)$ is locally $(\beta_1,\beta_2,\beta_3)$-H\"{o}lder continuous in $(s,t,y)$ for all $\omega\in \Omega_0$. Let $U$ be the unique strong solution of the Cauchy problem (\ref{3.1}). We set $\Phi(t,x)=x+U(t,x)$ and $\Psi(t,\cdot)=\Phi^{-1}(t,\cdot)$. Then
$$
\nabla \Phi\in L^\infty(0,T;{\mathcal C}_b^\theta({\mathbb R}^d;{\mathbb R}^{d\times d})) \ {\rm and} \ \partial_t\Phi \in L^q(0,T;{\mathcal C}({\mathbb R}^d;{\mathbb R}^d)), \ \ \forall \ \theta\in (0,1+\alpha-2/q).
$$
Observing that
$$
\nabla \Psi(t,x)=[\nabla \Phi(t,\Psi(t,x))]^{-1} \ {\rm and} \
\partial_t\Psi(t,x)=-[\nabla \Phi(t,\Psi(t,x))]^{-1}\partial_t\Phi(t,\Psi(t,x)),
$$
then
$$
\nabla\Psi\in L^\infty(0,T;{\mathcal C}_b^\theta({\mathbb R}^d;{\mathbb R}^{d\times d}))\ {\rm and} \ \partial_t\Psi \in L^q(0,T;{\mathcal C}({\mathbb R}^d;{\mathbb R}^d)), \ \ \forall \ \theta\in (0,1+\alpha-2/q).
$$
On the other hand, $X_{s,t}(x)=\Psi(t,Y_{s,t}(\Phi(s,x)))$, then for every $s<s^\prime$, $t<t^\prime$ and $\omega\in \Omega_0$, we have
\begin{eqnarray*}
&&|X_{s,t}(x)-X_{s^\prime,t^\prime}(x)|\nonumber \\
&\leq&|\Psi(t,Y_{s,t}(\Phi(s,x)))-
\Psi(t^\prime,Y_{s,t}(\Phi(s,x)))|+|\Psi(t^\prime,Y_{s,t}(\Phi(s,x)))-
\Psi(t^\prime,Y_{s^\prime,t^\prime}(\Phi(s^\prime,x)))|\nonumber \\
&\leq& C|t-t^\prime|^{1-\frac{1}{q}}+\|\nabla \Psi\|_{\infty,0}|Y_{s,t}(\Phi(s,x))-
Y_{s^\prime,t^\prime}(\Phi(s^\prime,x))|\nonumber \\
&\leq& C[|t-t^\prime|^{1-\frac{1}{q}}+|s-s^\prime|^{\beta_1}+|t-t^\prime|^{\beta_2}+ |\Phi(s,x)-\Phi(s^\prime,x)|^{\beta_3}]\nonumber \\
&\leq& C[|t-t^\prime|^{1-\frac{1}{q}}+|s-s^\prime|^{\beta_1}+|t-t^\prime|^{\beta_2}+ |s-s^\prime|^{(1-\frac{1}{q})\beta_3}]\leq C[|t-t^\prime|^{1-\frac{1}{q}}+|s-s^\prime|^{(1-\frac{1}{q})\beta_3}],
\end{eqnarray*}
where in the third and fifth lines we have used Sobolev imbedding theorem: $W^{1,q}(0,T;{\mathcal C}({\mathbb R}^d))$ $\hookrightarrow {\mathcal C}^{1-\frac{1}{q}}([0,T];{\mathcal C}({\mathbb R}^d))$ and in the last line we have chosen $\beta_1,\beta_2\in (1-1/q,1/2)$.
Therefore, $X_{s,t}(x)$ is $(\beta_4,1-1/q)$-H\"{o}lder continuous in $(s,t)$ for every $\beta_4\in (0,1-1/q)$.
\end{remark}

As an application, we consider the following stochastic transport equation with low regularity growing drift:
\begin{eqnarray}\label{1.14}
\left\{
  \begin{array}{ll}
\partial_tu(t,x)+b(t,x)\cdot\nabla u(t,x)
+\sum\limits_{i=1}^d\partial_{x_i}u(t,x)\circ\dot{W}_{i,t}=0, \ \ (t,x)\in(0,T]\times {\mathbb R}^d, \\
u(t,x)|_{t=0}=u_0(x), \quad  x\in{\mathbb R}^d,
  \end{array}
\right.
\end{eqnarray}
where the stochastic integral with the notation $\circ$ is interpreted in Stratonovich
sense, and others are interpreted in It\^{o}'s. The choice of Stratonovich integral in (\ref{1.14}) is motivated by mass conservation. In fact if $b$ is divergence free, we rewrite (\ref{1.14}) by
\begin{eqnarray*}
\partial_tu(t,x)+\div[(b(t,x)+\dot{W}_t)u(t,x)]=0, \
u(t,x)|_{t=0}=u_0(x),
\end{eqnarray*}
which implies
\begin{eqnarray*}
\int_{{\mathbb R}^d}u(t,x)dx=\int_{{\mathbb R}^d}u_0(x)dx, \quad \forall \ t\in [0,T], \quad {\mathbb P}-a.s..
\end{eqnarray*}

Firstly, we give the definitions of weak and strong solutions for the above Cauchy problem.
\begin{definition} \label{def1.5} Let $u_0\in L^r_{loc}({\mathbb R}^d)$, $b\in L^1(0,T;L^{r^\prime}_{loc}({\mathbb R}^d;{\mathbb R}^d))$  and ${\rm div}  b\in L^1(0,T;L^{r^\prime}_{loc}({\mathbb R}^d))$ with $r\in [1,+\infty]$ and $1/r+1/r^\prime=1$.
Let $u\in L^\infty(\Omega\times[0,T];L^r({\mathbb R}^d))$ be a random field.
We call $u$ a stochastic weak solution of (\ref{1.14}) if for every
$\varphi\in{\mathcal C}_0^\infty({\mathbb R}^d)$, $\int_{{\mathbb R}^d}u(t,x)\varphi(x)dx$
has a continuous modification which is an ${\mathcal F}_t$-semimartingale and
for every  $t\in [0,T]$,
\begin{eqnarray}\label{1.15}
\int_{{\mathbb R}^d}u(t,x)\varphi(x)dx&=&\int_{{\mathbb R}^d}u_0(x)\varphi(x)dx+
\int^t_0\int_{{\mathbb R}^d}u(\tau,x)\div (b(\tau,x)\varphi(x))dxd\tau\nonumber\\&&
+\sum_{i=1}^d\int^t_0\circ dW_{i,\tau}\int_{{\mathbb R}^d}u(\tau,x)\partial_{x_i}\varphi(x)dx,  \quad
{\mathbb P}-a.s..
\end{eqnarray}
Moreover, if the following additional estimates hold:
\begin{eqnarray}\label{1.16}
\left\{
  \begin{array}{ll}
\sup\limits_{0\leq t\leq T}{\mathbb E}\|\nabla u(t)\|_{L^r_{loc}({\mathbb R}^d)}^r<+\infty, \quad {\rm if} \ r<\infty, \\
\sup\limits_{0\leq t\leq T} {\mathbb E}\|\nabla u(t)\|_{L^\infty_{loc}({\mathbb R}^d)}^p<+\infty, \quad \forall \ p\in [2,+\infty), \quad {\rm if} \ r=+\infty,
  \end{array}
\right.
\end{eqnarray}
then $u$ is called  a stochastic strong solution.
\end{definition}

\begin{remark} \label{rem1.6}
 The stochastic integral in (\ref{1.15}) can be represented by the It\^{o} sense equivalently.
In fact, the last term in (\ref{1.15}) is equivalent to
\begin{eqnarray*}
&&\sum_{i=1}^d\int^t_0dW_{i,\tau}\int_{{\mathbb R}^d}u(\tau,x)\partial_{x_i}\varphi(x)dx+
\frac{1}{2}\sum_{i=1}^d\int_0^t \Big(d\int_{{\mathbb R}^d}u(\tau,x)\partial_{x_i}\varphi(x)dx\Big)\Big(dW_{i,\tau}\Big)
\nonumber \\&=&\sum_{i=1}^d\int^t_0dW_{i,\tau}\int_{{\mathbb R}^d}u(\tau,x)\partial_{x_i}\varphi(x)dx+
\frac{1}{2} \int_0^t\int_{{\mathbb R}^d}u(\tau,x)\Delta\varphi(x)dxd\tau.
\end{eqnarray*}
\end{remark}

Let us now state our second result.
\begin{theorem} \label{the1.7} \textbf{(Existence and uniqueness)} Suppose $r\in [1,+\infty]$ and $u_0\in L^r({\mathbb R}^d)$. Let  $b\in L^q(0,T;{\mathcal C}^{\frac{2}{q}-1}\cap{\mathcal C}^\alpha({\mathbb R}^d;{\mathbb R}^d))$ with $\alpha\in (0,1)$ and $q\in (4/(2+\alpha),2)$, and be divergence free. Then one has:

\medskip
(i) There exists a unique stochastic weak solution to the Cauchy problem (\ref{1.14}).  Moreover, the unique
stochastic weak solution can be represented by $u(t,x)=u_0(X^{-1}_t(x))$, where $X_t(x)$ is the unique
strong solution of the associated SDE (\ref{1.1}) with $\sigma=I_{d\times d}$ and $s=0$.

\medskip
(ii) The unique stochastic weak solution of (\ref{1.14}) is also the unique stochastic strong solution if $|\nabla u_0|\in L^r({\mathbb R}^d)$.
 \end{theorem}

\begin{remark} \label{rem1.8} Instead of a fixed $r$, if we  require the intersection of all $r\geq 1$'s, the well-posedness of stochastic strong solutions for (\ref{1.14}) has been established by Fedrizzi and Flandoli \cite[Theorem 1]{FF2} for the drift $b$ which is in the Krylov-R\"{o}ckner class (see (\ref{1.2}) and (\ref{1.3})). Precisely speaking, for every given initial data $u_0\in \mathop{\cap}\limits_{r\geq 1}W^{1,r}({\mathbb R}^d)$, they proved the existence and uniqueness of $\mathop{\cap}\limits_{r\geq 1}W^{1,r}_{loc}({\mathbb R}^d)$-solutions, i.e., for every $t\in [0,T]$,
\begin{eqnarray*}
{\mathbb P}\Big(u(t)\in \mathop{\cap}\limits_{r\geq 1}W^{1,r}_{loc}({\mathbb R}^d)\Big)=1.
\end{eqnarray*}
After the result of \cite[Theorem 1]{FF2}, it remains to open the question whether the solution is Lipschitz
continuous (or more) when $u_0\in W^{1,\infty}({\mathbb R}^d)$ (or more) for irregular drift.

\medskip
$\bullet$ When $d=1$, the answer for the above question is positive for certain discontinuous drift $b$, including for instance $b(x)=\mbox{sign}(x)$, see \cite{Att}.

\medskip
$\bullet$ When $b\in L^q(0,T;\mathcal{C}_b^\alpha({\mathbb R}^d;{\mathbb R}^d))$ with $\alpha\in (0,1)$ and $q>2/\alpha$, it is positive as well, see \cite{WDGL} (or see \cite{FGP1} if $q=+\infty$).

\medskip
$\bullet$ When $b$ is in the Krylov-R\"{o}ckner class, by virtue of \cite[Corollary 1.1]{Rezakhanlou}, we also get the local quasi-Lipschitz estimate for stochastic strong solutions, i.e.
\begin{eqnarray*}
{\mathbb P}\Big\{\sup_{x,y\in B_R,x\neq y}\frac{|u(t,x)-u(t,y)|}{|x-y|\exp\Big(C(d,T,p,q,R)\big(\log\frac{R}{|x-y|}
\big)^{\nu}
\Big)}<+\infty\Big\}=1,
\end{eqnarray*}
where $2\nu=1+d/p+2/q$. However, it is still unknown whether the solution is (locally) Lipschitz
continuous or not.

\medskip
$\bullet$ Now, when the drift is in a low regularity Lebesgue-H\"{o}lder space  $L^q(0,T;
{\mathcal C}^{\frac{2}{q}-1} \cap{\mathcal C}^{\alpha}({\mathbb R}^d;{\mathbb R}^d))$ with $\alpha\in (0,1)$ and $q\in (4/(2+\alpha),2)$, if $u_0\in W^{1,r}({\mathbb R}^d)$ for $r\in [1,+\infty]$, by (\ref{1.16}) we also get
\begin{eqnarray}\label{1.17}
{\mathbb P}\Big(u(t)\in W^{1,r}_{loc}({\mathbb R}^d)\Big)=1.
\end{eqnarray}
In particular, (\ref{1.17}) holds for $r=+\infty$, and so we give a positive answer for the above question for this low regularity growing drift.
\end{remark}

\noindent
\textbf{Notations.} The letter $C$  denotes a positive constant, whose values may change in different places. For a parameter or a function $\tilde{\nu}$, $C(\tilde{\nu})$ means the constant is only dependent on $\tilde{\nu}$, and we also write it as $C$ if there is no confusion.  We
use $\nabla$ to denote the gradient of a function with respect to the space variable.  As usual, ${\mathbb N}$ stands for
the set of all natural numbers.  a.s. is the abbreviation of almost surely. For every $R>0$, $B_R:=\{x\in{\mathbb R}^d:|x|\leq R\}$. For a given ${\mathbb R}^{n\times m}$ matrix-valued function $\Xi$ with $0<n,m\in {\mathbb N}$, $\Xi^{\top}$ and $\|\Xi\|$ represent its  transposition and Hilbert-Schmidt norm, respectively. If $n=m$, $tr(\Xi)$ stands for the trace of $\Xi$.

\section{Lebesgue-Schauder estimates for Kolmogorov equations}\label{sec2}
\setcounter{equation}{0}
Let $b: [0,T]\times{\mathbb R}^d\rightarrow{\mathbb R}^d$ and $f: [0,T]\times{\mathbb R}^d\rightarrow{\mathbb R}$ be Borel measurable functions. Consider the following Kolmogorov equation for $u:[0,T]\times{\mathbb R}^d\to{\mathbb R}$:
\begin{eqnarray}\label{2.1}
\left\{\begin{array}{ll}
\partial_{t}u(t,x)=\frac{1}{2}\Delta u(t,x)+b(t,x)\cdot \nabla u(t,x)
\\ \qquad\qquad \ \ -\lambda u(t,x)+f(t,x), \ (t,x)\in (0,T]\times {\mathbb R}^d, \\
u(t,x)|_{t=0}=0, \  x\in{\mathbb R}^d,  \end{array}\right.
\end{eqnarray}
where $\lambda>0$ is a real number. If $u\in
L^1(0,T;{\mathcal C}^2({\mathbb R}^d))\cap W^{1,1}(0,T;{\mathcal C}({\mathbb R}^d))$ such that (\ref{2.1}) holds true for almost all $(t,x)\in (0,T)\times {\mathbb R}^d$, then the unknown function $u$ is said to be a strong solution of (\ref{2.1}).

\medskip
Now let us establish the well-posedness of strong solutions for the Cauchy problem (\ref{2.1}).
\begin{lemma} \label{lem2.1} Assume that $b\in L^q(0,T;
{\mathcal C}^{\frac{2}{q}-1}\cap{\mathcal C}^\alpha({\mathbb R}^d;{\mathbb R}^d))$ and  $f\in L^q(0,T;{\mathcal C}^{\frac{2}{q}-1}\cap{\mathcal C}^\alpha({\mathbb R}^d))$ with $\alpha\in (0,1)$ and $q\in (2/(1+\alpha),2)$. Then it holds that:

\medskip
(i) \textbf{(Existence and uniqueness)} There is a unique strong solution $u$ to (\ref{2.1}).
Moreover, $u\in {\mathcal H}^{2,\theta}_{q,T}$ for every $\theta\in (0, 1+\alpha-2/q)$, where
\begin{eqnarray}\label{2.2}
{\mathcal H}^{2,\theta}_{q,T}&=&\{v\in L^\infty(0,T;{\mathcal C}^{\frac{2}{q}-1}({\mathbb R}^d)); \  \nabla v\in L^\infty(0,T;{\mathcal C}_b^\theta({\mathbb R}^d;{\mathbb R}^d)), \ \partial_tv\in L^q(0,T;{\mathcal C}({\mathbb R}^d)),\nonumber \\ && \quad  \nabla^2v\in L^2(0,T;{\mathcal C}_b^\theta({\mathbb R}^d;{\mathbb R}^{d\times d})) \ {\rm and} \ \|(1+|\cdot|^{\frac{2}{q}-1})^{-1}\partial_tv(\cdot,\cdot)\|_{q,0}<+\infty\}.
\end{eqnarray}
Furthermore, there is a real number $\varepsilon>0$ such that for large enough $\lambda>0$
\begin{eqnarray}\label{2.3}
\sup_{0\leq t\leq T}\|\nabla u(t)\|_{{\mathcal C}_b^\theta({\mathbb R}^d)}\leq C(d,T,[b]_{q,\frac{2}{q}-1},[b]_{q,\alpha})\lambda^{-\varepsilon}
([f]_{q,\frac{2}{q}-1}+[f]_{q,\alpha})
\end{eqnarray}
and
\begin{eqnarray}\label{2.4}
\|\nabla ^2u\|_{L^2(0,T;{\mathcal C}^\theta_b({\mathbb R}^d))}
\leq C(d,T,[b]_{q,\frac{2}{q}-1},[b]_{q,\alpha})\lambda^{-\varepsilon}
([f]_{q,\frac{2}{q}-1}+[f]_{q,\alpha}).
\end{eqnarray}

(ii) \textbf{(Stability)} Let $u_n$ be the unique strong solution of (\ref{2.1}) with $b$ and $f$ replaced by $b_n$ and $f_n$, where $b_n$ is given by (\ref{1.11}) and
$$
f_n(t,x)=\int_{{\mathbb R}^d}f(t,x-y)\rho_n(y)dy=:f\ast\rho_n(t,x).
$$
Then $u_n$ belongs to ${\mathcal H}^{2,\theta}_{q,T}$ for every $\theta\in (0, 1+\alpha-2/q)$ and satisfies (\ref{2.3})--(\ref{2.4}) uniformly in $n$. Furthermore, $u_n-u\in L^2(0,T;{\mathcal C}_b^{2,\theta}({\mathbb R}^d))$ $\cap L^\infty(0,T;{\mathcal C}_b^{1,\theta}({\mathbb R}^d))$, and if we assume in addition that $q>4/(2+\alpha)$, then for every $\theta\in (0, 1+\alpha-2/q)$ we have
\begin{eqnarray}\label{2.5}
\lim_{n\rightarrow+\infty}\Big[\sup_{0\leq t\leq T}\|u_n(t)-u(t)\|_{{\mathcal C}_b^{1,\theta}({\mathbb R}^d)} +\|\nabla^2u_n-\nabla^2u\|_{L^2(0,T;{\mathcal C}_b^\theta({\mathbb R}^d))}\Big]=0.
\end{eqnarray}
\end{lemma}
\noindent
\textbf{Proof.}  (i) Clearly, if $u\in
L^1(0,T;{\mathcal C}^2({\mathbb R}^d))\cap W^{1,1}(0,T;{\mathcal C}({\mathbb R}^d))$ solves the Cauchy problem (\ref{2.1}) for $\lambda\geq \lambda_0$ with some sufficiently large real number $\lambda_0$, then for all $\lambda>0$, $\tilde{u}(t,x)=u(t,x)e^{(\lambda_0-\lambda)t}\in
L^1(0,T;{\mathcal C}^2({\mathbb R}^d))\cap W^{1,1}(0,T;{\mathcal C}({\mathbb R}^d))$ solves  the following Cauchy problem:
\begin{eqnarray*}
\left\{\begin{array}{ll}
\partial_{t}\tilde{u}(t,x)=\frac{1}{2}\Delta \tilde{u}(t,x)+b(t,x)\cdot\nabla \tilde{u}(t,x)\\ \qquad\qquad \ \ -\lambda \tilde{u}(t,x)+\tilde{f}(t,x), \ \ (t,x)\in (0,T]\times{\mathbb R}^d , \\
\tilde{u}(t,x)|_{t=0}=0, \  x\in{\mathbb R}^d,
\end{array}\right.
\end{eqnarray*}
where $\tilde{f}(t,x)=f(t,x)e^{(\lambda_0-\lambda)t}$, and vice versa. So we just  need to prove the well-posedness of (\ref{2.1}) for some sufficiently large $\lambda$.

\medskip
On the other hand, if $u\in  {\mathcal H}^{2,\theta}_{q,T}$ for every $\theta\in (0, 1+\alpha-2/q)$ is a strong solution of (\ref{2.1}), then it has the following equivalent representation (see \cite[Lemma 2.1]{TDW}):
\begin{eqnarray}\label{2.6}
u(t,x)&=&\int_0^te^{-\lambda (t-\tau)}K(t-\tau,\cdot)\ast (b(\tau,\cdot)\cdot \nabla u(\tau,\cdot))(x)d\tau\nonumber\\ &&+\int_0^te^{-\lambda (t-\tau)}K(t-\tau,\cdot)\ast f(\tau,\cdot)(x)d\tau,
\end{eqnarray}
where $K(t-\tau,x)=(2\pi (t-\tau))^{-\frac{d}{2}}e^{-\frac{|x|^2}{2(t-\tau)}}$. Thus, it suffices to show the integral equation (\ref{2.6}) has a unique strong solution $u\in {\mathcal H}^{2,\theta}_{q,T}$.

\medskip
Firstly, let us prove the existence part and for simplicity sake, we set
\begin{eqnarray}\label{2.7}
\tilde{{\mathcal H}}^{2,\theta}_{q,T}&=&\{v\in L^\infty(0,T;{\mathcal C}^{\frac{2}{q}-1}({\mathbb R}^d)); \ \nabla v\in L^\infty(0,T;{\mathcal C}_b^\theta({\mathbb R}^d;{\mathbb R}^d)),\nonumber \\ &&\quad \nabla^2v\in L^2(0,T;{\mathcal C}_b^\theta({\mathbb R}^d;{\mathbb R}^{d\times d})).
\end{eqnarray}
The proof is divided into three steps.

\medskip
\textbf{Step 1} We assume that $b\in L^\infty(0,T;{\mathcal C}^{\infty}_b({\mathbb R}^d;{\mathbb R}^d))$ and $f\in L^\infty(0,T;{\mathcal C}^{\infty}_b({\mathbb R}^d))$. If $b\equiv 0$,
there is a unique strong solution $u\in L^\infty(0,T;{\mathcal C}_b^\infty({\mathbb R}^d))$ of (\ref{2.6}). For $b\neq 0$, we define a mapping on $L^\infty(0,T;{\mathcal C}_b^\infty({\mathbb R}^d))$ by
\begin{eqnarray}\label{2.8}
{\mathcal T}w(t,x)&=&\int_0^te^{-\lambda (t-\tau)}K(t-\tau,\cdot)\ast (b(\tau,\cdot)\cdot \nabla w(\tau,\cdot))(x)d\tau\nonumber\\ &&+\int_0^te^{-\lambda (t-\tau)}K(t-\tau,\cdot)\ast f(\tau,\cdot)(x)d\tau.
\end{eqnarray}
From (\ref{2.8}), for every fixed $0<k\in {\mathbb N}$ and every $w_1,w_2\in L^\infty(0,T;{\mathcal C}_b^\infty({\mathbb R}^d))$, then
\begin{eqnarray}\label{2.9}
&&\|{\mathcal T} w_1-{\mathcal T} w_2\|_{L^\infty(0,T;{\mathcal C}_b^k({\mathbb R}^d))}\nonumber\\&=&\Big\|\int_0^te^{-\lambda (t-\tau)} K(t-\tau,\cdot)\ast [b(\tau,\cdot)\cdot \nabla (w_1(\tau,\cdot)-w_2(\tau,\cdot))](x)d\tau\Big\|_{\infty,0}\nonumber\\ &&+\sum_{i=1}^k\Big\|\int_0^te^{-\lambda (t-\tau)}\nabla K(t-\tau,\cdot)\ast \nabla^{i-1}[b(\tau,\cdot)\cdot \nabla (w_1(\tau,\cdot)-w_2(\tau,\cdot))](x)d\tau\Big\|_{\infty,0}\nonumber\\ &\leq &C\lambda^{-\varepsilon}\|b\|_{L^\infty(0,T;{\mathcal C}_b^{k-1}({\mathbb R}^d;{\mathbb R}^d))}\|w_1- w_2\|_{L^\infty(0,T;{\mathcal C}_b^k({\mathbb R}^d))},
\end{eqnarray}
where $\varepsilon\in (0,1/2)$.

\medskip
By choosing $\lambda$ sufficiently large, from (\ref{2.9}) the mapping ${\mathcal T}$ on $L^\infty(0,T;{\mathcal C}_b^k({\mathbb R}^d))$ is contractive. With the aid of Banach's contraction mapping principle, there is a unique strong solution $u\in L^\infty(0,T;{\mathcal C}_b^k({\mathbb R}^d))$ of (\ref{2.6}). Since $k$ is arbitrary, we get $u\in L^\infty(0,T;{\mathcal C}_b^\infty({\mathbb R}^d))$.

\medskip
Let $x_0\in {\mathbb R}^d$. Consider the following differential equation:
\begin{eqnarray}\label{2.10}
\dot{x}_t=-b(t,x_0+x_t), \quad x_t|_{t=0}=0.
\end{eqnarray}
There exists a unique solution to (\ref{2.10}) for $b\in L^\infty(0,T;{\mathcal C}^{\infty}_b({\mathbb R}^d;{\mathbb R}^d))$.
By setting $\hat{u}(t,x):=u(t,x+x_0+x_t)$, $\hat{b}(t,x):=b(t,x+x_0+x_t)-b(t,x_0+x_t)$ and $\hat{f}(t,x):=f(t,x+x_0+x_t)$, then
\begin{eqnarray}\label{2.11}
\left\{\begin{array}{ll}
\partial_{t}\hat{u}(t,x)=\frac{1}{2}\Delta \hat{u}(t,x)+\hat{b}(t,x)\cdot \nabla \hat{u}(t,x) \\ \qquad\qquad \ \ -\lambda \hat{u}(t,x)+\hat{f}(t,x), \ (t,x)\in (0,T]\times {\mathbb R}^d, \\
\hat{u}(t,x)|_{t=0}=0, \  x\in{\mathbb R}^d,  \end{array}\right.
\end{eqnarray}
which also implies
\begin{eqnarray}\label{2.12}
\hat{u}(t,x) &=&\int_0^te^{-\lambda (t-\tau)}d\tau\int_{{\mathbb R}^d}K(t-\tau,x-y)[\hat{b}(\tau,y)\cdot \nabla \hat{u}(\tau,y)]dy\nonumber\\ &&+\int_0^te^{-\lambda (t-\tau)}d\tau\int_{{\mathbb R}^d}K(t-\tau,x-y) \hat{f}(\tau,y)dy.
\end{eqnarray}
For every $\theta\in (0, 1+\alpha-2/q)$, then
\begin{eqnarray}\label{2.13}
|\hat{b}(\tau,y)\cdot \nabla \hat{u}(\tau,y)|\leq [b(\tau)]_{\alpha-\theta} |y|^{\alpha-\theta} \|\nabla u(\tau)\|_0.
\end{eqnarray}
It follows from (\ref{2.12}) and (\ref{2.13}) that
\begin{eqnarray}\label{2.14}
|\nabla \hat{u}(t,0)|&\leq&\int_0^te^{-\lambda (t-\tau)}d\tau\int_{{\mathbb R}^d}|\nabla K(t-\tau,y)|[b(\tau)]_{\alpha-\theta} |y|^{\alpha-\theta} \|\nabla u(\tau)\|_0dy\nonumber\\ &&+\int_0^te^{-\lambda (t-\tau)}d\tau\int_{{\mathbb R}^d}|\nabla K(t-\tau,y)||\hat{f}(\tau,y)-\hat{f}(\tau,0)|dy\nonumber\\ &\leq& C(d)\int_0^te^{-\lambda (t-\tau)}(t-\tau)^{\frac{\alpha-\theta-1}{2}}\Big([b(\tau)]_{\alpha-\theta} \|\nabla u(\tau)\|_0+[f(\tau)]_{\alpha-\theta}\Big)
d\tau.
\end{eqnarray}
Since $x_0\in {\mathbb R}^d$ is arbitrary, we conclude
\begin{eqnarray}\label{2.15}
\sup_{0\leq t\leq T}\|\nabla u(t)\|_0\leq C(d)
\Big([b]_{q,\alpha-\theta}\sup_{0\leq \tau\leq T}\|\nabla u(\tau)\|_0 +[f]_{q,\alpha-\theta}\Big)
\Big(\int_0^Te^{-\lambda q^\prime \tau}\tau^{\frac{(\alpha-\theta-1)q^\prime}{2}}d\tau\Big)^{\frac{1}{q^\prime}},
\end{eqnarray}
where $q^\prime=q/(q-1)$.

\medskip
Notice that $q>2/(1+\alpha)$ and $\theta<1+\alpha-2/q$, we assert that $(1+\theta-\alpha)q^\prime/2<1$.
By applying the H\"{o}lder inequality, there exists some real number $\varepsilon>0$ such that
\begin{eqnarray}\label{2.16}
C(d)\Big(\int_0^Te^{-\lambda q^\prime \tau}\tau^{\frac{(\alpha-\theta-1)q^\prime}{2}}d\tau\Big)^{\frac{1}{q^\prime}}\leq C(d,T)\lambda^{-\varepsilon}.
\end{eqnarray}
We choose $\lambda$ large enough such that
\begin{eqnarray}\label{2.17}
C(d,T)[b]_{q,\alpha-\theta}\lambda^{-\varepsilon}\leq C(d,T)[b]_{q,\frac{2}{q}-1}^\nu[b]_{q,\alpha}^{1-\nu}\lambda^{-\varepsilon}\leq C(d,T)\max\{[b]_{q,\frac{2}{q}-1},[b]_{q,\alpha}\}\lambda^{-\varepsilon}<\frac{1}{2},
\end{eqnarray}
where $\nu=\theta/(1+\alpha-2/q)$, and in the first inequality we have used the fact $\alpha-\theta\in (2/q-1,\alpha)$ and the interpolation inequality, and in the second inequality we used the Young inequality
$$
[b]_{q,\frac{2}{q}-1}^\nu[b]_{q,\alpha}^{1-\nu}\leq \nu[b]_{q,\frac{2}{q}-1}+(1-\nu)[b]_{q,\alpha}.
$$

\medskip
Then by (\ref{2.15})--(\ref{2.17}), it yields that
\begin{eqnarray}\label{2.18}
\sup_{0\leq t\leq T}\|\nabla u(t)\|_0\leq C(d,T,[b]_{q,\alpha-\theta})\lambda^{-\varepsilon}
[f]_{q,\alpha-\theta},
\end{eqnarray}
which also implies that
\begin{eqnarray}\label{2.19}
\sup_{0\leq t\leq T}\|\nabla u(t)\|_0\leq C(d,T,[b]_{q,\frac{2}{q}-1},[b]_{q,\alpha})\lambda^{-\varepsilon}
([f]_{q,\frac{2}{q}-1}+[f]_{q,\alpha}),
\end{eqnarray}
if one uses the interpolation and Young inequalities again.

\medskip
Let
$$
b^1(t)=\sup_{x\in {\mathbb R}^d}\frac{|b(t,x)|}{1+|x|^{\frac{2}{q}-1}}, \ \ f^1(t)=\sup_{x\in {\mathbb R}^d}\frac{|f(t,x)|}{1+|x|^{\frac{2}{q}-1}}.
$$
Then $b^1,f^1\in L^q(0,T)$.

\medskip
Combining (\ref{2.6}) and (\ref{2.19}), it follows that
\begin{eqnarray}\label{2.20}
&&\sup_{0\leq t\leq T}\|(1+|\cdot|^{\frac{2}{q}-1})^{-1} u(t,\cdot)\|_0
\nonumber \\ &=& \sup_{0\leq t\leq T}\sup_{x\in {\mathbb R}^d}\frac{1}{1+|x|^{\frac{2}{q}-1}}\Big|\int_0^te^{-\lambda (t-\tau)}d\tau\int_{{\mathbb R}^d}K(t-\tau,y)b(\tau,x-y)\cdot \nabla u(\tau,x-y)dy\nonumber\\ &&+\int_0^te^{-\lambda (t-\tau)}d\tau\int_{{\mathbb R}^d}K(t-\tau,y)f(\tau,x-y)dy\Big|
\nonumber\\&\leq&
\sup_{0\leq t\leq T}\sup_{x\in {\mathbb R}^d}\frac{1}{1+|x|^{\frac{2}{q}-1}}\Big[\int_0^te^{-\lambda (t-\tau)}d\tau\int_{{\mathbb R}^d}K(t-\tau,y)b^1(\tau)(1+|x|^{\frac{2}{q}-1}+|y|^{\frac{2}{q}-1})
\nonumber\\ &&\quad
\times\|\nabla u\|_{\infty,0}dy+\int_0^te^{-\lambda (t-\tau)}d\tau\int_{{\mathbb R}^d}K(t-\tau,y)f^1(\tau)(1+|x|^{\frac{2}{q}-1}+|y|^{\frac{2}{q}-1})dy\Big]
\nonumber\\&\leq& C(d,T,[b]_{q,\frac{2}{q}-1},[b]_{q,\alpha})\lambda^{-\varepsilon}\Big[
([f]_{q,\frac{2}{q}-1}+[f]_{q,\alpha})\|b^1\|_q+ \|f^1\|_q\Big]\nonumber\\&\leq&  C(d,T,\|b\|_{L^q(0,T;{\mathcal C}^{\frac{2}{q}-1}\cap{\mathcal C}^\alpha({\mathbb R}^d;{\mathbb R}^d))})\lambda^{-\varepsilon}
\|f\|_{L^q(0,T;{\mathcal C}^{\frac{2}{q}-1}\cap{\mathcal C}^\alpha({\mathbb R}^d))},
\end{eqnarray}
where $\|b^1\|_q=\|b^1\|_{L^q(0,T)}, \ \|f^1\|_q=\|f^1\|_{L^q(0,T)}$.

\medskip
For $\nabla^2 u$, we estimate from (\ref{2.10}) to (\ref{2.13}) that
\begin{eqnarray}\label{2.21}
\|\nabla^2 u(t)\|_0\leq C(d)\int_0^te^{-\lambda (t-\tau)}(t-\tau)^{\frac{\alpha-\theta}{2}-1}\Big([b(\tau)]_{\alpha-\theta} \|\nabla u(\tau)\|_0+[f(\tau)]_{\alpha-\theta}\Big)
d\tau,
\end{eqnarray}
which leads to
\begin{eqnarray}\label{2.22}
\|\nabla^2 u\|_{2,0}&\leq&
 C\Big([b]_{q,\alpha-\theta}\sup_{0\leq \tau\leq T}\|\nabla u(\tau)\|_0+[f]_{q,\alpha-\theta}\Big)\Big[\int_0^Te^{-\frac{2\lambda q\tau}{3q-2}}\tau^{\frac{(\alpha-\theta-2)q}{3q-2}}d\tau \Big]^{\frac{3q-2}{2q}}\nonumber\\&\leq&
 C(d,T,[b]_{q,\alpha-\theta})[f]_{q,\alpha-\theta},
\end{eqnarray}
where in the last inequality we have used $(2+\theta-\alpha)q/(3q-2)<1$ since $\theta<1+\alpha-2/q$.

\medskip
Let $P_\xi$ be given by (\ref{1.5}). We set $v_\xi(t,x)=\partial_\xi P_\xi u(t,x)$, then
\begin{eqnarray}\label{2.23}
\left\{\begin{array}{ll}
\partial_{t}v_\xi(t,x)=\frac{1}{2}\Delta v_\xi(t,x)+b(t,x)\cdot\nabla v_\xi(t,x)
\\ \qquad\qquad\quad-\lambda v_\xi(t,x)+g_\xi(t,x), \ \ (t,x)\in (0,T]\times{\mathbb R}^d, \\
v_\xi(t,x)|_{t=0}=0, \  x\in{\mathbb R}^d,
\end{array}\right.
\end{eqnarray}
where
$$
g_\xi(t,x)=\partial_\xi P_\xi f(t,x)+\partial_\xi P_\xi(b(t,x)\cdot \nabla u(t,x))-b(t,x)\cdot \partial_\xi P_\xi \nabla u(t,x).
$$
Owing to (\ref{2.18}), for every fixed $\theta\in (0,1+\alpha-2/q)$, we achieve
\begin{eqnarray}\label{2.24}
\sup_{0\leq t\leq T}\|\nabla v_\xi(t)\|_0 &\leq&C\lambda^{-\varepsilon}
[g_\xi]_{q,\alpha-\theta}\nonumber\\&\leq& C\lambda^{-\varepsilon}
\Big([\partial_\xi P_\xi f]_{q,\alpha-\theta}+[\partial_\xi P_\xi(b\cdot \nabla u)-b\cdot  \partial_\xi P_\xi \nabla u]_{q,\alpha-\theta}\Big).
\end{eqnarray}
By \cite[Lemma 2.1]{CSZ}, for every $0<\beta\leq \alpha<1$, there exists a positive constant $C(\alpha,\beta,d)$ such that
\begin{eqnarray*}
[\partial_\xi P_\xi (h_1h_2)-h_1\partial_\xi P_\xi h_2]_{\alpha-\beta}\leq C(\alpha,\beta,d)[h_1]_\alpha\|h_2\|_0\xi^{\beta-1},
\end{eqnarray*}
if $[h_1]_\alpha$ and $\|h_2\|_0$ are finite, which implies
\begin{eqnarray}\label{2.25}
[\partial_\xi P_\xi f(t)]_{\alpha-\theta}\leq C(\alpha,\theta,d)[f(t)]_\alpha\xi^{\theta-1}
\end{eqnarray}
and
\begin{eqnarray}\label{2.26}
[\partial_\xi P_\xi(b(t)\cdot \nabla u(t))-b(t)\cdot  \partial_\xi P_\xi \nabla u(t)]_{\alpha-\theta}\leq C[b(t)]_\alpha\|\nabla u(t)\|_0\xi^{\theta-1}.
\end{eqnarray}

\medskip
Combining (\ref{2.24})--(\ref{2.26}) and (\ref{2.19}), for every $\theta\in (0, 1+\alpha-2/q)$ we deduce
\begin{eqnarray}\label{2.27}
\sup_{0\leq t\leq T}\|\nabla v_\xi(t)\|_0&\leq& C\lambda^{-\varepsilon}
\Big([f]_{q,\alpha}+[b]_{q,\alpha} \sup_{0\leq t\leq T}\|\nabla u(t)\|_0 \Big)\xi^{\theta-1} \nonumber\\&\leq& C\lambda^{-\varepsilon}([f]_{q,\frac{2}{q}-1}+[f]_{q,\alpha})\xi^{\theta-1}.
\end{eqnarray}
By  (\ref{2.19}), (\ref{2.27}) and (\ref{1.6}), then $\nabla u\in L^\infty(0,T;{\mathcal C}_b^\theta({\mathbb R}^d;{\mathbb R}^d))$ and (\ref{2.3}) holds.

\medskip
Similarly, by (\ref{2.21}), (\ref{2.23}) and (\ref{2.25})--(\ref{2.27}), we guarantee
\begin{eqnarray}\label{2.28}
\|\nabla^2  v_\xi(t)\|_0&\leq& C\int_0^te^{-\lambda (t-\tau)}(t-\tau)^{\frac{\alpha-\theta}{2}-1}\Big([b(\tau)]_{\alpha-\theta}\sup_{0\leq \tau\leq T}\|\nabla v_\xi(\tau)\|_0
\nonumber\\&&
+
[\partial_\xi P_\xi f(\tau)]_{\alpha-\theta}+[\partial_\xi P_\xi(b(\tau)\cdot \nabla u(\tau))-b(\tau)\cdot  \partial_\xi P_\xi \nabla u(\tau)]_{\alpha-\theta}\Big)d\tau
\nonumber\\ &\leq& C\int_0^te^{-\lambda (t-\tau)}(t-\tau)^{\frac{\alpha-\theta}{2}-1}\Big[([b(\tau)]_{\frac{2}{q}-1}+[b(\tau)]_{\alpha})
\nonumber\\ && \quad \times ([f]_{q,\frac{2}{q}-1}+[f]_{q,\alpha})+[f(\tau)]_{\alpha})d\tau\xi^{\theta-1}.
\end{eqnarray}
On account of (\ref{2.18}) and Hausdorff-Young's inequality, we get from (\ref{2.28}) that
\begin{eqnarray}\label{2.29}
&&\|\nabla ^2u\|_{L^2(0,T;{\mathcal C}^\theta_b({\mathbb R}^d))}
\nonumber\\&=&\Big[\|\nabla ^2u\|_{2,0}^2+\int_0^T\|\xi^{1-\theta}\nabla^2  v_\xi(t)\|_{L^\infty({\mathbb R}_+\times {\mathbb R}^d)}^2 dt\Big]^{\frac{1}{2}}
\nonumber\\&\leq&  C(d,T,[b]_{q,\frac{2}{q}-1},[b]_{q,\alpha})
([f]_{q,\frac{2}{q}-1}+[f]_{q,\alpha})\Big[\int_0^Te^{-\frac{2\lambda q\tau}{3q-2}}\tau^{\frac{(\alpha-\theta-2)q}{3q-2}}d\tau \Big]^{\frac{3q-2}{2q}}\nonumber\\&\leq& C(d,T,[b]_{q,\frac{2}{q}-1},[b]_{q,\alpha})\lambda^{-\varepsilon}
([f]_{q,\frac{2}{q}-1}+[f]_{q,\alpha}),
\end{eqnarray}
which implies (\ref{2.4}). Summing over (\ref{2.19}), (\ref{2.20}), (\ref{2.27}) and (\ref{2.29}), then $u\in \tilde{{\mathcal H}}^{2,\theta}_{q,T}$ for every $\theta\in (0, 1+\alpha-2/q)$.

\medskip
\textbf{Step 2.} We assume that $b\in L^q(0,T;{\mathcal C}^{\alpha}_b({\mathbb R}^d;{\mathbb R}^d))$ and $f\in L^q(0,T;{\mathcal C}^{\alpha}_b({\mathbb R}^d))$. Firstly, we extend them from $[0,T]$ to $(-\infty,T]$ and define them by
$$
b(t,x)=b(0,x), \quad f(t,x)=f(0,x), \quad {\rm if}  \ t<0.
$$
Let the regularizing kernel $\rho$ be given by (\ref{1.10}), and let $\varrho$ be another regularizing kernel:
\begin{eqnarray}\label{2.30}
0\leq \varrho \in {\mathcal C}^\infty_0({\mathbb R}) , \ \ \, {\rm supp}(\varrho)\subset [0,1], \ \int_{{\mathbb R}}\varrho(t)dt=1.
 \end{eqnarray}
For $n,m\in{\mathbb N}$, we set $\rho_n(x)=n^d\rho(nx)$ and $\varrho_m(t)=m\varrho(mt)$. We then smooth $b$ and $f$ by $\rho_n$ and $\varrho_m$:
$$
b_{n,m}(t,x)=(b(\cdot,\cdot)\ast\rho_n\ast \varrho_m)(t,x)=\int_{{\mathbb R}^{d+1}}b(t-\tau,x-y)\rho_n(y)\varrho_m(\tau)dyd\tau
$$
and
$$
f_{n,m}(t,x)=(f(\cdot,\cdot)\ast\rho_n\ast \varrho_m)(t,x)=\int_{{\mathbb R}^{d+1}}f(t-\tau,x-y)\rho_n(y)\varrho_m(\tau)dyd\tau.
$$
Then $b_{n,m}\in L^\infty(0,T;{\mathcal C}_b^\infty({\mathbb R}^d;{\mathbb R}^d))$ and $f_{n,m}\in L^\infty(0,T;{\mathcal C}_b^\infty({\mathbb R}^d))$. Moreover,
\begin{eqnarray}\label{2.31}
\left\{
  \begin{array}{ll}
\|(1+|\cdot|^{\frac{2}{q}-1})^{-1}b_{n,m}\|_{q,0}\leq 2\|(1+|\cdot|^{\frac{2}{q}-1})^{-1}b\|_{q,0},\\ [0.3cm]
\ [b_{n,m}]_{q,\alpha}\leq [b]_{q,\alpha}, \ [f_{n,m}]_{q,\frac{2}{q}-1}\leq [b]_{q,\frac{2}{q}-1},  \\ [0.3cm]
 \|(1+|\cdot|^{\frac{2}{q}-1})^{-1}f_{n,m}\|_{q,0}\leq 2\|(1+|\cdot|^{\frac{2}{q}-1})^{-1}f\|_{q,0},\\ [0.3cm]
\ [f_{n,m}]_{q,\alpha}\leq [f]_{q,\alpha}, \ [f_{n,m}]_{q,\frac{2}{q}-1}\leq [f]_{q,\frac{2}{q}-1}.   \end{array}
\right.
\end{eqnarray}

Notice that $b_{n,m}-b=b_{n,m}-b_n+b_n-b$ with $b_n=b(t,\cdot)\ast\rho_n(x)$, and
\begin{eqnarray}\label{2.32}
|b_n(t,x)-b(t,x)|\leq \int_{{\mathbb R}^d}|b(t,x-z)-b(t,x)|\rho_n(z)dz\leq [b(t)]_\alpha\int_{{\mathbb R}^d}|z|^\alpha\rho_n(z)dz,
\end{eqnarray} thus $b_n-b\rightarrow 0$ in $L^q(0,T;{\mathcal C}_b({\mathbb R}^d))$. This fact is also true for $f_n-f$. Therefore,
\begin{eqnarray}\label{2.33}
 \lim_{n\rightarrow+\infty}\lim_{m\rightarrow+\infty}[\|b_{n,m}-b\|_{L^q(0,T;{\mathcal C}_b^{\alpha^\prime}({\mathbb R}^d;{\mathbb R}^d))}+\|f_{n,m}-f\|_{L^q(0,T;{\mathcal C}_b^{\alpha^\prime}({\mathbb R}^d))}]=0,
\end{eqnarray}
for every $\alpha^\prime\in (0,\alpha)$.

\medskip
Consider the following Kolmogorov equation:
\begin{eqnarray}\label{2.34}
\left\{\begin{array}{ll}
\partial_{t}u_{n,m}(t,x)=\frac{1}{2}\Delta u_{n,m}(t,x)+b_{n,m}(t,x)\cdot \nabla u_{n,m}(t,x)
\\ \ \qquad\qquad\qquad -\lambda u_{n,m}(t,x)+f_{n,m}(t,x), \ (t,x)\in (0,T]\times {\mathbb R}^d, \\
u_{n,m}(t,x)|_{t=0}=0, \  x\in{\mathbb R}^d.
\end{array}\right.
\end{eqnarray}
By \textbf{Step 1}, there exists a unique strong solution $u_{n,m}$ to the Cauchy problem (\ref{2.34}). Moreover,
$u_{n,m}\in \tilde{{\mathcal H}}^{2,\theta}_{q,T}$ for every $\theta\in (0, 1+\alpha-2/q)$, and by (\ref{2.3}), (\ref{2.4}), (\ref{2.20}) and (\ref{2.31}), for large enough $\lambda>0$ we have
\begin{eqnarray}\label{2.35}
&&\sup_{0\leq t\leq T}\|(1+|\cdot|^{\frac{2}{q}-1})^{-1} u_{n,m}(t,\cdot)\|_0
\nonumber\\&\leq& C(d,T,\|b_{n,m}\|_{L^q(0,T;{\mathcal C}^{\frac{2}{q}-1}\cap{\mathcal C}^\alpha({\mathbb R}^d;{\mathbb R}^d))})\lambda^{-\varepsilon}
\|f_{n,m}\|_{L^q(0,T;{\mathcal C}^{\frac{2}{q}-1}\cap{\mathcal C}^\alpha({\mathbb R}^d))} \nonumber\\&\leq&  C(d,T,\|b\|_{L^q(0,T;{\mathcal C}^{\frac{2}{q}-1}\cap{\mathcal C}^\alpha({\mathbb R}^d;{\mathbb R}^d))})\lambda^{-\varepsilon}
\|f\|_{L^q(0,T;{\mathcal C}^{\frac{2}{q}-1}\cap{\mathcal C}^\alpha({\mathbb R}^d))}
\end{eqnarray}
and
\begin{eqnarray}\label{2.36}
&&\sup_{0\leq t\leq T}\|\nabla u_{n,m}(t)\|_{{\mathcal C}_b^\theta({\mathbb R}^d)} +\|\nabla ^2u_{n,m}\|_{L^2(0,T;{\mathcal C}^\theta_b({\mathbb R}^d))}\nonumber\\&\leq& C(d,T,[b_{n,m}]_{q,\alpha},[b_{n,m}]_{q,\frac{2}{q}-1})\lambda^{-\varepsilon}
([f_{n,m}]_{q,\alpha}+[f_{n,m}]_{q,\frac{2}{q}-1})
\nonumber\\&\leq& C(d,T,[b]_{q,\frac{2}{q}-1},[b]_{q,\alpha})\lambda^{-\varepsilon}
([f]_{q,\frac{2}{q}-1}+[f]_{q,\alpha}).
\end{eqnarray}

By (\ref{2.34})--(\ref{2.36}), then $u_{n,m}\in {\mathcal H}^{2,\theta}_{q,T}$ and there is a positive constant $C$ which depends on $d,T$, $\|b\|_{L^q(0,T;{\mathcal C}^{\alpha}\cap {\mathcal C}^{\frac{2}{q}-1}({\mathbb R}^d;{\mathbb R}^d))}$ and $\lambda$ such that
\begin{eqnarray}\label{2.37}
&&\|(1+|\cdot|^{\frac{2}{q}-1})^{-1}\partial_tu_{n,m}(t,\cdot)\|_{q,0}
\nonumber\\&\leq& \Big[\frac{1}{2}\|\nabla^2 u_{n,m}\|_{2,0}+\|\nabla u_{n,m}\|_{\infty,0}\|(1+|\cdot|^{\frac{2}{q}-1})^{-1}b_{n,m}(t,\cdot)\|_{q,0}\nonumber\\ &&+\lambda
\|(1+|\cdot|^{\frac{2}{q}-1})^{-1}u_{n,m}(t,\cdot)\|_{\infty,0}+ \|(1+|\cdot|^{\frac{2}{q}-1})^{-1}f_{n,m}(t,\cdot)\|_{q,0} \Big]
\nonumber\\ &
\leq& C\Big[\|b\|_{L^q(0,T;{\mathcal C}^{\frac{2}{q}-1}\cap{\mathcal C}^\alpha({\mathbb R}^d;{\mathbb R}^d))}+
\|f\|_{L^q(0,T;{\mathcal C}^{\frac{2}{q}-1}\cap{\mathcal C}^\alpha({\mathbb R}^d))}\Big].
\end{eqnarray}
On account of (\ref{2.35})--(\ref{2.37}), there exists a (unlabelled) subsequence $u_{n,m}$  and a measurable function $u\in  {\mathcal H}^{2,\theta}_{q,T}$ with $\theta\in (0,1+\alpha-2/q)$ such that $u_{n,m}(t,x)\rightarrow u(t,x)$ for a.e. $(t,x)\in [0,T]\times {\mathbb R}^d$ as $m$ and $n$ tend to infinity in turn. Moreover, for every fixed $\theta\in (0,1+\alpha-2/q)$, (\ref{2.35})--(\ref{2.37}) hold for $u$, and in particular (\ref{2.3}) and (\ref{2.4}) are true. Furthermore, since $u_{n,m}$ satisfies (\ref{2.34}), $u$ satisfies (\ref{2.1}).

\medskip
\textbf{Step 3.} For $b\in L^q(0,T;{\mathcal C}^{\frac{2}{q}-1}\cap{\mathcal C}^\alpha({\mathbb R}^d;{\mathbb R}^d))$ and $f\in L^q(0,T;{\mathcal C}^{\frac{2}{q}-1}\cap{\mathcal C}^\alpha({\mathbb R}^d))$, we define $b_R(t,x)=b(t,x\chi_R(x))$ and $f_R(t,x)=f(t,x\chi_R(x))$, where $R>0$, $\chi_R(x)=\chi(x/R)$ and
\begin{eqnarray}\label{2.38}
 \chi \in {\mathcal C}^\infty_0({\mathbb R}),  \ \ 0\leq \chi \leq 1, \ \ \chi^\prime\leq 2 \ \ {\rm and} \  \chi(x)= \left\{\begin{array}{ll}
1, \ \ {\rm if}\ x\in B_R,
\\ 0, \ \ {\rm if}\ x\in {\mathbb R}^d\setminus B_{2R}.
\end{array}\right.
 \end{eqnarray}
Then $b_R\in L^q(0,T;{\mathcal C}^{\alpha}_b({\mathbb R}^d;{\mathbb R}^d))$ and $f_R\in L^q(0,T;{\mathcal C}^{\alpha}_b({\mathbb R}^d))$. Moreover,
\begin{eqnarray}\label{2.39}
 \lim_{R\rightarrow+\infty}[|b_R(t,x)-b(t,x)|+|f_R(t,x)-f(t,x)|]=0, \ \ \forall \ (t,x)\in [0,T]\times{\mathbb R}^d.
\end{eqnarray}
On the other hand, we have
\begin{eqnarray}\label{2.40}
[b_R]_{q,\alpha}&\leq&[b]_{q,\alpha}\sup_{x,y\in{\mathbb R}^d, x\neq y}\frac{|x\chi_R(x)-y\chi_R(y)|^\alpha}{|x-y|^\alpha}\nonumber\\&\leq& [b]_{q,\alpha}\sup_{x,y\in{\mathbb R}^d, x\neq y,\tau\in [0,1]}\Big[\chi_R^\alpha(x)+|\chi_R^\prime(\tau x+(1-\tau)y)|^\alpha\Big]\leq 3[b]_{q,\alpha}
\end{eqnarray}
and
\begin{eqnarray}\label{2.41}
[f_R]_{q,\alpha}\leq [f]_{q,\alpha}\sup_{x,y\in{\mathbb R}^d, x\neq y}\frac{|x\chi_R(x)-y\chi_R(y)|^\alpha}{|x-y|^\alpha}\leq 3[f]_{q,\alpha}.
\end{eqnarray}
Similarly calculations also suggests that
\begin{eqnarray}\label{2.42}
[b_R]_{q,\frac{2}{q}-1}\leq 3[b]_{q,\frac{2}{q}-1} \ \ {\rm and} \ \ [f_R]_{q,\frac{2}{q}-1}\leq 3[f]_{q,\frac{2}{q}-1}.
\end{eqnarray}

By \textbf{Step 2}, there is a unique $u_R\in  {\mathcal H}^{2,\theta}_{q,T}$ with $\theta\in (0,1+\alpha-2/q)$ solving the following Cauchy problem:
\begin{eqnarray}\label{2.43}
\left\{\begin{array}{ll}
\partial_{t}u_R(t,x)=\frac{1}{2}\Delta u_R(t,x)+b_R(t,x)\cdot \nabla u_R(t,x)
\\ \ \qquad\qquad\qquad -\lambda u_R(t,x)+f_R(t,x), \ (t,x)\in (0,T]\times {\mathbb R}^d, \\
u_R(t,x)|_{t=0}=0, \  x\in{\mathbb R}^d.  \end{array}\right.
\end{eqnarray}
Furthermore, it yields that
\begin{eqnarray}\label{2.44}
&&\sup_{0\leq t\leq T}\|(1+|\cdot|^{\frac{2}{q}-1})^{-1} u_R(t,\cdot)\|_0+\|(1+|\cdot|)^{-1}\partial_tu_{n,m}(t,\cdot)\|_{q,0}
\nonumber\\&\leq& C(d,T,\|b_R\|_{L^q(0,T;{\mathcal C}^{\frac{2}{q}-1}\cap{\mathcal C}^\alpha({\mathbb R}^d;{\mathbb R}^d))},\lambda)
\|f_R\|_{L^q(0,T;{\mathcal C}^{\frac{2}{q}-1}\cap{\mathcal C}^\alpha({\mathbb R}^d))} \nonumber\\&\leq&  C(d,T,\|b\|_{L^q(0,T;{\mathcal C}^{\frac{2}{q}-1}\cap{\mathcal C}^\alpha({\mathbb R}^d;{\mathbb R}^d))},\lambda)
\|f\|_{L^q(0,T;{\mathcal C}^{\frac{2}{q}-1}\cap{\mathcal C}^\alpha({\mathbb R}^d))}
\end{eqnarray}
and
\begin{eqnarray}\label{2.45}
&&\sup_{0\leq t\leq T}\|\nabla u_R(t)\|_{{\mathcal C}_b^\theta({\mathbb R}^d)}+\|\nabla ^2u_R\|_{L^2(0,T;{\mathcal C}^\theta_b({\mathbb R}^d))}
\nonumber\\&\leq& C(d,T,\|b_R\|_{L^q(0,T;{\mathcal C}^{\frac{2}{q}-1}\cap{\mathcal C}^\alpha({\mathbb R}^d;{\mathbb R}^d))})\lambda^{-\varepsilon}
\|f_R\|_{L^q(0,T;{\mathcal C}^{\frac{2}{q}-1}\cap{\mathcal C}^\alpha({\mathbb R}^d))} \nonumber\\&\leq&  C(d,T,\|b\|_{L^q(0,T;{\mathcal C}^{\frac{2}{q}-1}\cap{\mathcal C}^\alpha({\mathbb R}^d;{\mathbb R}^d))})\lambda^{-\varepsilon}
\|f\|_{L^q(0,T;{\mathcal C}^{\frac{2}{q}-1}\cap{\mathcal C}^\alpha({\mathbb R}^d))}.
\end{eqnarray}
In view of (\ref{2.44}) and (\ref{2.45}), by letting $R$ tend to infinity in (\ref{2.43}) we get the desired result.

\medskip
Now let us prove the uniqueness. Observing that the equation is linear, it suffices to prove that $u\equiv 0$ if the nonhomogeneous term $f$ vanishes, and it is clear by the estimate (\ref{2.44}).

\medskip
(ii)  Let $u_n$ be the unique strong solution of (\ref{2.1}) with $b$ and $f$ replaced by  $b_n$ and $f_n$ respectively. Then $u_n$ lies in ${\mathcal H}^{2,\theta}_{q,T}$ for every $\theta\in (0, 1+\alpha-2/q)$ and satisfies (\ref{2.3}) and (\ref{2.4}) uniformly in $n$.  It remains to  check $u_n-u\in L^2(0,T;{\mathcal C}_b^{2,\theta}({\mathbb R}^d))\cap L^\infty(0,T;{\mathcal C}_b^{1,\theta}({\mathbb R}^d))$ and (\ref{2.5}).

\medskip
We set $v_n=u_n-u$, then $v_n$ satisfies
\begin{eqnarray}\label{2.46}
\left\{\begin{array}{ll}
\partial_{t}v_n(t,x)=\frac{1}{2}\Delta v_n(t,x)+b_n(t,x)\cdot \nabla v_n(t,x)
\\ \qquad\qquad\quad \ \ -\lambda v_n(t,x)+F_n(t,x), \ (t,x)\in (0,T]\times {\mathbb R}^d, \\
v_n(t,x)|_{t=0}=0, \  x\in{\mathbb R}^d,  \end{array}\right.
\end{eqnarray}
where $F_n(t,x)=f_n(t,x)-f(t,x)+(b_n(t,x)-b(t,x))\cdot\nabla u(t,x)$.

\medskip
Let $x_0\in {\mathbb R}^d$. Consider the following differential equation:
\begin{eqnarray}\label{2.47}
\dot{x}^n_t=-b_n(t,x_0+x^n_t), \quad x^n_t|_{t=0}=0.
\end{eqnarray}
There exists a unique solution to (\ref{2.47}).
By setting $\hat{v}_n(t,x):=v_n(t,x+x_0+x^n_t)$, $\hat{b}_n(t,x):=b_n(t,x+x_0+x^n_t)-b_n(t,x_0+x^n_t)$ and $\hat{F}_n(t,x):=F_n(t,x+x_0+x^n_t)$, then
\begin{eqnarray*}
\left\{\begin{array}{ll}
\partial_{t}\hat{v}_n(t,x)=\frac{1}{2}\Delta \hat{v}_n(t,x)+\hat{b}_n(t,x)\cdot \nabla \hat{v}_n(t,x) \\ \qquad\qquad \ \ -\lambda \hat{v}_n(t,x)+\hat{F}_n(t,x), \ (t,x)\in (0,T]\times {\mathbb R}^d, \\
\hat{v}_n(t,x)|_{t=0}=0, \  x\in{\mathbb R}^d.  \end{array}\right.
\end{eqnarray*}
Thus
\begin{eqnarray}\label{2.48}
\hat{v}_n(t,x) &=&\int_0^te^{-\lambda (t-\tau)}d\tau\int_{{\mathbb R}^d}K(t-\tau,x-y)[\hat{b}_n(\tau,y)\cdot \nabla \hat{v}_n(\tau,y)]dy\nonumber\\ &&+\int_0^te^{-\lambda (t-\tau)}d\tau\int_{{\mathbb R}^d}K(t-\tau,x-y) \hat{F}_n(\tau,y)dy.
\end{eqnarray}
Observing that
\begin{eqnarray}\label{2.49}
|\hat{b}_n(\tau,y)\cdot \nabla \hat{v}_n(\tau,y)|\leq [b_n(\tau)]_{\alpha} |y|^{\alpha} \|\nabla v_n(\tau)\|_0\leq [b(\tau)]_{\alpha} |y|^{\alpha} \sup_{0\leq \tau\leq T}\|\nabla v_n(\tau)\|_0
\end{eqnarray}
and
\begin{eqnarray}\label{2.50}
|\hat{F}_n(\tau,y)|&\leq&| f_n(\tau,y+x_0+x^n_\tau)-f(\tau,y+x_0+x^n_\tau)|\nonumber\\ &&+|(b_n(\tau,y+x_0+x^n_\tau)-b(\tau,y+x_0+x^n_\tau))\cdot\nabla u(\tau,y+x_0+x^n_\tau)|\nonumber\\ &\leq &\int_{{\mathbb R}^d}|f(\tau,y+x_0+x^n_\tau-z)-f(\tau,y+x_0+x^n_\tau)|\rho_n(z)dz\nonumber\\ &&+\int_{{\mathbb R}^d}|b(\tau,y+x_0+x^n_\tau-z)-b(\tau,y+x_0+x^n_\tau)|\rho_n(z)dz\|\nabla u\|_{\infty,0}
\nonumber\\ &\leq & \Big([f(\tau)]_\alpha+[b(\tau)]_\alpha\|\nabla u\|_{\infty,0}\Big)\int_{{\mathbb R}^d}|z|^\alpha\rho_n(z)dz\nonumber\\ &&
\nonumber\\ &\leq & C\Big([f(\tau)]_\alpha+[b(\tau)]_\alpha\Big)\int_{{\mathbb R}^d}|z|^\alpha\rho_n(z)dz,
\end{eqnarray}
we have
\begin{eqnarray}\label{2.51}
\sup_{0\leq t\leq T}\|v_n(t)\|_0
&=&\sup_{0\leq t\leq T}\sup_{x_0\in{\mathbb R}^d}|\hat{v}_n(t,0)|\nonumber\\ &\leq&C(d,T)([b]_{q,\alpha}+[f]_{q,\alpha})( \sup_{0\leq t\leq T}\|\nabla v_n(t)\|_0+\int_{{\mathbb R}^d}|z|^\alpha\rho_n(z)dz),
\end{eqnarray}
which suggests $v_n\in L^\infty(0,T;{\mathcal C}_b({\mathbb R}^d))$, and thus $v_n\in L^2(0,T;{\mathcal C}_b^{2,\theta}({\mathbb R}^d))\cap L^\infty(0,T;{\mathcal C}_b^{1,\theta}({\mathbb R}^d))$. Let us check (\ref{2.5}). By (\ref{2.51}) and the interpolation inequality, it suffices to show
\begin{eqnarray}\label{2.52}
\lim_{n\rightarrow+\infty}\Big\{\sup_{0\leq t\leq T}\|\nabla v_n(t)\|_0 +\|\nabla^2v_n\|_{2,0}\Big\}=0.
\end{eqnarray}
For every fixed $\theta\in (0,1+\alpha-2/q)$, there exists $2/q-1<\alpha^\prime<\alpha$ such that $\theta\in (0,1+\alpha^\prime-2/q)$. For these $\theta$ and $\alpha^\prime$, by analogue calculations of (\ref{2.18}) and (\ref{2.22}), we arrive at
\begin{eqnarray}\label{2.53}
&&\sup_{0\leq t\leq T}\|\nabla v_n(t)\|_0+\|\nabla^2 v_n\|_{2,0}\nonumber \\ &\leq & C(d,T)
[f_n-f+(b_n-b)\cdot\nabla u]_{q,\alpha^\prime-\theta}\nonumber \\ &\leq& C(d,T)
\Big([f_n-f]_{q,\alpha^\prime-\theta}+\|b_n-b\|_{q,0}[\nabla u]_{\infty,\alpha^\prime-\theta}+[b_n-b]_{q,\alpha^\prime-\theta}\|\nabla u\|_{\infty,0}\Big).
\end{eqnarray}
Since $\nabla u\in L^\infty(0,T;{\mathcal C}_b^\iota({\mathbb R}^d;{\mathbb R}^d))$ for every $\iota\in (0, 1+\alpha-2/q)$ and $q>4/(2+\alpha)$, we can choose proper $\alpha^\prime\in (0,\alpha)$ and $\theta\in (0,1+\alpha^\prime-2/q)$ such that $2/q-1<\alpha^\prime-\theta<1+\alpha-2/q$.
For these fixed $\alpha^\prime$ and $\theta$, from (\ref{2.53}) and (\ref{2.3}) we achieve
\begin{eqnarray}\label{2.54}
\sup_{0\leq t\leq T}\|\nabla v_n(t)\|_0+\|\nabla^2 v_n\|_{2,0}\leq C
\Big([f_n-f]_{q,\alpha^\prime-\theta}+\|b_n-b\|_{q,0}+[b_n-b]_{q,\alpha^\prime-\theta}\Big).
\end{eqnarray}
Noting that
\begin{eqnarray}\label{2.55}
[f_n-f]_{q,\alpha}\leq 2[f]_{q,\alpha}, \ \ [b_n-b]_{q,\alpha}\leq 2[b]_{q,\alpha}
\end{eqnarray}
and
\begin{eqnarray}\label{2.56}
\|f_n-f\|_{q,0}\leq [f]_{q,\alpha}\int_{{\mathbb R}^d}|z|^\alpha\rho_n(z)dz, \ \ \|b_n-b\|_{q,0}\leq [b]_{q,\alpha}\int_{{\mathbb R}^d}|z|^\alpha\rho_n(z)dz,
\end{eqnarray}
if one takes $n$ to infinity in (\ref{2.54}), we get (\ref{2.52}). $\Box$

\medskip
We now extend the constant coefficients equation (\ref{2.1}) to a variable coefficients equation and establish an analogue of Lemma \ref{lem2.1}. To be precise, let us consider the following Kolmogorov equation:
\begin{eqnarray}\label{2.57}
\left\{\begin{array}{ll}
\partial_{t}u(t,x)=\frac{1}{2}\sum\limits_{i,j=1}^da_{i,j}(t)\partial^2_{x_i,x_j} u(t,x)+b(t,x)\cdot\nabla u(t,x)\\ \qquad\qquad \ \ -\lambda u(t,x)+f(t,x), \ \ (t,x)\in (0,T]\times {\mathbb R}^d, \\
u(t,x)|_{t=0}=0, \  x\in{\mathbb R}^d,
\end{array}\right.
\end{eqnarray}
where $a_{i,j}(t)$ are Borel bounded measurable functions, which satisfies condition (\ref{1.7}).

\medskip
For $0\leq s<t\leq T$, let
\begin{eqnarray*}
A_{s,t}:=\int_s^ta(\tau)d\tau, \ B_{s,t}=A_{s,t}^{-1}.
\end{eqnarray*}
For every $\vartheta\in{\mathbb R}^d$, it is obvious that
\begin{eqnarray*}
\Theta^{-1}(t-s)|\vartheta|^2\leq \vartheta^\top A_{s,t}\vartheta \leq \Theta(t-s)|\vartheta|^2
\end{eqnarray*}
and
\begin{eqnarray*}
\Theta^{-1}(t-s)^{-1}|\vartheta|^2\leq \vartheta^\top B_{s,t}\vartheta \leq \Theta(t-s)^{-1}|\vartheta|^2.
\end{eqnarray*}Let
\begin{eqnarray}\label{2.58}
\hat{K}(s,t,x)=(2\pi)^{-\frac{d}{2}}\det(B_{s,t})^{\frac{1}{2}}\exp\Big\{-\frac{(B_{s,t}x,x)}{2}\Big\}.
\end{eqnarray}
Then for every $0\leq s<t\leq T$ and $x\in {\mathbb R}^d$, there exist positive constants $C(\Theta)$ and $\mu(\Theta)$ such that
\begin{eqnarray}\label{2.59}
|\nabla^k\hat{K}(s,t,x)|\leq C(t-s)^{-\frac{d+k}{2}}e^{-\frac{\mu |x-y|^2}{t-s}}, \ k=0,1,2.
\end{eqnarray}
By (\ref{2.59}), all calculations used in Lemma \ref{lem2.1} are applicable for the Cauchy problem:
\begin{eqnarray}\label{2.60}
\left\{\begin{array}{ll}
\partial_{t}u(t,x)=\frac{1}{2}\sum\limits_{i,j=1}^da_{i,j}(t)\partial^2_{x_i,x_j} u(t,x) +f(t,x), \ \ (t,x)\in (0,T]\times {\mathbb R}^d, \\
u(t,x)|_{t=0}=0, \  x\in{\mathbb R}^d.
\end{array}\right.
\end{eqnarray}
Using the same arguments for Lemma \ref{lem2.1}, we obtain that
\begin{theorem} \label{the2.2}  Let $b\in L^q(0,T;
{\mathcal C}^{\frac{2}{q}-1}\cap{\mathcal C}^\alpha({\mathbb R}^d;{\mathbb R}^d))$ and
$f\in L^q(0,T;{\mathcal C}^{\frac{2}{q}-1}\cap{\mathcal C}^\alpha({\mathbb R}^d))$ with $\alpha\in (0,1)$ and $q\in (2/(1+\alpha),2)$. Let $(a_{i,j})_{d\times d}$ be a
symmetric $d\times d$ matrix-valued bounded function, which satisfies condition (\ref{1.7}). Then one has:

\medskip
(i) \textbf{(Existence and uniqueness)} There is a unique strong solution $u$ to the Cauchy problem (\ref{2.57}).
Moreover, $u\in {\mathcal H}^{2,\theta}_{q,T}$ for every $\theta\in (0, 1+\alpha-2/q)$ and there is a real number $\varepsilon>0$ such that for large enough $\lambda>0$
\begin{eqnarray}\label{2.61}
\sup_{0\leq t\leq T}\|\nabla u(t)\|_{{\mathcal C}_b^\theta({\mathbb R}^d)}\leq C(d,T,\Theta,[b]_{q,\frac{2}{q}-1},[b]_{q,\alpha})\lambda^{-\varepsilon}
([f]_{q,\frac{2}{q}-1}+[f]_{q,\alpha}).
\end{eqnarray}

(ii) \textbf{(Stability)} Let $u_n$ be the unique strong solution of (\ref{2.57}) with $b$ and $f$ replaced by  $b_n=b\ast \rho_n$ and $f_n=f\ast \rho_n$, respectively.  Then $u_n$ belongs to  ${\mathcal H}^{2,\theta}_{q,T}$  for every $\theta\in (0, 1+\alpha-2/q)$ and satisfies (\ref{2.3})--(\ref{2.4}) uniformly in $n$. Furthermore, $u_n-u\in L^2(0,T;{\mathcal C}_b^{2,\theta}({\mathbb R}^d))$ $\cap L^\infty(0,T;{\mathcal C}_b^{1,\theta}({\mathbb R}^d))$, and if we assume in addition that $q>4/(2+\alpha)$, then for every $\theta\in (0, 1+\alpha-2/q)$ we have
\begin{eqnarray}\label{2.62}
\lim_{n\rightarrow+\infty}\Big[\sup_{0\leq t\leq T}\|u_n(t)-u(t)\|_{{\mathcal C}_b^{1,\theta}({\mathbb R}^d)} +\|\nabla ^2u_n-\nabla^2u\|_{L^2(0,T;{\mathcal C}_b^\theta({\mathbb R}^d))}\Big]=0.
\end{eqnarray}
\end{theorem}

\begin{remark}\label{rem2.3} For the Cauchy problem (\ref{2.57}), when $f\in L^q(0,T;{\mathcal C}_b^{\alpha}({\mathbb R}^d))$ for $q\in (1,+\infty]$ and
$b\in L^\infty(0,T;{\mathcal C}_b^{\alpha}({\mathbb R}^d;{\mathbb R}^d))$, the unique strong solvability of $L^q(0,T;{\mathcal C}_b^{2+\alpha}({\mathbb R}^d))\cap W^{1,q}(0,T;{\mathcal C}_b^\alpha({\mathbb R}^d))$ solutions has been proved by Krylov
\cite{Kry02}, and when $f$ has polynomially or exponentially growth with H\"{o}lder norms, the well-posedness was also established by Lorenzi \cite{Lor}. Recently, Tian, Ding and Wei \cite{TDW} generalized Krylov and Lorenzi's results to the case of $b\in L^2(0,T;{\mathcal C}_b^\alpha({\mathbb R}^d;{\mathbb R}^d))$. When the coefficients are bounded in time variable and locally $\alpha$-H\"{o}lder continuous in space variable (uniformly in time), the $L^\infty({\mathcal C}^{2+\alpha})$-Schauder estimate for solutions was derived by Krylov and and Priola \cite{KP} as well for the following parabolic PDE (also see \cite{CMP} for fractional PDE):
\begin{eqnarray}\label{2.63}
\left\{\begin{array}{ll}
\partial_{t}u(t,x)+\sum_{i,j=1}^da_{i,j}(t,x)\partial^2_{x_i,x_j} u(t,x)+b(t,x)\cdot\nabla u(t,x)\\ \qquad\qquad -c(t,x)u(t,x)=f(t,x), \ \ (t,x)\in (T,S)\times {\mathbb R}^d, \\ u(t,x)|_{t=S}=g(x), \quad |f(t,x)|\leq F_0c(t,x).
\end{array}\right.
\end{eqnarray}
Here, we only assume that
$b\in L^q(0,T;{\mathcal C}^{\frac{2}{q}-1}\cap{\mathcal C}^\alpha({\mathbb R}^d;{\mathbb R}^d))$
with $q\in (2/(1+\alpha),2)$,
so we extend the existing results. This result plays a key role in proving the strong well-posedness for SDE (\ref{1.1}).
\end{remark}

\section{Proof of Theorem \ref{the1.3}}
\label{sec3}\setcounter{equation}{0}
Let $\lambda>0$ be a large enough real number. Consider the following vector-valued Cauchy problem:
\begin{eqnarray}\label{3.1}
\left\{\begin{array}{ll}\partial_tU(t,x)+\frac{1}{2}\sum\limits_{i,j=1}^da_{i,j}(t)\partial^2_{x_i,x_j} U(t,x)+b(t,x)\cdot\nabla U(t,x)\\ \qquad\qquad=\lambda U(t,x)-b(t,x), \ \ (t,x)\in [0,T)\times {\mathbb R}^d,\\
U(t,x)|_{t=T}=0, \  x\in{\mathbb R}^d,
\end{array}\right.
\end{eqnarray}
where $(a_{i,j})_{d\times d}=\sigma\sigma^{\top}$. Since
$b\in L^q(0,T;{\mathcal C}^{\frac{2}{q}-1}\cap{\mathcal C}^\alpha({\mathbb R}^d;{\mathbb R}^d))$ and
$\sigma\in L^\infty(0,T)$, in view of Theorem \ref{the2.2}, there is a unique
$U\in ({\mathcal H}^{2,\theta}_{q,T})^d$ (see (\ref{2.2}))
solving the Cauchy problem (\ref{3.1}) for every $\theta\in (0,1+\alpha-2/q)$. Moreover, by (\ref{2.61}) there is a real number $\varepsilon>0$ such that for large enough $\lambda>0$,
\begin{eqnarray}\label{3.2}
\sup_{0\leq t\leq T}\|\nabla U(t)\|_{{\mathcal C}_b^\theta({\mathbb R}^d)}\leq C(d,T,\Theta,[b]_{q,\frac{2}{q}-1},[b]_{q,\alpha})\lambda^{-\varepsilon}
([b]_{q,\frac{2}{q}-1}+[b]_{q,\alpha})<\frac{1}{2}.
\end{eqnarray}
We set $\Phi(t,x)=x+U(t,x)$, by (\ref{3.2}) then $\Phi$ forms a non-singular diffeomorphism of class ${\mathcal C}^{1,\theta}$
uniformly in $t\in [0,T]$ and
\begin{eqnarray}\label{3.3}
\frac{1}{2}<\sup_{0\leq t\leq T}\|\nabla\Phi(t)\|_0 <\frac{3}{2},
\quad  \frac{2}{3}<\sup_{0\leq t\leq T}\|\nabla\Psi(t)\|_0<2,
\end{eqnarray}
where $\Psi(t,\cdot)=\Phi^{-1}(t,\cdot)$.

\medskip
For $0<\epsilon<1$ and $t\in [0,T]$, define
\begin{eqnarray*}
U_\epsilon(t,x)=\frac{1}{\epsilon}\int_t^{t+\epsilon}U(\tau,x)d\tau=\int_0^1U(t+\tau\epsilon,x)d\tau
\end{eqnarray*}
and $\Phi_\epsilon(t,x)=x+U_\epsilon(t,x)$, where $U(\tau,x):=U(T,x)=0$ when $\tau>T$. Notice that $\Phi_\epsilon\in {\mathcal C}^1([0,T];{\mathcal C}({\mathbb R}^d))\cap {\mathcal C}([0,T];{\mathcal C}^2({\mathbb R}^d))$, if $X_{s,t}$ is a strong solution of SDE (\ref{1.1}), in light of It\^{o}'s formula, we derive
\begin{eqnarray}\label{3.4}
\Phi_\epsilon(t,X_{s,t}(x))&=&\Phi_\epsilon(s,x)+\int_s^t\partial_\tau U_\epsilon(\tau,X_{s,\tau}(x))d\tau+\int_s^tb(\tau,X_{s,\tau}(x))\cdot\nabla U_\epsilon(\tau,X_{s,\tau}(x))d\tau
\nonumber\\&&+\frac{1}{2}\int_s^t tr(\nabla^2 U_\epsilon(\tau,X_{s,\tau}(x))\sigma(\tau)\sigma^{\top}(\tau))d\tau+\int_s^tb(\tau,X_{s,\tau}(x))d\tau
\nonumber\\&&+ \int_s^t[I+\nabla U_\epsilon(\tau,X_{s,\tau}(x))]\sigma(\tau)dW_\tau.
\end{eqnarray}
Since $U\in ({\mathcal H}^{2,\theta}_{q,T})^d$, if one lets $\epsilon$ tend to $0$, we obtain
\begin{eqnarray}\label{3.5}
\left\{\begin{array}{ll}
U_\epsilon(t,X_{s,t}(x)) \longrightarrow U(t,X_{s,t}(x)), & \forall \ t\in [s,T], \  {\mathbb P}-a.s.,
\\
\partial_\tau U_\epsilon(\tau,X_{s,\tau}(x))\longrightarrow \partial_\tau U(\tau,X_{s,\tau}(x)), & a.e. \ \tau\in [s,T], \  {\mathbb P}-a.s.,
\\
\nabla U_\epsilon(\tau,X_{s,\tau}(x))\longrightarrow \nabla U(\tau,X_{s,\tau}(x)), & a.e. \ \tau\in [s,T], \ {\mathbb P}-a.s.,
\\
\nabla^2U_\epsilon(\tau,X_{s,\tau}(x))\longrightarrow \nabla^2U_\epsilon(\tau,X_{s,\tau}(x)), & a.e. \ \tau\in [s,T], \ {\mathbb P}-a.s..
\end{array}\right.
\end{eqnarray}
Combining (\ref{3.4}) and (\ref{3.5}), by employing Lebesgue's dominated convergence theorem, then
\begin{eqnarray}\label{3.6}
\Phi(t,X_{s,t}(x))&=&\Phi(s,x)+\int_s^t\partial_\tau U(\tau,X_{s,\tau}(x))d\tau+\int_s^tb(\tau,X_{s,\tau}(x))\cdot\nabla U(\tau,X_{s,\tau}(x))d\tau
\nonumber\\&&+
\frac{1}{2}\int_s^t tr(\nabla^2 U(\tau,X_{s,\tau}(x))\sigma(\tau)\sigma^{\top}(\tau))d\tau
+\int_s^tb(\tau,X_{s,\tau}(x))d\tau
\nonumber\\&&+\int_s^t[I+\nabla U(\tau,X_{s,\tau}(x))]\sigma(\tau)dW_\tau
\nonumber\\ &=&\Phi(s,x)+\lambda \int_s^tU(\tau,X_{s,\tau}(x))d\tau+\int_s^t(I+\nabla U(\tau,X_{s,\tau}(x)))\sigma(\tau)dW_\tau,
\end{eqnarray}
where in the second identity we have used the fact that $U(t,x)$ satisfies the Cauchy problem (\ref{3.1}).

\medskip
Denote $Y_{s,t}=\Phi(t,X_{s,t})$,  it follows from (\ref{3.6}) that
\begin{eqnarray}\label{3.7}
\left\{
  \begin{array}{ll}
  dY_{s,t}=\lambda U(t,\Psi(t,Y_{s,t}))dt+(I+\nabla U(t,\Psi(t,Y_{s,t})))\sigma(t)dW_t\\ \quad\quad \ =:
\tilde{b}(t,Y_{s,t})dt+\tilde{\sigma}(t,Y_{s,t})dW_t,\ t\in(s,T], \\
Y_{s,t}|_{t=s}=y=\Phi(s,x).
  \end{array}
\right.
\end{eqnarray}
Conversely, if $Y_{s,t}$ is a strong solution of SDE (\ref{3.7}), then $X_{s,t}=\Psi(t,Y_{s,t})$ satisfies SDE (\ref{1.1}). Therefore SDEs (\ref{1.1}) and (\ref{3.7}) are equivalent. Observing that $\partial_tU \in L^q(0,T;{\mathcal C}({\mathbb R}^d;{\mathbb R}^d))$, one has  $U\in {\mathcal C}^{1-\frac{1}{q}}([0,T];{\mathcal C}({\mathbb R}^d;{\mathbb R}^d))$ and $\Phi\in {\mathcal C}^{1-\frac{1}{q}}([0,T];{\mathcal C}({\mathbb R}^d;{\mathbb R}^d))$  by using the Sobolev imbedding theorem. This, together with the fact that $\Phi(t,\cdot)$ forms a non-singular diffeomorphism of class ${\mathcal C}^{1,\theta}$ uniformly in $t$, implies that we only need to prove the conclusions (i) and (\ref{1.8}) for $Y_{s,t}$ and $Y_{s,t}^{-1}$.  For conclusions (\ref{1.9}) and (iii), we first prove for $Y_{s,t}$ and $Y_{s,t}^{-1}$ then to prove for $X_{s,t}$ and $X_{s,t}^{-1}$. Since the calculations from $Y$ to $X$ are similar, we only prove  the conclusion (iii) for $Y_{s,t}$ and $Y_{s,t}^{-1}$, and for (\ref{1.9}) we give the complete proof details.

\medskip
(i) We divide the proof of stochastic flow of homeomorphisms into two parts.

\medskip
$\bullet$ \textbf{The unique strong solvability.} By the regularity of $U$ and assumptions on $\sigma$, we have $\tilde{b}\in L^\infty(0,T;{\mathcal C}^{\frac{2}{q}-1}({\mathbb R}^d;{\mathbb R}^d))$, $\nabla\tilde{b}\in L^\infty(0,T;{\mathcal C}_b^\theta({\mathbb R}^d;{\mathbb R}^{d\times d}))$ and $\tilde{\sigma}\in L^2(0,T;{\mathcal C}_b^{1,\theta}({\mathbb R}^d;{\mathbb R}^{d\times d}))\cap L^\infty(0,T;{\mathcal C}_b^\theta({\mathbb R}^d;{\mathbb R}^{d\times d}))$. Owing to Cauchy-Lipschitz's theorem, there exists a unique strong solution $Y_{s,t}(y)$ to (\ref{3.7}). Moreover, an application of  the It\^{o} formula to $|Y_{s,t}|^p$ yields that
\begin{eqnarray*}
d|Y_{s,t}(y)|^p&\leq& p|Y_{s,t}(y)|^{p-1} |\tilde{b}(t,Y_{s,t}(y))|dt+\frac{p(p-1)}{2}|Y_{s,t}(y)|^{p-2} \|\tilde{\sigma}(t,Y_{s,t}(y))\|^2dt \nonumber\\ &&+p|Y_{s,t}(y)|^{p-2}\langle Y_{s,t}(y), \tilde{\sigma}(t,Y_{s,t}(y))dW_t \rangle
\nonumber\\ &\leq& C[1+|Y_{s,t}(y)|^p]dt+p|Y_{s,t}(y)|^{p-2}\langle Y_{s,t}(y), \tilde{\sigma}(t,Y_{s,t}(y))dW_t\rangle.
\end{eqnarray*}
Observing that for every $t>s$, $\int_s^t|Y_{s,\tau}(y)|^{p-2}\langle Y_{s,\tau}(y), \tilde{\sigma}(\tau,Y_{s,\tau}(y))dW_\tau\rangle$ is a martingale. Then with the help of Gr\"{o}nwall's inequality, we get
\begin{eqnarray}\label{3.8}
\sup_{s\leq t\leq T}{\mathbb E}|Y_{s,t}(y)|^p\leq C(1+|y|^p).
\end{eqnarray}

\medskip
$\bullet$ \textbf{Stochastic flow of homeomorphisms.} Due to \cite[Lemmas II.2.4, II.4.1 and II.4.2]{Kun84}, we should prove that: for every $x,x^\prime,y,y^\prime\in {\mathbb R}^d$ ($x\neq y,\, x^\prime\neq y^\prime$) and every $s,t,s^\prime,t^\prime\in [0,T]$ ($s<t, \, s^\prime<t^\prime$),
\begin{eqnarray}\label{3.9}
\sup_{s\leq t\leq T}{\mathbb E}|Y_{s,t}(x)-Y_{s,t}(y)|^{2\varsigma}\leq C|x-y|^{2\varsigma},\quad \forall  \ \varsigma<0,
\end{eqnarray}
and
\begin{eqnarray}\label{3.10}
&&{\mathbb E}|\eta_{s,t}(x,x^\prime)-\eta_{s^\prime,t^\prime}(y,y^\prime)|^p
\nonumber \\ &\leq& C\delta^{-2p}\Big\{(1+|x|^p+|x^\prime|^p+|y|^p+|y^\prime|^p)[|s-s^\prime|^{\frac{p}{2}}+
|t-t^\prime|^{\frac{p}{2}}]\nonumber \\ && +|x-y|^p+|x^\prime-y^\prime|^p\Big\}, \ \forall \ p\geq 2, \ \ |x-x^\prime|\geq \delta>0, \ \ |y-y^\prime|\geq \delta>0,
\end{eqnarray}
and
\begin{eqnarray}\label{3.11}
{\mathbb E}|\eta_{s,t}(\hat{x})-\eta_{s^\prime,t^\prime}(\hat{y})|^p\leq C[|\hat{x}-\hat{y}|^p+|s-s^\prime|^{\frac{p}{2}}+|t-t^\prime|^{\frac{p}{2}}],\ \ \forall \ p>0,
\end{eqnarray}
where
\begin{eqnarray*}
\eta_{s,t}(x,x^\prime)=\frac{1}{|Y_{s,t}(x)-Y_{s,t}(x^\prime)|} \ \ {\rm and}  \ \  \eta_{s,t}(\hat{x})=\left\{
  \begin{array}{ll}
\frac{1}{1+|Y_{s,t}(x)|},& {\rm if} \ \hat{x}=|x|^{-2}x\in {\mathbb R}^d,\\
\quad 0,  & {\rm if} \ \hat{x}=|x|^{-2}x=\infty.
  \end{array}
\right.
\end{eqnarray*}
For every $p\geq 2$, we have
\begin{eqnarray*}
&&|\eta_{s,t}(x,x^\prime)-\eta_{s^\prime,t^\prime}(y,y^\prime)|^p\\ \nonumber &\leq& 2^{p-1}|\eta_{s,t}(x,y)|^p|\eta_{s^\prime,t^\prime}(x^\prime,y^\prime)|^p[|Y_{s,t}(x)-Y_{s^\prime,t^\prime}(y)|^p
+|Y_{s,t}(x^\prime)-Y_{s^\prime,t^\prime}(y^\prime)|^p]
\end{eqnarray*}
and
\begin{eqnarray*}
|\eta_{s,t}(\hat{x})-\eta_{s^\prime,t^\prime}(\hat{y})|^p\leq |\eta_{s,t}(\hat{x})|^p|\eta_{s^\prime,t^\prime}(\hat{y})|^p|Y_{s,t}(x)-Y_{s^\prime,t^\prime}(y)|^p.
\end{eqnarray*}
These, together with  H\"{o}lder's inequality, the element inequality $|a_1+b_1|^{\tilde{p}}\leq \max\{2^{\tilde{p}-1}, 1\}[a_1^{\tilde{p}}+b_1^{\tilde{p}}]$ $(a_1,b_1,{\tilde{p}}\in {\mathbb R}_+$), (\ref{3.9}), and the estimate
\begin{eqnarray*}
{\mathbb E}\sup_{s\leq t\leq T}(1+|Y_{s,t}(y)|)^\varsigma\leq C(1+|y|)^\varsigma, \quad \forall \ \varsigma<0,
\end{eqnarray*}
which can be proved clearly since $\tilde{b}$ is Lipschitz continuous and $\tilde{\sigma}$ is bounded, implies the estimates (\ref{3.10}) and (\ref{3.11}) if one shows the following inequality
\begin{eqnarray}\label{3.12}
{\mathbb E}|Y_{s,t}(x)-Y_{s^\prime,t^\prime}(y)|^p
\leq \Big\{|x-y|^p+(1+|x|^p+|y|^p)[|s-s^\prime|^{\frac{p}{2}}+|t-t^\prime|^{\frac{p}{2}}]\Big\},\ \forall \ p\geq 2.
\end{eqnarray}
Thus, it suffices to show (\ref{3.9}) and (\ref{3.12}).

\medskip
Let us prove (\ref{3.9}) first. For $\epsilon>0$, we
choose $F_\epsilon(x)=f^\varsigma_\epsilon(x)=(\epsilon+|x|^2)^\varsigma$ and
$Y_{s,t}(x,y):=Y_{s,t}(x)-Y_{s,t}(y)$. Thanks to the It\^{o} formula, then
\begin{eqnarray}\label{3.13}
&&F_\epsilon(Y_{s,t}(x,y))\nonumber\\
&=&F_\epsilon(x-y)+2\varsigma\int_s^tf^{\varsigma-1}_\epsilon(Y_{s,\tau}(x,y))\langle Y_{s,\tau}(x,y), \tilde{b}(\tau,Y_{s,\tau}(x))-\tilde{b}(\tau,Y_{s,\tau}(y))\rangle d\tau\nonumber\\&&+
2\varsigma\int_s^tf^{\varsigma-1}_\epsilon(Y_{s,\tau}(x,y))\langle Y_{s,\tau}(x,y), (\tilde{\sigma}(\tau,Y_{s,\tau}(x))-\tilde{\sigma}(\tau,Y_{s,\tau}(y)))dW_\tau\rangle
\nonumber\\&&+\varsigma(\varsigma-1)\sum_{i,j,k=1}^d\int_s^tf^{\varsigma-2}_\epsilon(Y_{s,\tau}(x,y))
[f_\epsilon(Y_{s,\tau}(x,y))\delta_{i,j}+
2Y^i_{s,\tau}(x,y)Y^j_{s,\tau}(x,y)]
\nonumber\\&&\quad\times[\tilde{\sigma}_{i,k}(\tau,Y_{s,\tau}(x))-
\tilde{\sigma}_{i,k}(\tau,Y_{s,\tau}(y))][\tilde{\sigma}_{j,k}(\tau,Y_{s,\tau}(x))
-\tilde{\sigma}_{j,k}(\tau,Y_{s,\tau}(y))]d\tau
\nonumber\\&\leq& F_\epsilon(x-y)+C|\varsigma|\int_s^tF_\epsilon(Y_{s,\tau}(x,y))d\tau+C|\varsigma(\varsigma-1)
|\int_s^t\kappa^2(\tau)
F_\epsilon(Y_{s,\tau}(x,y))d\tau\nonumber\\&&+
2\varsigma\int_s^tf^{\tau-1}_\epsilon(Y_{s,\tau}(x,y))\langle Y_{s,\tau}(x,y), (\tilde{\sigma}(\tau,Y_{s,\tau}(x))-\tilde{\sigma}(\tau,Y_{s,\tau}(y)))dW_\tau\rangle,
\end{eqnarray}
where $\kappa(\tau)=\|\nabla^2 U(\tau)\|_0\in L^2(0,T)$ for $\nabla U\in L^2(0,T;{\mathcal C}_b^{1,\theta}({\mathbb R}^d;{\mathbb R}^{d\times d}))$ and
\begin{eqnarray*}
Y_{s,\tau}=(Y^1_{s,\tau},Y^2_{s,\tau},\ldots,Y^d_{s,\tau}), \ \ \delta_{i,j}=\left\{\begin{array}{ll}
1, \ \ \mbox{if} \ \ i=j, \\
0, \  \ \mbox{if} \ \ i\neq j.
\end{array}\right.
\end{eqnarray*}
By (\ref{3.13}) and the Gr\"{o}nwall inequality, we derive
\begin{eqnarray*}
\sup_{s\leq t\leq T}{\mathbb E}[\epsilon+|Y_{s,t}(x)-Y_{s,t}(y)|^2]^\varsigma\leq C[\epsilon+|x-y|^2]^\varsigma.
\end{eqnarray*}
By letting $\epsilon\downarrow 0$, then (\ref{3.9}) holds.

\medskip
Let $x,y,\in {\mathbb R}^d$ ($x\neq y$) and $s,t,s^\prime,t^\prime\in [0,T]$ ($s<t, \, s^\prime<t^\prime$). Without loss of generality, we assume $s<s^\prime<t<t^\prime$. For every $p\geq 2$, then
\begin{eqnarray}\label{3.14}
&&|Y_{s,t}(x)-Y_{s^\prime,t^\prime}(y)|^p\nonumber\\&\leq& 3^{p-1}[|Y_{s,t}(x)-Y_{s,t}(y)|^p+|Y_{s,t}(y)-Y_{s^\prime,t}(y)|^p+
|Y_{s^\prime,t}(y)-Y_{s^\prime,t^\prime}(y)|^p].
\end{eqnarray}
By using the It\^{o} formula to $|Y_{s,t}(x)-Y_{s,t}(y)|^p$ to get
\begin{eqnarray}\label{3.15}
&&{\mathbb E}|Y_{s,t}(x)-Y_{s,t}(y)|^{p}\nonumber\\&\leq& |x-y|^p+p(p-1) \Big\{{\mathbb E}\int_s^t|Y_{s,\tau}(x)-Y_{s,\tau}(y)|^{p-1}|\tilde{b}(\tau,Y_{s,\tau}(x))
-\tilde{b}(\tau,Y_{s,\tau}(y))|d\tau
\nonumber\\&& +
{\mathbb E}\int_s^t|Y_{s,\tau}(x)-Y_{s,\tau}(y)|^{p-2}
\|\tilde{\sigma}(\tau,Y_{s,\tau}(x))-\tilde{\sigma}(\tau,Y_{s,\tau}(y))\|^2d\tau\Big\}
\nonumber\\&\leq& |x-y|^p+C(d,T,\Theta,p) \int_s^t[1+\kappa^2(\tau)]{\mathbb E}|Y_{s,\tau}(x)-Y_{s,\tau}(y)|^{p}d\tau,
\end{eqnarray}
which yields
\begin{eqnarray}\label{3.16}
\sup_{s\leq t\leq T}{\mathbb E}|Y_{s,t}(x)-Y_{s,t}(y)|^p\leq C|x-y|^p.
\end{eqnarray}
Similarly, by It\^{o}'s formula and BDG's inequality, for every $p\geq 2$, then
\begin{eqnarray}\label{3.17}
&&\sup_{0\leq s\leq T}{\mathbb E}\sup_{s\leq t \leq T}|Y_{s,t}(x)-Y_{s,t}(y)|^p
\nonumber\\&\leq& |x-y|^p+C\sup_{0\leq s\leq T}{\mathbb E}\int_s^T[1+\kappa^2(\tau)]|Y_{s,\tau}(x)-Y_{s,\tau}(y)|^pd\tau\nonumber\\&&+C\Big[\sup_{0\leq s\leq T}{\mathbb E}\int_s^T\kappa^2(\tau)|Y_{s,\tau}(x)-Y_{s,\tau}(y)|^{2p}d\tau\Big]^{\frac{1}{2}}\leq C|x-y|^p,
\end{eqnarray}
where in the last inequality we have used (\ref{3.16}).

\medskip
For $|Y_{s,t}(y)-Y_{s^\prime,t}(y)|^p$, by employing the It\^{o} formula again, one ascertains
\begin{eqnarray*}
&&{\mathbb E}|Y_{s,t}(y)-Y_{s^\prime,t}(y)|^p\nonumber\\&\leq&{\mathbb E}|Y_{s,s^\prime}(y)-y|^p
+C{\mathbb E}\int_{s^\prime}^t|Y_{s,\tau}(y)-Y_{s^\prime,\tau}(y)|^{p-1}|\tilde{b}(\tau,Y_{s,\tau}(y))
-\tilde{b}(\tau,Y_{s,\tau}(x))|d\tau
\nonumber\\&&+C{\mathbb E}\int_{s^\prime}^t|Y_{s,\tau}(y)-Y_{s^\prime,\tau}(y)|^{p-2}
\|\tilde{\sigma}(\tau,Y_{s,\tau}(y))-\tilde{\sigma}(\tau,Y_{s,\tau}(y))\|^2d\tau
\nonumber\\&\leq&{\mathbb E}|Y_{s,s^\prime}(y)-y|^p
+C{\mathbb E}\int_{s^\prime}^t[1+\kappa^2(\tau)]|Y_{s,\tau}(y)-Y_{s^\prime,\tau}(y)|^pd\tau,
\end{eqnarray*}
where $\kappa$ is given in (\ref{3.13}). This, together with the Gr\"{o}nwall, Minkowski and BGD inequalities, leads to
\begin{eqnarray}\label{3.18}
{\mathbb E}|Y_{s,t}(y)-Y_{s^\prime,t}(y)|^p&\leq&C{\mathbb E}|Y_{s,s^\prime}(y)-y|^p
\nonumber \\ &=&C{\mathbb E}\Big|\int_s^{s^\prime}\tilde{b}(\tau,Y_{s,\tau}(y))d\tau+
\int_s^{s^\prime}\tilde{\sigma}(\tau,Y_{s,\tau}(y))dW_\tau\Big|^p\nonumber \\
&\leq&C\Big|\int_s^{s^\prime}[{\mathbb E}|\tilde{b}(\tau,Y_{s,\tau}(y))|^p]^{\frac{1}{p}}d\tau\Big|^p+C{\mathbb E}\Big[
\int_s^{s^\prime}\|\tilde{\sigma}(\tau,Y_{s,\tau}(y))\|^2d\tau\Big]^{\frac{p}{2}}
\nonumber \\  &\leq& C[1+\sup_{s\leq \tau\leq T}{\mathbb E}|Y_{s,\tau}(y))|^p]|s-s^\prime|^p+C|s-s^\prime|^{\frac{p}{2}}
\nonumber \\ &\leq& C[(1+|y|^p)|s-s^\prime|^p+|s-s^\prime|^{\frac{p}{2}}]\leq C(1+|y|^p)|s-s^\prime|^{\frac{p}{2}},
\end{eqnarray}
where in the fourth line we have used the fact $\tilde{b}$ is Lipschitz continuous uniformly in time variable and $\tilde{\sigma}$ is bounded,  and in the fifth line we have used (\ref{3.8}).

\medskip
For the term $|Y_{s^\prime,t}(y)-Y_{s^\prime,t^\prime}(y)|^p$, then
\begin{eqnarray}\label{3.19}
{\mathbb E}|Y_{s^\prime,t}(y)-Y_{s^\prime,t^\prime}(y)|^p&=& {\mathbb E}\Big|\int_t^{t^\prime}\tilde{b}(\tau,Y_{s^\prime,\tau}(y))d\tau +\int_t^{t^\prime}\tilde{\sigma}(\tau,Y_{s^\prime,\tau}(y))dW_\tau\Big|^p
\nonumber\\&\leq& C(1+|y|^p)|t-t^\prime|^{\frac{p}{2}}.
\end{eqnarray}
Summing over (\ref{3.14}), (\ref{3.16}) and  (\ref{3.18}) and (\ref{3.19}), we obtain (\ref{3.12}). Thus $Y_{s,t}(\cdot)$ forms a homeomorphism. Observing that $Y_{s,t}$ satisfies equation (\ref{3.7}), then
\begin{eqnarray*}
Y_{s,t}(Y^{-1}_{s,t}(y))=Y^{-1}_{s,t}(y)+\int_s^t\tilde{b}(\tau,Y_{s,\tau}(Y^{-1}_{s,t}(y)))d\tau+
\int_s^t\tilde{\sigma}(\tau,Y_{s,\tau}(Y^{-1}_{s,t}(y)))dW_\tau.
\end{eqnarray*}
Noting that $Y_{s,\tau}(Y^{-1}_{s,t}(y))=Y^{-1}_{\tau,t}(y)$, thus
\begin{eqnarray}\label{3.20}
Y^{-1}_{s,t}(y)=y-\int_s^t\tilde{b}(\tau,Y^{-1}_{\tau,t}(y))d\tau-
\int_s^t\tilde{\sigma}(\tau,Y^{-1}_{\tau,t}(y))dW_\tau.
\end{eqnarray}
Hence $Y^{-1}_{s,t}(y)$ is continuous in $(s,t,y)$, almost surely in $\omega$, and $\{Y_{s,t}(x), \ t\in [s,T ]\}$ forms a stochastic flow of homeomorphisms to SDE
(\ref{3.7}).

\medskip
(ii) We now turn to show the  weak differentiability, gradient and H\"{o}lder estimates. Observing that the inverse flow $Y_{s,t}^{-1}$ satisfies SDE (\ref{3.20}), which has the same form as the original one (\ref{3.7}) (only the drift and diffusion have opposite sign), the proof of the weak differentiability, gradient and H\"{o}lder estimates for $Y_{s,t}^{-1}$ are similar to that of  $Y_{s,t}$ after taking into consideration the backward character of the equation. For the weak differentiability and gradient estimates, it suffices to show the conclusions for $Y_{s,t}$. For the H\"{o}lder estimate, we first prove the conclusion for $Y_{s,t}$, then to prove the estimates for $X_{s,t}$.

\medskip
$\bullet$ \textbf{Weak differentiability and gradient estimates.} Let $e\in {\mathbb R}^d$ with $|e|=1$. If
\begin{eqnarray}\label{3.21}
\lim_{\delta\rightarrow 0}\frac{Y_{s,t}(y+\delta e)-Y_{s,t}(y)}{\delta}
\end{eqnarray}
exists in $L^2(\Omega)$ uniformly  in $s$ and $t$, and the limit is continuous in $y$, then we complete the proof. The continuity of the limit for space variable will be proved in the next step, we only show (\ref{3.21}) and the moment estimates for the weak derivatives.

\medskip
Let $Y_{s,t}$ be the unique strong solution of SDE (\ref{3.7}). Consider the following SDE:
\begin{eqnarray}\label{3.22}
d \zeta_{s,t}(y)&=&\lambda \nabla U(t,\Psi(t,Y_{s,t}(y)))\nabla\Psi(t,Y_{s,t}(y))\zeta_{s,t}(y)dt\nonumber\\&&+\nabla^2
U(t,\Psi(t,Y_{s,t}(y)))\nabla\Psi(t,Y_{s,t}(y))\sigma(t)\zeta_{s,t}(y)dW_t\nonumber\\
&=&:\hat{b}(t,Y_{s,t}(y))\zeta_{s,t}(y)dt+\hat{\sigma}(t,Y_{s,t}(y))\zeta_{s,t}(y)dW_t, \quad t\in (s,T],
\end{eqnarray}
with $\zeta_{s,t}(y)|_{t=s}=e$. Since the equation is linear, $\hat{b}\in L^\infty(0,T;{\mathcal C}_b^\theta({\mathbb R}^d))$ and $\hat{\sigma}\in L^2(0,T;{\mathcal C}_b^\theta({\mathbb R}^d))$ for every $\theta\in (0,1+\alpha-2/q)$, there exists a unique strong solution $\zeta_{s,t}(y)$ of (\ref{3.22}). Moreover, for every $p\geq 2$ there exists a positive constant $C(d,T,\Theta,p,\|b\|_{L^q(0,T;{\mathcal C}^{\frac{2}{q}-1}\cap{\mathcal C}^{\alpha}({\mathbb R}^d;{\mathbb R}^d))})$ such that
\begin{eqnarray}\label{3.23}
\sup_{y\in{\mathbb R}^d}\sup_{0\leq s\leq
T}{\mathbb E}\sup_{s\leq t \leq T}\|\zeta_{s,t}(y)\|^p\leq C.
\end{eqnarray}

\medskip
For $\delta\in {\mathbb R}$, set $Y_{s,t}^\delta(y):=Y_{s,t}(y+\delta e)-Y_{s,t}(y)$, then by (\ref{3.7})
\begin{eqnarray}\label{3.24}
Y_{s,t}^\delta(y)&=&\delta e+\int_s^t[\tilde{b}(\tau,Y_{s,\tau}(y+\delta e))-\tilde{b}(\tau,Y_{s,\tau}(y))]d\tau\nonumber\\&&+\int_s^t[\tilde{\sigma}(\tau,Y_{s,\tau}(y+\delta e))-\tilde{\sigma}(\tau,Y_{s,\tau}(y))]dW_\tau\nonumber\\&=&\delta e+\int_0^1\int_s^t\nabla\tilde{b}
(\tau,\iota Y_{s,\tau}(y+\delta e)+(1-\iota)Y_{s,\tau}(y))Y_{s,\tau}^\delta(y)d\iota d\tau
\nonumber\\&&+\int_0^1\int_s^t\nabla\tilde{\sigma}(\tau,\iota Y_{s,\tau}(y+\delta e)+(1-\iota)Y_{s,\tau}(y))Y_{s,\tau}^\delta(y)d\iota dW_\tau .
\end{eqnarray}
By virtue of It\^{o}'s formula and H\"{o}lder's inequality, we achieve from (\ref{3.24}) that
\begin{eqnarray*}
{\mathbb E}|Y_{s,t}^\delta(y)|^p\leq |\delta|^p+C{\mathbb E}\int_s^t[1+\kappa^2(\tau)]|Y_{s,\tau}^\delta(y)|^pd\tau,
\end{eqnarray*}
for $\nabla \tilde{b}\in L^\infty(0,T;{\mathcal C}_b^\theta({\mathbb R}^d;{\mathbb R}^{d\times d}))$ and $\tilde{\sigma}\in L^2(0,T;{\mathcal C}_b^{1,\theta}({\mathbb R}^d;{\mathbb R}^{d\times d}))$, which suggests that
\begin{eqnarray}\label{3.25}
\sup_{s\leq t\leq T}{\mathbb E}\Big|\frac{|Y_{s,t}^\delta(y)|}{\delta}\Big|^p\leq C.
\end{eqnarray}
We rewrite (\ref{3.24}) by
\begin{eqnarray}\label{3.26}
d\frac{Y_{s,t}^\delta(y)}{\delta}&=&\int_0^1\nabla\tilde{b}
(t,\iota Y_{s,t}(y+\delta e)+(1-\iota)Y_{s,t}(y))d\iota\frac{Y_{s,t}^\delta(y)}{\delta}dt
\nonumber\\&&+\int_0^1\nabla\tilde{\sigma}(t,\iota Y_{s,t}(y+\delta e)+(1-\iota)Y_{s,t}(y))d\iota\frac{Y_{s,t}^\delta(y)}{\delta}dW_t\nonumber\\&=&:
\hat{b}_\delta(t,Y_{s,t}(y,\delta e))\frac{Y_{s,t}^\delta(y)}{\delta}dt+\hat{\sigma}_\delta(t,Y_{s,t}(y,\delta e))\frac{Y_{s,t}^\delta(y)}{\delta}dW_t.
\end{eqnarray}
By (\ref{3.22}), (\ref{3.26}) and It\^{o}'s formula, then
\begin{eqnarray}\label{3.27}
&&{\mathbb E}\Big|\frac{Y_{s,t}^\delta(y)}{\delta}-\zeta_{s,t}(y)\Big|^2
\nonumber\\&=&2{\mathbb E}\int_s^t\Big\langle\frac{Y_{s,\tau}^\delta(y)}{\delta}-\zeta_{s,\tau}(y),\hat{b}_\delta(\tau,Y_{s,\tau}(y,
\delta e))\frac{Y_{s,\tau}^\delta(y)}{\delta}-\hat{b}(\tau,Y_{s,\tau}(y))\zeta_{s,\tau}(y)\Big\rangle d\tau\nonumber\\&&+{\mathbb E}\int_s^t\Big\langle\hat{\sigma}_\delta(\tau,Y_{s,\tau}(y,\delta e))\frac{Y_{s,\tau}^\delta(y)}{\delta}-\hat{\sigma}(\tau,Y_{s,\tau}(y))\zeta_{s,\tau}(y),\nonumber\\&& \quad\hat{\sigma}_\delta(\tau,Y_{s,\tau}(y,\delta e))\frac{Y_{s,\tau}^\delta(y)}{\delta}-\hat{\sigma}(\tau,Y_{s,\tau}(y))\zeta_{s,\tau}(y)\Big\rangle d\tau
\nonumber\\&\leq& C{\mathbb E}\int_s^t[1+\kappa^2(\tau)]\Big|\frac{Y_{s,\tau}^\delta(y)}{\delta}-\zeta_{s,\tau}(y)\Big|^2d\tau
\nonumber\\&&+2{\mathbb E}\int_s^t\Big|\frac{Y_{s,\tau}^\delta(y)}{\delta}-\zeta_{s,\tau}(y)\Big|\Big|\frac{Y_{s,\tau}^\delta(y)}{\delta}\Big||\hat{b}_\delta(\tau,Y_{s,\tau}(y,
\delta e))-\hat{b}(\tau,Y_{s,\tau}(y))|d\tau\nonumber\\&&+C{\mathbb E}\int_s^t\Big|\frac{Y_{s,\tau}^\delta(y)}{\delta}\Big|^2\|\hat{\sigma}_\delta(\tau,Y_{s,\tau}(y,\delta e))-\hat{\sigma}(\tau,Y_{s,\tau}(y))\|^2d\tau
\nonumber\\&\leq& C{\mathbb E}
\int_s^t[1+\kappa^2(\tau)]\Big|\frac{Y_{s,\tau}^\delta(y)}{\delta}-\zeta_{s,\tau}(y)\Big|^2d\tau
\nonumber\\&&+
C\int_s^t\Big[{\mathbb E}\Big(\Big|\frac{Y_{s,\tau}^\delta(y)}{\delta}\Big|^4+|\zeta_{s,\tau}(y)|^4\Big)
\Big]^{\frac{1}{2}}
\Big[{\mathbb E}
|\hat{b}_\delta(\tau,Y_{s,\tau}(y,
\delta e))-\hat{b}(\tau,Y_{s,\tau}(y))|^2\Big]^{\frac{1}{2}}d\tau
\nonumber\\&&+
C\int_s^t\Big[{\mathbb E}\Big|\frac{Y_{s,\tau}^\delta(y)}{\delta}\Big|^4\Big]^{\frac{1}{2}}
\Big[{\mathbb E}|\hat{\sigma}_\delta(\tau,Y_{s,\tau}(y,\delta e))-\hat{\sigma}(\tau,Y_{s,\tau}(y))|^4\Big]^{\frac{1}{2}}d\tau.
\end{eqnarray}
From (\ref{3.27}), by using the Gr\"{o}nwall inequality to get
\begin{eqnarray}\label{3.28}
&&\sup_{s\leq t\leq T}{\mathbb E}\Big|\frac{Y_{s,t}^\delta(y)}{\delta}-\zeta_{s,t}(y)\Big|^2
\nonumber\\ &\leq&  C
\int_s^T\Big[{\mathbb E}\Big(\Big|\frac{Y_{s,\tau}^\delta(y)}{\delta}\Big|^4+|\zeta_{s,\tau}(y)|^4\Big)\Big]^{\frac{1}{2}}
\Big[{\mathbb E}
|\hat{b}_\delta(\tau,Y_{s,\tau}(y,
\delta e))-\hat{b}(\tau,Y_{s,\tau}(y))|^2\Big]^{\frac{1}{2}}d\tau
\nonumber\\&&+
C\int_s^T\Big[{\mathbb E}\Big|\frac{Y_{s,\tau}^\delta(y)}{\delta}\Big|^4\Big]^{\frac{1}{2}}
\Big[{\mathbb E}|\hat{\sigma}_\delta(\tau,Y_{s,\tau}(y,\delta e))-\hat{\sigma}(\tau,Y_{s,\tau}(y))|^4\Big]^{\frac{1}{2}}d\tau
\nonumber\\&\leq&  C
\int_s^T\Big[{\mathbb E}
|\hat{b}_\delta(\tau,Y_{s,\tau}(y,
\delta e))-\hat{b}(\tau,Y_{s,\tau}(y))|^2\Big]^{\frac{1}{2}}d\tau
\nonumber\\&&+
C\int_s^T\Big[{\mathbb E}|\hat{\sigma}_\delta(\tau,Y_{s,\tau}(y,\delta e))-\hat{\sigma}(\tau,Y_{s,\tau}(y))|^4\Big]^{\frac{1}{2}}d\tau,
\end{eqnarray}
where in the second inequality we have used (\ref{3.23}) and (\ref{3.25}).

\medskip
Observing that $\hat{b}_\delta\in L^\infty(0,T;{\mathcal C}_b^\theta({\mathbb R}^d))$ and $\hat{\sigma}_\delta\in L^2(0,T;{\mathcal C}_b^\theta({\mathbb R}^d))$ for every $\theta\in (0,1+\alpha-2/q)$, and
\begin{equation}\label{3.29}
\begin{cases}
\hat{b}_\delta(t,Y_{s,t}(y,\delta e)) \longrightarrow \hat{b}(t,Y_{s,t}(y)), \ \ t\in [s,T],
\ {\mathbb P}-a.s., \\
\hat{\sigma}_\delta(t,Y_{s,t}(y,\delta e))\longrightarrow \hat{\sigma}(t,Y_{s,t}(y)), \ t\in [s,T],
\ {\mathbb P}-a.s.,
\end{cases}
\end{equation}
then by the Lebesgue dominated convergence theorem, we deduce from (\ref{3.28}) and (\ref{3.29}) that
\begin{eqnarray}\label{3.30}
\lim_{\delta\rightarrow 0}\sup_{s\leq t\leq T}{\mathbb E}\Big|\frac{Y_{s,t}^\delta(y)}{\delta}-\zeta_{s,t}(y)\Big|^2
=0.
\end{eqnarray}
So, (\ref{3.21}) exists in $L^2(\Omega)$.  Moreover, the weak derivative satisfies SDE (\ref{3.22}) and moment estimate (\ref{3.23}).

\medskip
$\bullet$ \textbf{H\"{o}lder estimates.} For ease of notations, we write 
$\zeta_{s,t}(x)-\zeta_{s,t}(y)$, $U(t,\Psi(t,Y_{s,t}(x)))$ and $U(t,\Psi(t,Y_{s,t}(y)))$ by 
$\zeta_{s,t}(x,y),U(\Psi(Y(x)))$ and $U(\Psi(Y(y)))$, respectively. Then for every $\tilde{p}\geq 2$, an application of It\^o's  formula yields
\begin{eqnarray}\label{3.31}
&&d |\zeta_{s,t}(x,y)|^{\tilde{p}}\nonumber\\
&=&\tilde{p}\lambda|\zeta_{s,t}(x,y)|^{\tilde{p}-2}\langle \zeta_{s,t}(x,y), \nabla U(\Psi(Y(x)))\nabla\Psi(Y(x))\zeta_{s,t}(x)\nonumber\\&&-\nabla U(\Psi(Y(y)))\nabla\Psi(Y(y))\zeta_{s,t}(y)\rangle dt\nonumber\\&&+ \frac{1}{2}\tilde{p}(\tilde{p}-1)|\zeta_{s,t}(x,y)|^{\tilde{p}-2}\langle\nabla^2 U(\Psi(Y(x)))\nabla\Psi(Y(x))\sigma(t)\zeta_{s,t}(x)-\nabla^2 U(\Psi(Y(y))) \nonumber\\&&\cdot\nabla\Psi(Y(y))\sigma(t)\zeta_{s,t}(y),\nabla^2 U(\Psi(Y(x)))\nabla\Psi(Y(x))\sigma(t)\zeta_{s,t}(x)\nonumber\\&&-\nabla^2 U(\Psi(Y(y)))\nabla\Psi(Y(y))\sigma(t)\zeta_{s,t}(y)\rangle dt \nonumber\\&&+\tilde{p}|\zeta_{s,t}(x,y)|^{\tilde{p}-2}\langle \zeta_{s,t}(x,y), [\nabla^2 U(\Psi(Y(x)))\nabla\Psi(Y(x))\zeta_{s,t}(x)\nonumber\\&&\quad-\nabla^2 U(\Psi(Y(y)))\nabla\Psi(Y(y))\zeta_{s,t}(y)]\sigma(t)dW_t\rangle \nonumber\\
&\leq& C(\tilde{p})\Big[|\zeta_{s,t}(x,y)|^{\tilde{p}-1}|\nabla U(\Psi(Y(x)))\nabla\Psi(Y(x))\zeta_{s,t}(x)-\nabla U(\Psi(Y(y)))\nabla\Psi(Y(y))\zeta_{s,t}(y)|\nonumber\\&&+|\zeta_{s,t}(x,y)|^{\tilde{p}-2}|\nabla^2 U(\Psi(Y(x)))\nabla\Psi(Y(x))\sigma(t)\zeta_{s,t}(x)-\nabla^2 U(\Psi(Y(y)))\nabla\Psi(Y(y))
\nonumber\\&&\quad \cdot
\sigma(t)\zeta_{s,t}(y)|^2\Big]dt+\tilde{p}|\zeta_{s,t}(x,y)|^{\tilde{p}-2}\langle \zeta_{s,t}(x,y), [\nabla^2 U(\Psi(Y(x)))\nabla\Psi(Y(x))\sigma(t)\zeta_{s,t}(x)\nonumber\\&&\quad-\nabla^2 U(\Psi(Y(y)))\nabla\Psi(Y(y))\sigma(t)\zeta_{s,t}(y)]dW_t\rangle.
\end{eqnarray}
Noting that $U\in (\tilde{{\mathcal H}}^{2,\theta}_{q,T})^d$ (see (\ref{2.7})) for every $\theta\in (0,1+\alpha-2/q)$ and $\Phi(t,x)=x+U(t,x)$. Additionally $\Psi(t,\cdot)=\Phi^{-1}(t,\cdot)$ and (\ref{3.3}) hold, we then have $\nabla \Phi, \nabla \Psi\in L^\infty(0,T;{\mathcal C}_b^{\theta}({\mathbb R}^d;{\mathbb R}^{d\times d}))$. On the other hand,
\begin{eqnarray*}
&&\nabla U(\Psi(Y(x)))\nabla\Psi(Y(x))\zeta_{s,t}(x)-\nabla U(\Psi(Y(y)))\nabla\Psi(Y(y))\zeta_{s,t}(y)\\ &=&\nabla U(\Psi(Y(x)))\nabla\Psi(Y(x))\zeta_{s,t}(x)-\nabla U(\Psi(Y(y)))\nabla\Psi(Y(x))\zeta_{s,t}(x)\\&&+\nabla U(\Psi(Y(y)))\nabla\Psi(Y(x))\zeta_{s,t}(x)-\nabla U(\Psi(Y(y)))\nabla\Psi(Y(y))\zeta_{s,t}(x)\\&&+\nabla U(\Psi(Y(y)))\nabla\Psi(Y(y))\zeta_{s,t}(x)-\nabla U(\Psi(Y(y)))\nabla\Psi(Y(y))\zeta_{s,t}(y).
\end{eqnarray*}
Therefore,
\begin{eqnarray}\label{3.32}
&&|\nabla U(\Psi(Y(x)))\nabla\Psi(Y(x))\zeta_{s,t}(x)-\nabla U(\Psi(Y(y)))\nabla\Psi(Y(y))\zeta_{s,t}(y)|\nonumber\\&\leq &
\|\nabla^2 U(t)\|_0\|\nabla \Psi\|^2_{\infty,0}|Y_{s,t}(x,y)||\zeta_{s,t}(x)|+ \|\nabla U\|_{\infty,0}[\nabla\Psi]_{\infty,\theta}|Y_{s,t}(x,y)|^\theta|\zeta_{s,t}(x)|\nonumber\\&&+ \|\nabla U\|_{\infty,0}\|\nabla \Psi\|_{\infty,0}|\zeta_{s,t}(x,y)|\nonumber\\&\leq & C\Big[\kappa(t)|Y_{s,t}(x,y)||\zeta_{s,t}(x)|+|Y_{s,t}(x,y)|^\theta|\zeta_{s,t}(x)|+|\zeta_{s,t}(x,y)|\Big],
\end{eqnarray}
where $\hat{\kappa}(t)=[\nabla^2 U(t)]_\theta\in L^2(0,T)$ and $\kappa(t)$ is given in (\ref{3.13}).

\medskip
Similarly, we obtain
\begin{eqnarray}\label{3.33}
&&|\nabla^2 U(\Psi(Y(x)))\nabla\Psi(Y(x))\sigma(t)\zeta_{s,t}(x)-\nabla^2 U(\Psi(Y(y)))\nabla\Psi(Y(y))\sigma(t)\zeta_{s,t}(y)|
\nonumber\\&\leq&|\nabla^2 U(\Psi(Y(x)))\nabla\Psi(Y(x))\sigma(t)\zeta_{s,t}(x)-\nabla^2 U(\Psi(Y(y)))\nabla\Psi(Y(x))\sigma(t)\zeta_{s,t}(x)|\nonumber\\&&+|\nabla^2 U(\Psi(Y(y)))\nabla\Psi(Y(x))\sigma(t)\zeta_{s,t}(x)-\nabla^2 U(\Psi(Y(y)))\nabla\Psi(Y(y))\sigma(t)\zeta_{s,t}(x)|\nonumber\\&&+|\nabla^2 U(\Psi(Y(y)))\nabla\Psi(Y(y))\sigma(t)\zeta_{s,t}(x)-\nabla^2 U(\Psi(Y(y)))\nabla\Psi(Y(y))\sigma(t)\zeta_{s,t}(y)|
\nonumber\\&\leq &\|\sigma\|_{L^\infty(0,T)}\Big[ [\nabla^2 U(t)]_\theta\|\nabla \Psi\|^{1+\theta}_{\infty,0}|Y_{s,t}(x,y)|^\theta|\zeta_{s,t}(x)|
\nonumber\\&&
+ \|\nabla^2 U(t)\|_0[\nabla \Psi]_{\infty,\theta}|Y_{s,t}(x,y)|^\theta|\zeta_{s,t}(x)|+ \|\nabla^2 U(t)\|_0\|\nabla \Psi\|_{\infty,0}|\zeta_{s,t}(x,y)|\Big]\nonumber\\&\leq & C\Big[(\kappa(t)+\hat{\kappa}(t))|Y_{s,t}(x,y)|^\theta |\zeta_{s,t}(x)| +\kappa(t)|\zeta_{s,t}(x,y)|\Big].
\end{eqnarray}

\medskip
Summing over (\ref{3.31})--(\ref{3.33}), for $t\in [s,T]$, we arrive at
\begin{eqnarray}\label{3.34}
&&{\mathbb E}|\zeta_{s,t}(x,y)|^{\tilde{p}}\nonumber\\ &\leq& C {\mathbb E}\int_s^t[1+\kappa^2(\tau)]|\zeta_{s,\tau}(x,y)|^{\tilde{p}}d\tau+ C{\mathbb E}\int_s^t\kappa(\tau)|\zeta_{s,\tau}(x,y)|^{\tilde{p}-1}|Y_{s,\tau}(x,y)|
|\zeta_{s,\tau}(x)|d\tau\nonumber\\&& +C{\mathbb E}\int_s^t|\zeta_{s,\tau}(x,y)|^{\tilde{p}-1}|Y_{s,\tau}(x,y)|^\theta
|\zeta_{s,\tau}(x)|d\tau
\nonumber\\&& +C{\mathbb E}\int_s^t\tilde{\kappa}^2(\tau)|\zeta_{s,\tau}(x,y)|^{\tilde{p}-2}
|Y_{s,\tau}(x,y)|^{2\theta}|\zeta_{s,\tau}(x)|^2d\tau
\nonumber\\ &\leq&C\int_s^t[1+\tilde{\kappa}^2(\tau)] {\mathbb E}|\zeta_{s,\tau}(x,y)|^{\tilde{p}}d\tau + C\int_s^t[1+\tilde{\kappa}^2(\tau)][{\mathbb E}|Y_{s,\tau}(x,y)|^{2\tilde{p}\theta}]^{\frac{1}{2}}
[{\mathbb E}|\zeta_{s,\tau}(x)|^{2\tilde{p}}]^{\frac{1}{2}}d\tau\nonumber\\&& +C\int_s^t\kappa(\tau)[{\mathbb E}|Y_{s,\tau}(x,y)|^{2\tilde{p}}]^{\frac{1}{2}}
[{\mathbb E}|\zeta_{s,\tau}(x)|^{2\tilde{p}}]^{\frac{1}{2}}d\tau,
\end{eqnarray}
where $\tilde{\kappa}=\kappa+\hat{\kappa}\in L^2(0,T)$.

\medskip
By (\ref{3.17}), (\ref{3.23}), (\ref{3.34}) and the Gr\"{o}nwall inequality, we conclude
\begin{eqnarray*}
\sup_{s\leq t\leq T}{\mathbb E}|\zeta_{s,t}(x,y)|^{\tilde{p}}\leq C[|x-y|^{\tilde{p}}+|x-y|^{\tilde{p}\theta}].
\end{eqnarray*}
This, together with the BDG inequality, suggests for every $p\geq 2$ that
\begin{eqnarray}\label{3.35}
&&{\mathbb E}\sup_{s\leq t\leq T}|\zeta_{s,t}(x,y)|^{p}
\nonumber\\&\leq&C \Bigg\{\int_s^T[1+\tilde{\kappa}^2(\tau)]
\Big\{
{\mathbb E}|\zeta_{s,\tau}(x,y)|^{p}+[{\mathbb E}|Y_{s,\tau}(x,y)|^{2p\theta}]^{\frac{1}{2}} +
[{\mathbb E}|Y_{s,\tau}(x,y)|^{2p}]^{\frac{1}{2}}\Big\}d\tau
\nonumber\\&& + \Big[{\mathbb E}\int_s^T\tilde{\kappa}^2(\tau)|\zeta_{s,\tau}(x,y)|^{2p-2}
[|Y_{s,\tau}(x,y)|^{2\theta}\zeta^2_{s,t}(x)+|\zeta_{s,\tau}(x,y)|^2]
d\tau\Big]^{\frac{1}{2}}\Bigg\}\nonumber\\&\leq&C\int_s^T[1+\tilde{\kappa}^2(\tau)]
\Big\{
{\mathbb E}|\zeta_{s,\tau}(x,y)|^{p}+[{\mathbb E}|Y_{s,\tau}(x,y)|^{2p\theta}]^{\frac{1}{2}}+
[{\mathbb E}|Y_{s,\tau}(x,y)|^{2p}]^{\frac{1}{2}}\Big\}d\tau
\nonumber\\&& +C\Big\{\int_s^T\tilde{\kappa}^2(\tau)\Big[{\mathbb E}|\zeta_{s,\tau}(x,y)|^{2p}+
[{\mathbb E}|Y_{s,\tau}(x,y)|^{4p\theta}]^{\frac{1}{2}}\Big]d\tau\Big\}^{\frac{1}{2}}
\nonumber\\&\leq& C[|x-y|^p+|x-y|^{p\theta}].
\end{eqnarray}
By (\ref{3.35}) and the Kolmogorov-Chentsov continuity criterion (or see \cite[Theorem 1.1]{WL}, \cite[Theorem 2.1]{Kuo}), for every $0\leq s\leq t\leq T$,  $\zeta_{s,t}(\cdot)$ has a continuous realization (denoted by itself), which is locally $\beta$-H\"{o}lder continuous in $y$ for every $\beta\in (0,1+\alpha-2/q)$ uniformly in $t$ and satisfies \begin{eqnarray}\label{3.36}
\sup_{0\leq s\leq
T}{\mathbb E}\Big[\sup_{s\leq t \leq T}\Big(\sup_{x,y\in B_R,x\neq y}\frac{|\zeta_{s,t}(x)-\zeta_{s,t}(y)|}{|x-y|^\beta}\Big)^p\Big]<+\infty, \ \ \forall \ p\geq 2, \ \forall \ R>0.
\end{eqnarray}
Then, together with the existence of (\ref{3.21}) in $L^2(\Omega)$, implies that $Y_{s,t}(\cdot)$ is weakly differentiable, which also implies the weak differentiability of $X_{s,t}(\cdot)$. Moreover,  if one differentiates $Y_{s,t}$ with respect to the initial data and denotes the derivative by $\tilde{\zeta}_{s,t}(y)$, then
\begin{eqnarray}\label{3.37}
d \tilde{\zeta}_{s,t}(y)&=&\lambda \nabla U(t,\Psi(t,Y_{s,t}(y)))\nabla\Psi(t,Y_{s,t}(y))\tilde{\zeta}_{s,t}(y)dt\nonumber\\&&+\nabla^2
U(t,\Psi(t,Y_{s,t}(y)))\nabla\Psi(t,Y_{s,t}(y))\sigma(t)\tilde{\zeta}_{s,t}(y)dW_t\nonumber\\
&=&:\hat{b}(t,Y_{s,t}(y))\tilde{\zeta}_{s,t}(y)dt+\hat{\sigma}(t,Y_{s,t}(y))\tilde{\zeta}_{s,t}(y)dW_t, \quad t\in (s,T],
\end{eqnarray}
with $\tilde{\zeta}_{s,t}(y)|_{t=s}=I_{d\times d}$. Additionally, the gradient estimate (\ref{3.23}) holds for $\tilde{\zeta}_{s,t}(y)$. Furthermore, by repeating the above calculations to SDE (\ref{3.37}), (\ref{3.35}) is true for $\tilde{\zeta}_{s,t}(y)$ as well, i.e. for every $p\geq 2$,
\begin{eqnarray}\label{3.38}
\sup_{0\leq s\leq T}{\mathbb E}\sup_{s\leq t \leq T}|\tilde{\zeta}_{s,t}(x,y)|^p \leq C[|x-y|^p+|x-y|^{p\theta}].
\end{eqnarray}
Observing that $X_{s,t}(x)=\Psi(t,Y_{s,t}(\Phi(s,x)))$, we have
\begin{eqnarray*}
\nabla X_{s,t}(x)=\nabla\Psi(t,Y_{s,t}(\Phi(s,x)))\nabla Y_{s,t}(\Phi(s,x)) \nabla\Phi(s,x),
\end{eqnarray*}
which leads to
\begin{eqnarray}\label{3.39}
&&\|\nabla X_{s,t}(x)-\nabla X_{s,t}(y)\|\nonumber\\&=&\|\nabla\Psi(t,Y_{s,t}(\Phi(s,x)))\nabla Y_{s,t}(\Phi(s,x)) \nabla\Phi(s,x)
\nonumber\\&&-\nabla\Psi(t,Y_{s,t}(\Phi(s,y)))\nabla Y_{s,t}(\Phi(s,y)) \nabla\Phi(s,y)\|
\nonumber\\&\leq&\|\nabla \Psi\|_{L^\infty(0,T;{\mathcal C}_b^\theta({\mathbb R}^d))}\|\nabla \Phi\|_{\infty,0}\Big[|Y_{s,t}(\Phi(s,x))-Y_{s,t}(\Phi(s,y))|^\theta\|\nabla Y_{s,t}(\Phi(s,x))\|\nonumber\\&&+
\|\nabla Y_{s,t}(\Phi(s,x))-\nabla Y_{s,t}(\Phi(s,y))\|\Big]\nonumber\\&&+\|\nabla \Psi\|_{\infty,0}\|\nabla Y_{s,t}(\Phi(s,y))\|\| \nabla\Phi(s,x)- \nabla\Phi(s,y)\|,
\end{eqnarray}
for every $x,y\in {\mathbb R}^d$.

\medskip
Noting that $\tilde{\zeta}_{s,t}(x,y)=\nabla Y_{s,t}(x)-\nabla Y_{s,t}(y)$, then  we get
from (\ref{3.17}), (\ref{3.38}) and (\ref{3.39}) that
\begin{eqnarray}\label{3.40}
&&\sup_{0\leq s\leq T}{\mathbb E}\sup_{s\leq t \leq T}\|\nabla X_{s,t}(x)-\nabla X_{s,t}(y)\|^p
\nonumber\\&\leq& C\Big[\sup_{0\leq s\leq T}\Big({\mathbb E}\sup_{s\leq t \leq T}|Y_{s,t}(\Phi(s,x))-Y_{s,t}(\Phi(s,y))|^{2p\theta}\Big)^{\frac{1}{2}}\Big({\mathbb E}\sup_{s\leq t \leq T}\|\nabla Y_{s,t}(\Phi(s,x))\|^{2p}\Big)^{\frac{1}{2}}\nonumber\\&& +\sup_{0\leq s\leq T}{\mathbb E}\sup_{s\leq t \leq T}| \nabla Y_{s,t}(\Phi(s,x))-\nabla Y_{s,t}(\Phi(s,y))|^p
\nonumber\\&&
+\sup_{0\leq s\leq T}{\mathbb E}\sup_{s\leq t \leq T}\|\nabla Y_{s,t}(\Phi(s,y))\|^p|x-y|^{p\theta}\Big]
\nonumber\\&\leq& C[|x-y|^p+|x-y|^{p\theta}].
\end{eqnarray}
By (\ref{3.40}) and the Kolmogorov-Chentsov continuity criterion, for every $0\leq s\leq t\leq T$,  $\nabla X_{s,t}(\cdot)$ has a continuous realization (denoted by itself), which is locally $\beta$-H\"{o}lder continuous in $x$ for every $\beta\in (0,1+\alpha-2/q)$ uniformly in $t$ such that (\ref{1.9}) holds.

\medskip
(iii) Observing that the inverse flows $Y^{n,-1}_{s,t}$ and $Y_{s,t}^{-1}$ satisfy equations which have the same forms as the original ones
beyond the drift and diffusion have opposite signs, it suffices to show the stability for $Y_{s,t}$. Let $U_n$ be the unique strong solution of (\ref{3.1}) with $b$ replaced by  $b_n$. By Theorem \ref{the2.2}, then $U_n$ belongs to $({\mathcal H}^{2,\theta}_{q,T})^d$ for every $\theta\in (0, 1+\alpha-2/q)$ and satisfies (\ref{2.3}) and (\ref{2.4}) uniformly in $n$. Moreover, $U_n-U\in L^2(0,T;{\mathcal C}_b^{2,\theta}({\mathbb R}^d;{\mathbb R}^d))\cap L^\infty(0,T;{\mathcal C}_b^{1,\theta}({\mathbb R}^d;{\mathbb R}^d))$ and
\begin{eqnarray}\label{3.41}
\lim_{n\rightarrow+\infty}\Big[\sup_{0\leq t\leq T}\|U_n(t)-U(t)\|_{{\mathcal C}_b^{1,\theta}({\mathbb R}^d)} +\|\nabla^2U_n-\nabla^2U\|_{L^2(0,T;{\mathcal C}_b^\theta({\mathbb R}^d))}\Big]=0,
\end{eqnarray}
for $q>4/(2+\alpha)$.

\medskip
Let $\Phi_n(t,x)=x+U_n(t,x)$, then $\{\Phi_n\}_{n\geq 1}$ form non-singular diffeomorphisms of class ${\mathcal C}^{1,\theta}$
uniformly in $(t,n)\in [0,T]\times{\mathbb N}$, and
\begin{eqnarray*}
\frac{1}{2}<\sup_{0\leq t\leq T}\|\nabla\Phi_n(t)\|_0 <\frac{3}{2},
\quad  \frac{2}{3}<\sup_{0\leq t\leq T}\|\nabla\Psi_n(t)\|_0<2,
\end{eqnarray*}
where $\Psi_n(t,\cdot)=\Phi^{-1}_n(t,\cdot)$. Notice that
\begin{eqnarray*}
|\Psi_n(t,x)-\Psi(t,x)|&\leq & \sup_{0\leq t\leq T}\|[\nabla\Phi_n(t)]^{-1}\|_0|\Phi_n(t,\Psi_n(t,x))-\Phi_n(t,\Psi(t,x))|\nonumber\\ &\leq & 2|\Phi_n(t,\Psi_n(t,x))-\Phi_n(t,\Psi(t,x))|\nonumber\\&=&2|\Phi(t,\Psi(t,x))-\Phi_n(t,\Psi(t,x))|\leq 2\sup_{0\leq t\leq T}\|U_n(t)-U(t)\|_0
\end{eqnarray*}
and
$$
\nabla \Psi_n(t,x)=[\nabla \Phi_n(t,\Psi_n(t,x))]^{-1},
$$
then
\begin{eqnarray}\label{3.42}
\lim_{n\rightarrow+\infty}\sup_{0\leq t\leq T}\|\Psi_n(t)-\Psi(t)\|_{{\mathcal C}_b^{1,\theta}({\mathbb R}^d)}=0.
\end{eqnarray}
Therefore, it is sufficient to show for every $p\geq 2$,
\begin{eqnarray}\label{3.43}
\lim _{n\rightarrow+\infty}\sup_{y\in {\mathbb R}^d} \sup_{0\leq s\leq
T}{\mathbb E}[\sup_{s\leq t \leq T}|Y^n_{s,t}(y)-Y_{s,t}(y)|^p]=0
\end{eqnarray}
and
\begin{eqnarray}\label{3.44}
\lim_{n\rightarrow+\infty}\sup_{y\in {\mathbb R}^d} \sup_{0\leq s\leq
T}{\mathbb E}[\sup_{s\leq t \leq T}\|\nabla Y^n_{s,t}(y)-\nabla Y_{s,t}(y)\|^p]=0,
\end{eqnarray}
where $Y_{s,t}^n=\Phi_n(t,X_{s,t}^n)$ satisfies the following equation:
\begin{eqnarray}\label{3.45}
\left\{
  \begin{array}{ll}
  dY_{s,t}^n=\lambda U_n(t,\Psi_n(t,Y_{s,t}^n))dt+(I+\nabla U_n(t,\Psi_n(t,Y_{s,t}^n)))\sigma(t)dW_t\\ \quad\quad \ =:
\tilde{b}_n(t,Y_{s,t}^n)dt+\tilde{\sigma}_n(t,Y_{s,t}^n)dW_t,\ t\in(s,T], \\
Y_{s,t}^n|_{t=s}=y=\Phi_n(s,x).
  \end{array}
\right.
\end{eqnarray}

\medskip
For every $\tilde{p}\geq 2$, in view of It\^o's formula
\begin{eqnarray}\label{3.46}
&&d |Y_{s,t}^n-Y_{s,t}|^{\tilde{p}}
\nonumber\\&=&\tilde{p}\lambda|Y_{s,t}^n-Y_{s,t}|^{\tilde{p}-2}\langle Y_{s,t}^n-Y_{s,t}, U_n(t,\Psi_n(t,Y_{s,t}^n))-U(t,\Psi(t,Y_{s,t}))\rangle dt\nonumber\\&&+ \frac{1}{2}\tilde{p}(\tilde{p}-1)|Y_{s,t}^n-Y_{s,t}|^{\tilde{p}-2}tr([\nabla U_n(t,\Psi_n(t,Y_{s,t}^n))-\nabla U(t,\Psi(t,Y_{s,t}))]\sigma(t) \nonumber\\&&\quad\times\sigma^\top(t)[\nabla U_n(t,\Psi_n(t,Y_{s,t}^n))-\nabla U(t,\Psi(t,Y_{s,t}))]^\top)dt \nonumber\\&&+\tilde{p}|Y_{s,t}^n-Y_{s,t}|^{\tilde{p}-2}\langle Y_{s,t}^n-Y_{s,t}, (\nabla U_n(t,\Psi_n(t,Y_{s,t}^n))-\nabla U(t,\Psi(t,Y_{s,t})))\sigma(t)dW_t\rangle
\nonumber\\&\leq&C(\tilde{p})|Y_{s,t}^n-Y_{s,t}|^{\tilde{p}-1}|U_n(t,\Psi_n(t,Y_{s,t}^n))-U(t,\Psi(t,Y_{s,t}))|dt
\nonumber\\&&+C(\tilde{p})|Y_{s,t}^n-Y_{s,t}|^{\tilde{p}-2}|\nabla U_n(t,\Psi_n(t,Y_{s,t}^n))-\nabla U(t,\Psi(t,Y_{s,t}))|^2dt\nonumber\\&&+\tilde{p}|Y_{s,t}^n-Y_{s,t}|^{\tilde{p}-2}\langle Y_{s,t}^n-Y_{s,t}, (\nabla U_n(t,\Psi_n(t,Y_{s,t}^n))-\nabla U(t,\Psi(t,Y_{s,t})))\sigma(t)dW_t\rangle
\nonumber\\&\leq& C(\tilde{p})|Y_{s,t}^n-Y_{s,t}|^{\tilde{p}-1}\Big[\|U_n-U\|_{\infty,0}+\|\nabla U\|_{\infty,0}\|\Psi_n-\Psi\|_{\infty,0}\nonumber\\&&+
\|\nabla U\|_{\infty,0}\|\nabla\Psi\|_{\infty,0}|Y_{s,t}^n-Y_{s,t}|\Big]dt+
C(\tilde{p})|Y_{s,t}^n-Y_{s,t}|^{\tilde{p}-2}\Big[\|\nabla U_n-\nabla U\|^2_{\infty,0}\nonumber\\&&+\|\nabla^2 U(t)\|^2_0\|\Psi_n-\Psi\|^2_{\infty,0}+\|\nabla^2 U(t)\|^2_0\|\nabla\Psi\|^2_{\infty,0}|Y_{s,t}^n-Y_{s,t}|^2\Big]dt
\nonumber\\&&\!+\tilde{p}|Y_{s,t}^n-Y_{s,t}|^{\tilde{p}-2}\langle Y_{s,t}^n\!-\!Y_{s,t}, (\nabla U_n(t,\Psi_n(t,Y_{s,t}^n))-\nabla U(t,\Psi(t,Y_{s,t})))\sigma(t)dW_t\rangle.
\end{eqnarray}
By using H\"{o}lder's inequality, we deduce from (\ref{3.46}) that
\begin{eqnarray}\label{3.47}
&&{\mathbb E}|Y_{s,t}^n-Y_{s,t}|^{\tilde{p}}\nonumber\\ &\leq& C(d,T,\Theta,\tilde{p},[b]_{q,\frac{2}{q}-1},[b]_{q,\alpha}) {\mathbb E}\int_s^t[1+\kappa^2(\tau)]|Y_{s,\tau}^n-Y_{s,\tau}|^{\tilde{p}}d\tau
\nonumber\\&&+C(d,T,\Theta,\tilde{p},[b]_{q,\frac{2}{q}-1},[b]_{q,\alpha})
\Big[\|U_n-U\|^{\tilde{p}}_{L^\infty(0,T;{\mathcal C}_b^1({\mathbb R}^d;{\mathbb R}^d))}+
\|\Psi_n-\Psi\|^{\tilde{p}}_{\infty,0}\Big],
\end{eqnarray}
where $\kappa$ is given in (\ref{3.13}).

\medskip
From (\ref{3.41}), (\ref{3.42}) and (\ref{3.47}), it follows that
\begin{eqnarray}\label{3.48}
&&\lim _{n\rightarrow+\infty}\sup_{y\in {\mathbb R}^d} \sup_{0\leq s\leq
T}\sup_{s\leq t\leq
T}{\mathbb E}[|Y^n_{s,t}(y)-Y_{s,t}(y)|^{\tilde{p}}]=0.
\end{eqnarray}
By (\ref{3.45}), (\ref{3.47}), (\ref{3.48}) and BDG's inequality to get for every $p\geq 2$ that
\begin{eqnarray}\label{3.49}
&&\lim _{n\rightarrow+\infty}\sup_{y\in {\mathbb R}^d}\sup_{0\leq s\leq T}{\mathbb E}\sup_{s\leq t \leq T}| Y_{s,t}^n(y)-Y_{s,t}(y)|^p \nonumber\\ &\leq& C
\lim _{n\rightarrow+\infty}\Big[\|U_n-U\|^p_{L^\infty(0,T;{\mathcal C}_b^1({\mathbb R}^d))}+
\|\Psi_n-\Psi\|^p_{\infty,0}\int_0^T[1+\kappa^2(\tau)]d\tau\Big]
\nonumber\\&&+
C\lim _{n\rightarrow+\infty}\Big[\sup_{0\leq s\leq T}{\mathbb E}\int_s^T\|Y_{s,t}^n(y)-Y_{s,t}(y)|^{2p-2}
\nonumber\\&&\qquad \times\|\nabla U_n(t,\Psi_n(t,Y_{s,t}^n))-\nabla U(t,\Psi(t,Y_{s,t}))\|^2 dt\Big]^{\frac{1}{2}}\nonumber\\ &\leq& C
\lim _{n\rightarrow+\infty}\Big[\|U_n-U\|^p_{L^\infty(0,T;{\mathcal C}_b^1({\mathbb R}^d;{\mathbb R}^d))}+
\|\Psi_n-\Psi\|^p_{\infty,0}\Big]
\nonumber\\&&+
C\lim _{n\rightarrow+\infty}\Big[\sup_{0\leq s\leq T}{\mathbb E}\int_s^T\|Y_{s,t}^n(y)-Y_{s,t}(y)|^{2p-2}dt\Big]^{\frac{1}{2}}=0,
\end{eqnarray}
which implies (\ref{3.43}).

\medskip
For the gradient, we get an analogue of (\ref{3.46}) that
\begin{eqnarray}\label{3.50}
&&d \|\tilde{\zeta}_{s,t}^n(y)-\tilde{\zeta}_{s,t}(y)\|^{\tilde{p}}\nonumber\\&\leq& C\Big[\|
\tilde{\zeta}_{s,t}^n-\tilde{\zeta}_{s,t}\|^{\tilde{p}-1}\|\nabla U_n(\Psi_n(Y_{s,t}^n))\nabla\Psi_n(Y_{s,t}^n)\tilde{\zeta}_{s,t}^n-\nabla U(\Psi(Y_{s,t}))\nabla\Psi(Y_{s,t})\tilde{\zeta}_{s,t}\|\nonumber\\&&+\|\tilde{\zeta}_{s,t}^n-\tilde{\zeta}_{s,t}\|^{\tilde{p}-2}\|\nabla^2 U_n(\Psi_n(Y_{s,t}^n))\nabla\Psi_n(Y_{s,t}^n)\tilde{\zeta}_{s,t}^n-\nabla^2 U(\Psi(Y_{s,t}))\nabla\Psi(Y_{s,t})\tilde{\zeta}_{s,t}\|^2\Big]dt \nonumber\\&&+\tilde{p}\|\tilde{\zeta}_{s,t}^n-\tilde{\zeta}_{s,t}\|^{\tilde{p}-2}\langle \tilde{\zeta}_{s,t}^n-\tilde{\zeta}_{s,t}, [\nabla^2 U_n(\Psi_n(Y_{s,t}^n))\nabla\Psi_n(Y_{s,t}^n)\tilde{\zeta}_{s,t}^n\nonumber\\&&\quad-\nabla^2 U(\Psi(Y_{s,t}))\nabla\Psi(Y_{s,t})\tilde{\zeta}_{s,t}]\sigma(t)dW_t\rangle,
\end{eqnarray}
where $\tilde{\zeta}_{s,t}^n(y)|_{t=s}=I_{d\times d}$, $\tilde{\zeta}_{s,t}^n(y)=\nabla Y_{s,t}^n(y)$ and $\tilde{\zeta}_{s,t}(y) =\nabla Y_{s,t}(y)$.

\medskip
Moreover, by  analogue calculations of (\ref{3.32}) and (\ref{3.33}), we get
\begin{eqnarray}\label{3.51}
&&\|\nabla U_n(\Psi_n(Y_{s,t}^n))\nabla\Psi_n(Y_{s,t}^n)\tilde{\zeta}_{s,t}^n-\nabla U(\Psi(Y_{s,t}))\nabla\Psi(Y_{s,t})\tilde{\zeta}_{s,t}\|
\nonumber\\&\leq &\|\nabla U_n(\Psi_n(Y_{s,t}^n))\nabla\Psi_n(Y_{s,t}^n)\tilde{\zeta}_{s,t}^n-\nabla U(\Psi_n(Y_{s,t}^n))\nabla\Psi_n(Y_{s,t}^n)\tilde{\zeta}_{s,t}^n\|\nonumber\\&&+\|\nabla U(\Psi_n(Y_{s,t}^n))\nabla\Psi_n(Y_{s,t}^n)\tilde{\zeta}_{s,t}^n-\nabla U(\Psi(Y_{s,t}^n))\nabla\Psi_n(Y_{s,t}^n)\tilde{\zeta}_{s,t}^n\|\nonumber\\&&+\|\nabla U(\Psi(Y_{s,t}^n))\nabla\Psi_n(Y_{s,t}^n)\tilde{\zeta}_{s,t}^n-\nabla U(\Psi(Y_{s,t}))\nabla\Psi_n(Y_{s,t}^n)\tilde{\zeta}_{s,t}^n\|\nonumber\\&&+\|\nabla U(\Psi(Y_{s,t}))\nabla\Psi_n(Y_{s,t}^n)\tilde{\zeta}_{s,t}^n-\nabla U(\Psi(Y_{s,t}))\nabla\Psi(Y_{s,t}^n)\tilde{\zeta}_{s,t}^n\|\nonumber\\&&+\|\nabla U(\Psi(Y_{s,t}))\nabla\Psi(Y_{s,t}^n)\tilde{\zeta}_{s,t}^n-\nabla U(\Psi(Y_{s,t}))\nabla\Psi(Y_{s,t})\tilde{\zeta}_{s,t}^n\|\nonumber\\&&+\|\nabla U(\Psi(Y_{s,t}))\nabla\Psi(Y_{s,t})\tilde{\zeta}_{s,t}^n-\nabla U(\Psi(Y_{s,t}))\nabla\Psi(Y_{s,t})\tilde{\zeta}_{s,t}\|
\nonumber\\
&\leq &
\|\nabla U_n-\nabla U\|_{\infty,0}\|\nabla \Psi_n\|_{\infty,0}\|\tilde{\zeta}_{s,t}^n\|
+\|\nabla^2 U(t)\|_0\|\Psi_n-\Psi\|_{\infty,0}\|\nabla \Psi_n\|_{\infty,0}\|\tilde{\zeta}_{s,t}^n\|\nonumber\\&&
+\|\nabla^2 U(t)\|_0\|\nabla\Psi\|_{\infty,0}\|\nabla \Psi_n\|_{\infty,0}\|\tilde{\zeta}_{s,t}^n\||Y^n_{s,t}- Y_{s,t}|
+\|\nabla U\|_{\infty,0}\|\nabla \Psi_n-\nabla\Psi\|_{\infty,0}
\|\tilde{\zeta}_{s,t}^n\|\nonumber\\&&
+\|\nabla U\|_{\infty,0}[\nabla \Psi]_{\infty,\theta}|Y^n_{s,t}- Y_{s,t}|^\theta\|\tilde{\zeta}_{s,t}^n\|
+\|\nabla U\|_{\infty,0}\|\nabla\Psi\|_{\infty,0}\|\tilde{\zeta}_{s,t}^n-\tilde{\zeta}_{s,t}\|
\nonumber\\
&\leq &
C \Big[\|\nabla U_n-\nabla U\|_{\infty,0}
+(1+\kappa(t))\|\Psi_n-\Psi\|_{L^\infty(0,T;{\mathcal C}_b^1({\mathbb R}^d))}\nonumber\\&&
+\kappa(t)|Y^n_{s,t}- Y_{s,t}|
+|Y^n_{s,t}- Y_{s,t}|^\theta\Big]\|\tilde{\zeta}_{s,t}^n\|
+C\|\tilde{\zeta}_{s,t}^n-\tilde{\zeta}_{s,t}\|
\end{eqnarray}
and
\begin{eqnarray}\label{3.52}
&&\|\nabla^2 U_n(\Psi_n(Y_{s,t}^n))\nabla\Psi_n(Y_{s,t}^n)\tilde{\zeta}_{s,t}^n-\nabla^2 U(\Psi(Y_{s,t}))\nabla\Psi(Y_{s,t})\tilde{\zeta}_{s,t}\|
\nonumber\\&\leq &\|\nabla^2 U_n(\Psi_n(Y_{s,t}^n))\nabla\Psi_n(Y_{s,t}^n)\tilde{\zeta}_{s,t}^n-\nabla^2 U(\Psi_n(Y_{s,t}^n))\nabla\Psi_n(Y_{s,t}^n)\tilde{\zeta}_{s,t}^n\|\nonumber\\&&+\|\nabla^2 U(\Psi_n(Y_{s,t}^n))\nabla\Psi_n(Y_{s,t}^n)\tilde{\zeta}_{s,t}^n-\nabla^2 U(\Psi(Y_{s,t}^n))\nabla\Psi_n(Y_{s,t}^n)\tilde{\zeta}_{s,t}^n\|\nonumber\\&&+\|\nabla^2 U(\Psi(Y_{s,t}^n))\nabla\Psi_n(Y_{s,t}^n)\tilde{\zeta}_{s,t}^n-\nabla^2 U(\Psi(Y_{s,t}))\nabla\Psi_n(Y_{s,t}^n)\tilde{\zeta}_{s,t}^n\|\nonumber\\&&+\|\nabla^2 U(\Psi(Y_{s,t}))\nabla\Psi_n(Y_{s,t}^n)\tilde{\zeta}_{s,t}^n-\nabla^2 U(\Psi(Y_{s,t}))\nabla\Psi(Y_{s,t}^n)\tilde{\zeta}_{s,t}^n\|\nonumber\\&&+\|\nabla^2 U(\Psi(Y_{s,t}))\nabla\Psi(Y_{s,t}^n)\tilde{\zeta}_{s,t}^n-\nabla^2 U(\Psi(Y_{s,t}))\nabla\Psi(Y_{s,t})\tilde{\zeta}_{s,t}^n\|\nonumber\\&&+\|\nabla^2 U(\Psi(Y_{s,t}))\nabla\Psi(Y_{s,t})\tilde{\zeta}_{s,t}^n-\nabla^2 U(\Psi(Y_{s,t}))\nabla\Psi(Y_{s,t})\tilde{\zeta}_{s,t}\|
\nonumber\\&\leq&C\Big[\|\nabla^2 U_n(t)-\nabla^2U(t)\|_0
+\hat{\kappa}(t)\|\Psi_n-\Psi\|^\theta_{\infty,0}
+\tilde{\kappa}(t)|Y^n_{s,t}- Y_{s,t}|^\theta
\nonumber\\&&+\kappa(t)\|\nabla\Psi_n-\nabla\Psi\|_{\infty,0}\Big]\|\tilde{\zeta}_{s,t}^n\|
+C\kappa(t)\|\tilde{\zeta}_{s,t}^n-\tilde{\zeta}_{s,t}\|,
\end{eqnarray}
where $\tilde{\kappa}$ is given in (\ref{3.34}). Combining  (\ref{3.50}), (\ref{3.51}) and (\ref{3.52}), one asserts
\begin{eqnarray}\label{3.53}
&&{\mathbb E}\|\tilde{\zeta}_{s,t}^n(y)-\tilde{\zeta}_{s,t}(y)\|^{\tilde{p}}\nonumber\\&\leq& C\int_s^t[1+\tilde{\kappa}^2(\tau)]{\mathbb E}\|\tilde{\zeta}_{s,\tau}^n-\tilde{\zeta}_{s,\tau}\|^{\tilde{p}}d\tau
+C\Big[\|\nabla U_n-\nabla U\|^{\tilde{p}}_{\infty,0}+
\|\Psi_n-\Psi\|^{\tilde{p}}_{L^\infty(0,T;{\mathcal C}_b^1({\mathbb R}^d;{\mathbb R}^d))}\nonumber\\&&+
\|\Psi_n-\Psi\|^{\tilde{p}\theta}_{\infty,0}\Big]
\int_s^t[1+\kappa(\tau)]{\mathbb E}\|\tilde{\zeta}_{s,\tau}^n\|^{\tilde{p}}d\tau
\nonumber\\&&+
C\int_s^t[1+\tilde{\kappa}^2(\tau)]{\mathbb E}\Big[(|Y^n_{s,\tau}-Y_{s,\tau}|^{\tilde{p}}+|Y^n_{s,\tau}
-Y_{s,\tau}|^{\tilde{p}\theta})\|\tilde{\zeta}_{s,\tau}^n\|^{\tilde{p}}\Big]d\tau
\nonumber\\&&
+C\int_s^t{\mathbb E}\Big[\|\tilde{\zeta}_{s,\tau}^n-\tilde{\zeta}_{s,\tau}\|^{\tilde{p}-2}\|\nabla^2 U_n(\tau)-\nabla^2 U(\tau)\|_0^2\|\tilde{\zeta}_{s,\tau}^n\|^2\Big]d\tau.
\end{eqnarray}

Observing that for every $\tilde{q}\geq 2$, there is a positive constant $C$ such that
\begin{eqnarray}\label{3.54}
\sup_{n}\sup_{y\in {\mathbb R}^d}\sup_{0\leq s\leq T}{\mathbb E}\sup_{s\leq t\leq T}[\|\tilde{\zeta}_{s,t}^n(y)\|^{\tilde{q}}+\|\tilde{\zeta}_{s,t}(y)\|^{\tilde{q}}]\leq C(d,T,\Theta,\tilde{q},\|b\|_{L^q(0,T;{\mathcal C}^{\frac{2}{q}-1}\cap{\mathcal C}^{\alpha}({\mathbb R}^d;{\mathbb R}^d))}),
\end{eqnarray}
then by (\ref{3.53}) and (\ref{3.54}), it leads to
\begin{eqnarray*}
&&{\mathbb E}\|\tilde{\zeta}_{s,t}^n(y)-\tilde{\zeta}_{s,t}(y)\|^{\tilde{p}}\nonumber\\&\leq& C\int_s^t[1+\tilde{\kappa}^2(\tau)]{\mathbb E}\|\tilde{\zeta}_{s,\tau}^n-\tilde{\zeta}_{s,\tau}\|^{\tilde{p}}d\tau
+C\Big[\|\nabla U_n-\nabla U\|^{\tilde{p}}_{\infty,0}+
\|\Psi_n-\Psi\|^{\tilde{p}}_{L^\infty(0,T;{\mathcal C}_b^1({\mathbb R}^d;{\mathbb R}^d))}\nonumber\\&&+
\|\Psi_n-\Psi\|^{\tilde{p}\theta}_{\infty,0}\Big]+
C\Big[\Big(\sup_{0\leq s\leq T}\sup_{s\leq \tau\leq T}{\mathbb E}[|Y^n_{s,\tau}-Y_{s,\tau}|^{2\tilde{p}}
+|Y^n_{s,\tau}-Y_{s,\tau}|^{2\tilde{p}\theta}]\Big)^{\frac{1}{2}}
\nonumber\\&&+\|\nabla^2U_n-\nabla^2 U\|_{2,0}^2\Big].
\end{eqnarray*}
Consequently,
\begin{eqnarray}\label{3.55}
&&\sup_{0\leq s\leq T}\sup_{s\leq t\leq T}{\mathbb E}\|\tilde{\zeta}_{s,t}^n(y)-\tilde{\zeta}_{s,t}(y)\|^{\tilde{p}}\nonumber\\&\leq& C
\Big[\|\nabla U_n-\nabla U\|^{\tilde{p}}_{\infty,0}+
\|\Psi_n-\Psi\|^{\tilde{p}}_{L^\infty(0,T;{\mathcal C}_b^1({\mathbb R}^d;{\mathbb R}^d))}+
\|\Psi_n-\Psi\|^{\tilde{p}\theta}_{\infty,0}\Big]+
\nonumber\\&&
C\Big[\Big(\sup_{0\leq s\leq T}\sup_{s\leq \tau\leq T}{\mathbb E}[|Y^n_{s,\tau}-Y_{s,\tau}|^{2\tilde{p}}
+|Y^n_{s,\tau}-Y_{s,\tau}|^{2\tilde{p}\theta}]\Big)^{\frac{1}{2}}+\|\nabla^2U_n-\nabla^2 U\|_{2,0}^2\Big].
\end{eqnarray}
By (\ref{3.41}), (\ref{3.42}), (\ref{3.48}) and (\ref{3.55}), we get
\begin{eqnarray}\label{3.57}
\lim_{n\rightarrow+\infty}\sup_{0\leq s\leq T}\sup_{s\leq t\leq T}{\mathbb E}\|\tilde{\zeta}_{s,t}^n(y)-\tilde{\zeta}_{s,t}(y)\|^{\tilde{p}}=0.
\end{eqnarray}
Repeating a similar calculation as (\ref{3.49}), by (\ref{3.57}) we end up with (\ref{3.44}). $\Box$

\section{Proof of Theorem \ref{the1.7}}
\label{sec4}\setcounter{equation}{0}
(i) Firstly, we prove that $u(t,x)=u_0(X^{-1}_t(x))$ is a stochastic weak solution of (\ref{1.14}). Observing that $\int_{{\mathbb R}^d}u(t,x)\varphi(x)dx$ is ${\mathcal F}_t$-adapted for every $\varphi\in {\mathcal C}_0^\infty({\mathbb R}^d)$, it is sufficient to show $u\in L^\infty(\Omega\times[0,T];L^r({\mathbb R}^d))$ and (\ref{1.15}) holds.

\medskip
Since $u_0\in L^r({\mathbb R}^d)$, if $r=+\infty$ we clearly have $u\in L^\infty(\Omega\times[0,T]\times{\mathbb R}^d)$, and if $r<+\infty$, with the help of Euler's identity, it follows that
\begin{eqnarray}\label{4.1}
\int_{{\mathbb R}^d}|u_0(X^{-1}_t(x))|^rdx&=&\int_{{\mathbb R}^d}|u_0(x)|^r{\rm det}(\nabla X_t(x))dx
\nonumber\\&=&\int_{{\mathbb R}^d}|u_0(x)|^r\exp\Big(\int^t_0\div  b(\tau,X_\tau(x))d\tau\Big)dx
\nonumber\\&=&\int_{{\mathbb R}^d}|u_0(x)|^rdx.
\end{eqnarray}
Therefore, $u\in L^\infty(\Omega\times[0,T];L^r({\mathbb R}^d))$.

\medskip
Let $b_n$ be given in (iii) of Theorem \ref{the1.3}, and let $X^n_t$ be the unique solution of (\ref{1.1}) with $\sigma=I_{d\times d}$, $s=0$ and $b$ replaced by $b_n$. Let $X^{n,-1}_\cdot$ be the inverse of $X^n_\cdot$. Since $b_n$ is smooth in space variable, in view of the characteristic lines and It\^{o}'s formula,  $u_n(t,x)=u_0(X^{n,-1}_t(x))$ is the unique stochastic weak solution of the following Cauchy problem:
\begin{eqnarray}\label{4.2}
\left\{
  \begin{array}{ll}
\partial_tu_n(t,x)+b_n(t,x)\cdot\nabla u_n(t,x)
+\sum\limits_{i=1}^d\partial_{x_i}u_n(t,x)\circ\dot{W}_{i,t}=0, \ (t,x)\in(0,T]\times {\mathbb R}^d, \\
u_n(t,x)|_{t=0}=u_0(x), \  x\in{\mathbb R}^d.
  \end{array}
\right.
\end{eqnarray}
Observing that $\div b_n=(\div b)\ast\rho_n=0$, for every
$\varphi\in{\mathcal C}_0^\infty({\mathbb R}^d)$ and
 every  $t\in [0,T]$, then
\begin{eqnarray}\label{4.3}
\int_{{\mathbb R}^d}u_0(X^{n,-1}_t(x))\varphi(x)dx&=&\int_{{\mathbb R}^d}u_0(x)\varphi(x)dx+
\int^t_0\int_{{\mathbb R}^d}u_0(X^{n,-1}_\tau(x))\div(b_n(\tau,x)\varphi(x))dxd\tau\nonumber\\&&
+\sum_{i=1}^d\int^t_0\circ dW_{i,\tau}\int_{{\mathbb R}^d}u_0(X^{n,-1}_\tau(x))\partial_{x_i}\varphi(x)dx\nonumber\\&=&\int_{{\mathbb R}^d}u_0(x)\varphi(x)dx+
\int^t_0\int_{{\mathbb R}^d}u_0(X^{n,-1}_\tau(x))b_n(\tau,x)\cdot\nabla\varphi(x)dxd\tau\nonumber\\&&
+\sum_{i=1}^d\int^t_0\circ dW_{i,\tau}\int_{{\mathbb R}^d}u_0(X^{n,-1}_\tau(x))\partial_{x_i}\varphi(x)dx \quad
{\mathbb P}-a.s..
\end{eqnarray}
Thanks to Theorem \ref{the1.3}, $\div b=0$ and the following fact (see (\ref{2.55}) and (\ref{2.56}))
\begin{eqnarray}\label{4.4}
 \lim\limits_{n\rightarrow+\infty}\|b_n-b\|_{L^q(0,T;{\mathcal C}_b^{\alpha^\prime}({\mathbb R}^d;{\mathbb R}^d))}=0, \quad \forall \ \alpha^\prime\in (0,\alpha),
\end{eqnarray}
if one lets $n$ tend to infinity in (\ref{4.3}), then (\ref{1.15}) holds for $u(t,x)=u_0(X^{-1}_t(x))$. Thus $u(t,x)=u_0(X^{-1}_t(x))$ is a stochastic weak solution of (\ref{1.14}).

\medskip
Secondly, we prove the uniqueness of stochastic weak solutions. Observing that the equation is linear, it suffices to prove that $u\equiv 0$ a.s. if the initial data vanishes. Let $u(t,x)$ be a stochastic weak solution of (\ref{1.14}) with $u_0=0$. Let $\rho$ be given by (\ref{1.10}). For $m,n\in {\mathbb N}$, we set $\rho_n(x)=n^d\rho(nx)$ and $\rho_m(x)=m^d\rho(mx)$. Then
\begin{eqnarray}\label{4.5}
\int_{{\mathbb R}^d}u(t,y)\rho_m(x-y)dy&=&
-\int^t_0\int_{{\mathbb R}^d}u(\tau,y)b(\tau,y)\cdot\nabla\rho_m(x-y)dyd\tau\nonumber\\&&
-\sum_{i=1}^d\int^t_0\circ dW_{i,\tau}\int_{{\mathbb R}^d}u(\tau,y)\partial_{x_i}\rho_m(x-y)dy, \quad
{\mathbb P}-a.s..
\end{eqnarray}

\medskip
For every $\varphi\in{\mathcal C}_0^\infty({\mathbb R}^d)$, if we set $v_n(t,x)=\varphi(X^{n,-1}_t(x))$, it satisfies
\begin{eqnarray}\label{4.6}
\left\{
  \begin{array}{ll}
\partial_tv_n(t,x)+b_n(t,x)\cdot\nabla v_n(t,x)
+\sum\limits_{i=1}^d\partial_{x_i}v_n(t,x)\circ\dot{W}_{i,t}=0, \ (t,x)\in(0,T]\times {\mathbb R}^d, \\
v_n(t,x)|_{t=0}=\varphi(x), \  x\in{\mathbb R}^d.
  \end{array}
\right.
\end{eqnarray}
Let $u_m=u\ast\rho_m, b_n=b\ast\rho_n$ and $b_m=b\ast\rho_m$. In view of the fact that $\div b=\div b_n=\div b_m=0$, then
\begin{eqnarray}\label{4.7}
&&\int_{{\mathbb R}^d}u_m(t,X^n_t(x))\varphi(x)dx\nonumber\\&=&\int_{{\mathbb R}^d}u_m(t,x)
\varphi(X^{n,-1}_t(x))\det(X^{n,-1}_t(x))dx\nonumber\\&=&\int_{{\mathbb R}^d}u_m(t,x)
\varphi(X^{n,-1}_t(x))\exp\Big(-\int^t_0\div  b_n(\tau,X_\tau^n(X^{n,-1}_t(x)))d\tau\Big)dx\nonumber\\&=&\int_{{\mathbb R}^d}u_m(t,x)v_n(t,x)dx, \quad
{\mathbb P}-a.s.,
\end{eqnarray}
where in the second identity we have used the following Euler identity:
\begin{eqnarray*}
\det(\nabla X^{n,-1}_t(x))=[\det(\nabla X^n_t(X^{n,-1}_t(x)))]^{-1}=
\exp\Big(-\int_0^t\div b_n(\tau,X_\tau^n(X^{n,-1}_t(x)))d\tau\Big).
\end{eqnarray*}
Observing that $ X^{n,-1}_t(x)$ is continuous in $(t,x)$, then for every fixed $\omega$, $\varphi(X^{n,-1}_t(\cdot,\omega))\in {\mathcal C}_0^\infty({\mathbb R}^d)$ uniformly in $t\in [0,T]$.  By (\ref{4.5})--(\ref{4.7}), we derive
\begin{eqnarray}\label{4.8}
\int_{{\mathbb R}^d}u_m(t,X^n_t(x))\varphi(x)dx&=&
-\int^t_0\int_{{\mathbb R}^d}u_m(\tau,x)b_n(\tau,x)\cdot\nabla_x \varphi(X^{n,-1}_\tau(x))dxd\tau
\nonumber\\&&-\sum_{i=1}^d\int^t_0\circ dW_{i,\tau}\int_{{\mathbb R}^d}u_m(\tau,x)\partial_{x_i}\varphi(X^{n,-1}_\tau(x))dx\nonumber\\&&-
\int^t_0\int_{{\mathbb R}^d}\varphi(X^{n,-1}_\tau(x))dx
\int_{{\mathbb R}^d}u(\tau,y)b(\tau,y)\cdot\nabla \rho_m(x-y)dyd\tau\nonumber\\&&
-\sum_{i=1}^d\int^t_0\circ dW_{i,\tau}\int_{{\mathbb R}^d}\varphi(X^{n,-1}_\tau(x))dx
\int_{{\mathbb R}^d}u(\tau,y)\partial_{x_i}\rho_m(x-y)dy
\nonumber\\&=&
\int^t_0\int_{{\mathbb R}^d}\Big[\int_{{\mathbb R}^d}u(\tau,y)b(\tau,y) \rho_m(x-y)dy-u_m(\tau,x)b_n(\tau,x)\Big]\nonumber\\&&\qquad\qquad\cdot \nabla_x\varphi(X^{n,-1}_\tau(x))dxd\tau,\quad
{\mathbb P}-a.s.,
\end{eqnarray}
where in the second identity we have used integration by parts.

\medskip
Letting $m$ tend to infinity and $n\in {\mathbb N}$ be fixed in (\ref{4.8}), we obtain
\begin{eqnarray}\label{4.9}
&&\int_{{\mathbb R}^d}u(t,X^n_t(x))\varphi(x)dx\nonumber\\&=&
\int^t_0\int_{{\mathbb R}^d}u(\tau,x)[b(\tau,x)-b_n(\tau,x)]\cdot\nabla_x \varphi(X^{n,-1}_\tau(x))dxd\tau\nonumber\\&=&
\int^t_0\int_{{\mathbb R}^d}u(\tau,x)[b(\tau,x)-b_n(\tau,x)]\cdot[\nabla \varphi(X^{n,-1}_\tau(x))\nabla X^{n,-1}_\tau(x)]dxd\tau,\quad
{\mathbb P}-a.s..
\end{eqnarray}
By (\ref{1.12}), (\ref{1.13}), (\ref{4.4}) and the Lebesgue dominated convergence theorem, taking the limit as $n\rightarrow +\infty$ in (\ref{4.9}) yields that
\begin{eqnarray*}
\int_{{\mathbb R}^d}u(t,X_t(x))\varphi(x)dx=0, \quad {\mathbb P}-a.s.,
\end{eqnarray*}
which implies $u(t,X_t(x))=0$ for almost everywhere $x\in {\mathbb R}^d$ and almost all $\omega\in\Omega$. Since $X_t(x)$ is a stochastic flow of homeomorphisms associated with (\ref{1.1}) with $\sigma=I_{d\times d}$ and $s=0$, we have $u(t,x)=0$ for almost everywhere $x\in {\mathbb R}^d$ and almost all $\omega\in\Omega$.

\medskip
(ii) Now let us show that (\ref{1.16}) holds. Since $X_t^{-1}(x)$ is differentiable in $x$, we have the following chain rule
\begin{eqnarray}\label{4.10}
\nabla_xu_0(X^{-1}_t(x))=\nabla u_0(X^{-1}_t(x))
\nabla X^{-1}_t(x).
\end{eqnarray}
Let $R>0$ be a real number. Recalling (\ref{1.8}) and (\ref{1.9}), one can see that for every $p\in [1,+\infty)$
\begin{eqnarray}\label{4.11}
\sup_{0\leq t\leq T}{\mathbb E}\sup_{x\in B_R}\|\nabla X^{-1}_t(x)\|^p\leq C(d,T,\Theta,p,\|b\|_{L^q(0,T;{\mathcal C}^{\frac{2}{q}-1}\cap{\mathcal C}^\alpha({\mathbb R}^d;{\mathbb R}^d))})<+\infty.
\end{eqnarray}
Combining (\ref{4.10}) and (\ref{4.11}), we get for $r<+\infty$
\begin{eqnarray}\label{4.12}
&&\sup_{0\leq t\leq T}{\mathbb E}\int_{B_R}|\nabla_xu_0(X^{-1}_t(x))|^rdx\nonumber\\&\leq&\sup_{0\leq t\leq T}{\mathbb E}\Big[\int_{{\mathbb R}^d}|\nabla u_0(X^{-1}_t(x))|^rdx\sup_{x\in B_R}\|\nabla X^{-1}_t(x)\|^r\Big]\nonumber\\&=&\int_{{\mathbb R}^d}|\nabla u_0(x)|^rdx\sup_{0\leq t\leq T}{\mathbb E}\sup_{x\in B_R}\|\nabla X^{-1}_t(x)\|^r<+\infty.
\end{eqnarray}
Moreover, if $r=+\infty$, one has
\begin{eqnarray}\label{4.13}
\sup_{0\leq t\leq T}{\mathbb E}\sup_{x\in B_R}\|\nabla_xu_0(X^{-1}_t(x))\|^p\leq
\|\nabla u_0\|^p_{L^\infty({\mathbb R}^d)}\sup_{0\leq t\leq T}{\mathbb E}\sup_{x\in B_R}\|\nabla X^{-1}_t(x)\|^p<+\infty,
\end{eqnarray}
for every $p\in [2,+\infty)$. Then the required assertion follows from (\ref{4.12}) and (\ref{4.13}). $\Box$

\section*{Acknowledgements.} Jinlong Wei was partly supported by the National Natural Science Foundation of China grant 12171247. Junhao Hu was partially supported by the National Natural Science Foundation of China grant 62373383 and the Fundamental Research Funds for the Central Universities of South-Central Minzu University grants KTZ20051 and CZT20020.

\section*{Conflicts of Interest}
The authors declare that they have no competing interests.

\section*{Data Availability Statements}
Data sharing is not applicable to this article as no new data were created or analyzed in this study.

\end{document}